\documentclass[a4paper]{amsart}
\usepackage{amsmath,amsfonts,amssymb,latexsym,epic,eepic}
\usepackage[overload]{empheq}
\usepackage[dvips]{graphicx,graphics,epsfig,color}
\usepackage{ifthen}
\usepackage[top=2cm, bottom=2cm, left=2cm, right=2cm]{geometry}
\usepackage{xcolor}

\usepackage{amsthm}
\usepackage{mathrsfs}
\usepackage{multirow}
\usepackage{array}
\usepackage{enumerate}

\def\xR{{\mathbb R}}

\newcommand{\N}{\mathbb{N}}
\newcommand{\R}{\mathbb{R}}
\newcommand{\C}{\mathbb{R}}

\DeclareMathOperator{\arcsinh}{arcsinh}

\theoremstyle{plain}
\newtheorem{thm}{Theorem}[section]

\newtheorem*{conj}{Conjecture}

\theoremstyle{remark}
\newtheorem{rem}[thm]{Remark}

\begin{document}

\title[]{A numerical study of stability for solitary waves of a quasi-linear Schrödinger equation}

\author{Meriem Bahhi}
\address{Laboratoire J.A. Dieudonné, Université Côte d'Azur, CNRS UMR  7351,\\
Parc Valrose, 28, avenue Valrose, 06108 Nice Cedex 2, France\\
E-mail 
Meriem.BAHHI@univ-cotedazur.fr}

\author{Jonas Lampart}
\address{Laboratoire Interdisciplinaire Carnot de Bourgogne, CNRS UMR 6303\\
                Universit\'e Bourgogne Europe, 9 avenue Alain Savary, 21078 Dijon
                Cedex, France\\
    E-mail jonas.lampart@u-bourgogne.fr}

\author{Christian Klein}
\address{Institut de Math\'ematiques de Bourgogne,CNRS UMR 5584\\
Institut Universitaire de France\\
                Universit\'e Bourgogne Europe, 9 avenue Alain Savary, 21078 Dijon
                Cedex, France\\
    E-mail Christian.Klein@u-bourgogne.fr}

\author{Simona Rota Nodari}
\address{Laboratoire J.A. Dieudonné, Université Côte d'Azur, CNRS UMR  7351\\
Parc Valrose, 28, avenue Valrose, 06108 Nice Cedex 2, France\\
Institut Universitaire de France\\
E-mail simona.rotanodari@univ-cotedazur.fr}
\date{\today}

\thanks{This work was funded by the EIPHI Graduate School (ANR-17-EURE-0002) and the Bourgogne-Franche-Comté Region through the project EEQuaR}

\begin{abstract}
We discuss the (in)stability of solitary waves for a quasi-linear 
Schrödinger equation. The equation contains a quasi-linear term, 
responsible for a saturation effect, as well as a power nonlinearity.
For different exponents of the nonlinearity, we determine analytically the asymptotic behavior of the $L^2$-mass of the solution as a function of the frequency close to the critical frequencies, which leads to natural conjectures concerning their stability.
Depending on the exponent and the dimension, we expect all solitary waves to be stable, or the emergence of both a stable and an unstable branch of solutions.
We investigate our conjectures numerically, and find compatible results both for the mass-energy relation and the dynamics.
We observe that perturbations of solitary waves on the unstable branch
may converge dynamically to the stable solution of a similar mass, or disperse. More general initial conditions show a similar behavior.
\end{abstract}


\maketitle

\section{Introduction}

This paper is devoted to the study of solutions to the nonlinear Schrödinger (NLS) type equation 
\begin{equation}
\label{SSNLt}
i\partial_{t}\phi = -\nabla \cdot\left(\frac{\nabla\phi}{1-|\phi|^{2\alpha}} \right) +\alpha |\phi|^{2\alpha-2}\frac{|\nabla\phi|^2}{(1-|\phi|^{2\alpha})^2}
\phi - |\phi|^{2\alpha}\phi ,
\end{equation}
in $\R^d$, $d\ge 1$. From a physical point of view, this equation can be seen as a specific non-relativistic limit, proper to nuclear physics, and describes the dynamical evolution of a particle inside the atomic nucleus. Here, $\phi:\R_+\times \R^d\to \C$ is a function that describes the quantum state of a nucleon (a proton or a neutron) and $\alpha\in \N^*$.

We are interested in the (in)stability properties of solitary wave solutions, i.e., solutions of the form $\phi(t,x)=e^{i\omega t}\varphi(x)$ with $\varphi$ a real positive $H^1$ solution to the stationary equation
\begin{equation}
\label{SSNL}
-\omega\varphi = -\nabla \cdot\left(\frac{\nabla\varphi}{1-|\varphi|^{2\alpha}} \right) +\alpha |\varphi|^{2\alpha-2}\frac{|\nabla\varphi|^2}{(1-|\varphi|^{2\alpha})^2}
\varphi - |\varphi|^{2\alpha}\varphi .
\end{equation}

This stationary equation has a natural variational interpretation and solitary wave solutions can be seen as critical points of the energy functional 
\begin{equation}
	\mathcal E(\varphi)=\frac12 \int_{\xR^d} \frac{|\nabla\varphi|^2}{1-|\varphi|^{2\alpha}}  -  \frac{1}{2(\alpha+1)} \int_{\xR^d}|\varphi|^{2\alpha+2} 
	\label{ener}
\end{equation}
subject to the mass constraint $\int_{\R^d}|\varphi|^2=M>0$. In particular, we will call (normalized) ground state the solution of the minimization problem
\begin{equation*}
  I(M)=\inf_{\varphi\in X,\ \int_{\R^d}|\varphi|^2=M} \mathcal E(\varphi),
\end{equation*}
where the functional space $X$ is given by
\begin{equation*}
  \label{def:functional_space}
  X=\left\{\varphi\in L^2(\R^d),\int_{\xR^d} \frac{|\nabla\varphi|^2}{(1-|\varphi|^{2\alpha})_+}<+\infty\right\}\subset H^1(\R^d)
\end{equation*}
and $f_+$ denotes the positive part of any function $f$. This is the appropriate functional framework to define ground state solutions for this particular energy (see~\cite{EstRot-13}): $\mathcal E$ is indeed not bounded from below in the set $\left\{\varphi \in H^1(\R^d), \int_{\R^d}|\varphi|^2=M\right\}$. Moreover, it has been shown in~\cite{EstRot-13} that if $\varphi\in X$, then $|\varphi|^2\le 1$ \emph{a.e.} in $\R^d$. As a consequence, in what follows, we consider solitary waves with $0<\varphi<1$.

In section~\ref{sect:existence} we discuss the existence and uniqueness of solutions  to~\eqref{SSNL} based on results from \cite{KRN} (for $d=1$) and \cite{BahhiPhD-25} (for $d\geq 2$). In particular, for given $\omega$ there exists a unique positive, radial decreasing solution $\varphi_\omega$ to~\eqref{SSNL} if and only if $\omega\in (0,\omega^*)$, with $\omega^*=(1+\alpha)^{-1}$.

Based on a result of Kenig, Ponce, Vega~\cite{kenig2004cauchy} we expect that the Cauchy problem for the  time-dependent equation~\eqref{SSNLt} is (locally) well posed for an open set of sufficiently regular initial data, for which additionally the leading quasi-linear term has adequate dispersive properties. In particular,  we expect that this holds in a neighborhood of the positive solitary wave solutions, which are radial, decreasing and thus non-trapping in the sense of~\cite{kenig2004cauchy}.
Note that $\mathcal E$ defined by~\eqref{ener} is the Hamiltonian energy of the  equation~\eqref{SSNLt}, and thus conserved by the dynamics.

We expect (normalized) ground states to be stable under these dynamics, while perturbations of solitary waves that do not minimize the energy at their mass may decay into stable solitary waves of a different frequency plus some radiation, which disperses.
An important tool to study the stability is the function $\omega\mapsto M(\omega)$ assigning to the Lagrange multiplier $\omega$ the mass of the corresponding solitary wave (this function is well defined by the results of Section~\ref{sect:existence}).
 Indeed, one should be able to prove (see \cite{Weinstein-85}) that, $\varphi_\omega$ is a strict local minimum of $\mathcal E$ at fixed $L^2$ norm, if $M'(\omega)>0$, whereas it is not a local minimum if $M'(\omega)<0$.

The Grillakis-Shatah-Strauss theory of stability \cite{GriShaStr-87,GriShaStr-90,Weinstein-85,BieGenRot-15,BieRot-19} indicates that the solitary wave solution $\phi_\omega$ (more precisely, the family $e^{i \theta}\varphi_\omega(x-c)$, $\theta,c\in \R$) of the time-dependent equation~\eqref{SSNLt} is expected to be \emph{orbitally stable} when $M'(\omega)>0$ and \emph{unstable} when $M'(\omega)<0$.

As a consequence, it is important to identify the regions where $M$ is increasing, which correspond to stable solutions/normalized ground states, and those where $M$ is decreasing, corresponding to unstable solutions.
For $d=1$, the explicit form the the solution $\varphi_\omega$ is known, and we can determine completely  the monotonicity of $M(\omega)$, which is given in Theorem~\ref{thm:asymptotics1}.
For $d\geq 2$ we can only determine analytically the asymptotic behavior for $\omega\to \omega^*$ and   $\omega\to 0$, which is given in Theorem~\ref{thm:asymptoticsD}. For $\omega\to \omega^*$, $M(\omega)$ is always divergent, while for $\omega\to 0$ the limit can be zero, finite, or infinite, depending on the values of $\alpha$ and $d$.
Based on this analysis and our numerical experiments, we expect the following stability to hold.

\begin{conj}Let $1\leq d\leq 6$.
\begin{itemize}
 \item For $\alpha\leq 2/d$, the solution $\varphi_\omega$ minimizes $\mathcal{E}$ at mass $M=M(\omega)$ for all $\omega\in (0,\omega^*)$ and is orbitally stable;
 \item For $\alpha>2/d$ there exists an $\omega_\mathrm{c}\in(0,\omega^*)$ so that for $\omega\geq \omega_\mathrm{c}$,  $\varphi_\omega$ is a (normalized) ground state and  orbitally stable, while for $\omega<\omega_\mathrm{c}$, $\varphi_\omega$ is not a minimizer of $\mathcal{E}$ at fixed $M=M(\omega)$ and unstable.
\end{itemize}
\end{conj}

In Section~\ref{sect:N1D} (for $d=1$), and Section~\ref{sect:NdD} (mainly for $d= 3$) we verify this conjecture numerically.
We evaluate the mass-energy relation, and find that for $\alpha>2/d$ it has a cusp at $M(\omega_\mathrm{c})$, where the lower branch (corresponding to frequencies $\omega\geq \omega_\mathrm{c}$), and the upper branch ($M(\omega),E(\omega))$, $\omega\leq \omega_\mathrm{c}$ intersect.
Considering initial data of the form $\lambda \phi_\omega$ with $\lambda\approx 1$, we confirm numerically the stability for $\omega\geq \omega_\mathrm{c}$.
For $\omega<\omega_\mathrm{c}$ we also confirm the instability, and also observe a dichotomy in the long-time asymptotics. For $\lambda>1$, the final state resembles a different solitary wave $\phi_{\omega'}$ plus some radiation, while for $\lambda<1$ the solution seems to disperse completely.
We also consider Gaussian initial data in $d=3$ dimensions and find a similar dichotomy depending on the width of the initial data.
These observations suggest strengthening the conjecture to a type of asymptotic stability and soliton resolution:
For \textit{any} initial condition, the solution is asymptotic to a solitary wave (or zero), with an error that converges to zero in an appropriate norm (such as $L^\infty$).

For any $d\geq 3$, all natural numbers $\alpha$ are supercritical and one may expect that solitary waves with small $\omega$ are always unstable. However, for $d\geq 7$ we do not have much theoretical or numerical evidence for this, so we do not make a conjecture at this point (see Remark~\ref{rem_dimension7} below).


\section{Solitary waves: existence and qualitative properties}\label{sect:existence}

In this section,
we collect the known results on the existence of solutions to equation~\eqref{SSNL} and describe their qualitative properties.

To our knowledge, equation~\eqref{SSNL} was mathematically studied for the first time in~\cite{EstRot-12} for the cubic nonlinearity ($\alpha=1$) and in dimension $d=3$, see also~\cite{TreRot-13}. In particular, the authors prove that if $\omega \ge \frac{1}{2}$ the equation~\eqref{SSNL} does not admit non-trivial solutions such that $\varphi(x)\to 0$ as $|x|\to +\infty$. If $\omega < \frac{1}{2}$, they proved the existence of solutions which change sign a fixed number of times and converge exponentially to $0$ when $|x|\to +\infty$. In~\cite{LewRot-15}, Lewin and Rota Nodari tackled the question of uniqueness and non-degeneracy of solutions to~\eqref{SSNL} for $\alpha=1$ and in any dimension $d\ge 2$. They proved that~\eqref{SSNL} has no non-trivial  solutions $0<\varphi<1$ in $L^2(\R^d)$ when $\omega \ge \frac{1}{2}$. If $\omega < \frac{1}{2}$, equation~\eqref{SSNL} has a unique positive solution $0<\varphi<1$ which tends to $0$ at infinity, modulo translations. Moreover, the positive solution is radial, decreasing and non-degenerate. The proof of this result is based on the realization that equation~\eqref{SSNL} can be written in terms of $\varphi=\sin u$ as simpler semilinear Schrödinger equation of the form $-\Delta u=\sin(u)^3\cos(u)-\omega \sin(u)\cos(u)$.
For $d = 1$, the existence of solitary waves for all $\alpha \in \mathbb{N}^*$
has been established in~\cite[Theorem 1.1]{KRN} by Klein and Rota Nodari.

\begin{thm}[\cite{KRN}]\label{them:existencesol1d}
Let $\alpha \in \mathbb{N}^*$ and $d = 1$. Define $\omega^* = \frac{1}{\alpha + 1}$. Then:
\begin{itemize}
    \item If $\omega \geq \omega^*$, equation~\eqref{SSNL} has no non-trivial solution $0 \leq \varphi < 1$ such that $\lim_{x \to \pm\infty} \varphi(x) = 0$.
    \item If $0 < \omega <\omega^*$, equation~\eqref{SSNL} admits a unique (up to translations) solution $0 < \varphi < 1$ such that $\varphi(x) \to 0$ as $x \to \pm\infty$. This solution is explicitly given by:
    \begin{equation}\label{gsexplicit}
    \varphi_\omega(x) = \left(1+ \left(\frac{\omega^*}{\omega}-1\right)  \cosh^2\left(\alpha\sqrt{\omega}x\right) \right)^{-1/(2\alpha)}.
    \end{equation}
    Moreover, $\varphi_\omega$ is even, strictly decreasing and non-degenerate.
\end{itemize}
\end{thm}

\begin{rem}
  Note that in~\cite{EstRot-12, TreRot-13, LewRot-15,KRN}, equation~\eqref{SSNL} was given in the form
  \begin{equation}
\label{SSNLab}
-b\varphi = -\nabla\cdot \left(\frac{\nabla\varphi}{1-|\varphi|^{2\alpha}} \right) +\alpha |\varphi|^{2\alpha-2}\frac{|\nabla\varphi|^2}{(1-|\varphi|^{2\alpha})^2}
\varphi - a|\varphi|^{2\alpha}\varphi .
\end{equation} 
with $a,b>0$. A straightforward computation shows that if $\varphi$ is a solution to~\eqref{SSNLab} then $x\to \varphi(a^{-1/2}x)$ is a solution to~\eqref{SSNL} with $\omega=\frac{b}{a}$. Conversely, if $\varphi$ is a solution to~\eqref{SSNL} then $x\to \varphi(a^{1/2}x)$ is a solution to~\eqref{SSNLab} with $b= a\omega$. 
\end{rem}

In the higher-dimensional case $d \geq 2$, one can exploit the idea that
the quasilinear equation~\eqref{SSNL} can be written as a semilinear equation of the form $-\Delta u=f_{\omega}(u)$ thanks to a change of unknown of the form $\varphi=r(u)$ as in~\cite{LewRot-15}. The main difficulty is that now the change of unknown $s\to r(s)$ is not explicitly given but is defined as a solution of an ordinary differential equation. This problem has been studied in Bahhi's PhD thesis \cite{BahhiPhD-25} where the following result has been obtained.

\begin{thm}[\cite{BahhiPhD-25}]
\label{Thrm_existence_uniqueness_quasi}
Let $\alpha \in \mathbb{N}^*$ and $d \ge 2$. Define $\omega^* = \frac{1}{\alpha + 1}$. Then:
\begin{itemize}
    \item If $\omega \geq \omega^*$, equation \eqref{SSNL} admits no non-trivial solution $\varphi \in H^1(\mathbb{R}^d)$ with $0 \leq \varphi < 1$.
    \item If $0 < \omega < \omega^*$, then \eqref{SSNL} admits a unique (up to translations) positive solution $\varphi_\omega \in H^1(\mathbb{R}^d)$ such that $0 < \varphi < 1$ and $\varphi(x) \to 0$ as $|x| \to \infty$. Moreover, after a suitable translation $\varphi_\omega$ is spherically symmetric, strictly decreasing in $|x|$ and non-degenerate.
\end{itemize}
\end{thm}

From now on we will denote by $\varphi_\omega$, $\omega \in (0,\omega^*)$, the unique positive, radial solution $0<\varphi_\omega<1$, and its $L^2$-mass by
\begin{equation}
  \label{eq:def_mass}
  M(\omega)=\int_{\R^d} \varphi_\omega^2.
\end{equation}
%
%

\subsection{Bifurcation in the case $d=1$}

In dimension $d=1$, using the explicit form of $\varphi_\omega$, we obtain
\begin{align}
M(\omega) &= \frac{2}{\alpha \sqrt{\omega}} \,
\int_0^{+\infty} \left( 1 + \left( \frac{\omega^*}{\omega} - 1 \right) \cosh^2(y) \right)^{-1/\alpha} \, dy \nonumber\\
&= \frac{2}{\alpha \sqrt{\omega}}
\int_0^{+\infty} \left( \frac{\omega^*}{\omega} + \left( \frac{\omega^*}{\omega} - 1 \right) \sinh^2(y) \right)^{-1/\alpha} \, dy \nonumber\\
&= \frac{2}{\alpha} \, \frac{\omega^{1/\alpha - 1/2}}{(\omega^*)^{1/\alpha}}
\int_0^{+\infty} \left( 1 + z^2 \right)^{-1/\alpha}\left( \left( 1-\frac{\omega}{\omega^*} \right) +z^2 \right)^{-1/2} \, dz.
\end{align}

From this expression, we can deduce the asymptotic behavior of the mass as summarized below:

\begin{thm}\label{thm:asymptotics1}
Let $d = 1$ and $\alpha \geq 1$. Then the mass $M(\omega)$ satisfies:
\begin{enumerate}
    \item \textbf{As $\omega \to 0^+$}:
    \begin{equation}
     M(\omega) \underset{\omega \to 0}{\sim}\omega^{1/\alpha - 1/2}\frac{2}{\alpha \sqrt{\omega^*}}\int_0^{+\infty} \left( 1 + z^2 \right)^{-1/\alpha-1/2}\, dz,\notag
    \end{equation}
    and
    \begin{itemize}
        \item for $\alpha < 2$, $M(\omega)$ is strictly increasing with $ \lim_{\omega \to 0^+} M(\omega) = 0$,
        \item  for $\alpha = 2$, $M(\omega)$ is strictly increasing, convex with $\lim_{\omega \to 0^+} M(\omega)=\frac{\sqrt{3}\pi}{2}$,
        \item for $\alpha > 2$, $M(\omega)$ is convex, with $M'(\omega)$ vanishing exactly once in $(0, \omega^*)$ and $\lim_{\omega \to 0^+} M(\omega)=\infty$.
   \end{itemize}
    \item \textbf{As $\omega \to \omega^*$}:
    \begin{equation}
    M(\omega) \underset{\omega \to \omega^*}{\sim}-\frac{1}{\alpha \sqrt{\omega^*}}\ln(\omega^*-\omega).\notag
    \end{equation}
\end{enumerate}
\end{thm}

The proof of this theorem is given in the appendix.

\medskip

\begin{figure}[h!]
  \includegraphics[width=0.3\textwidth]{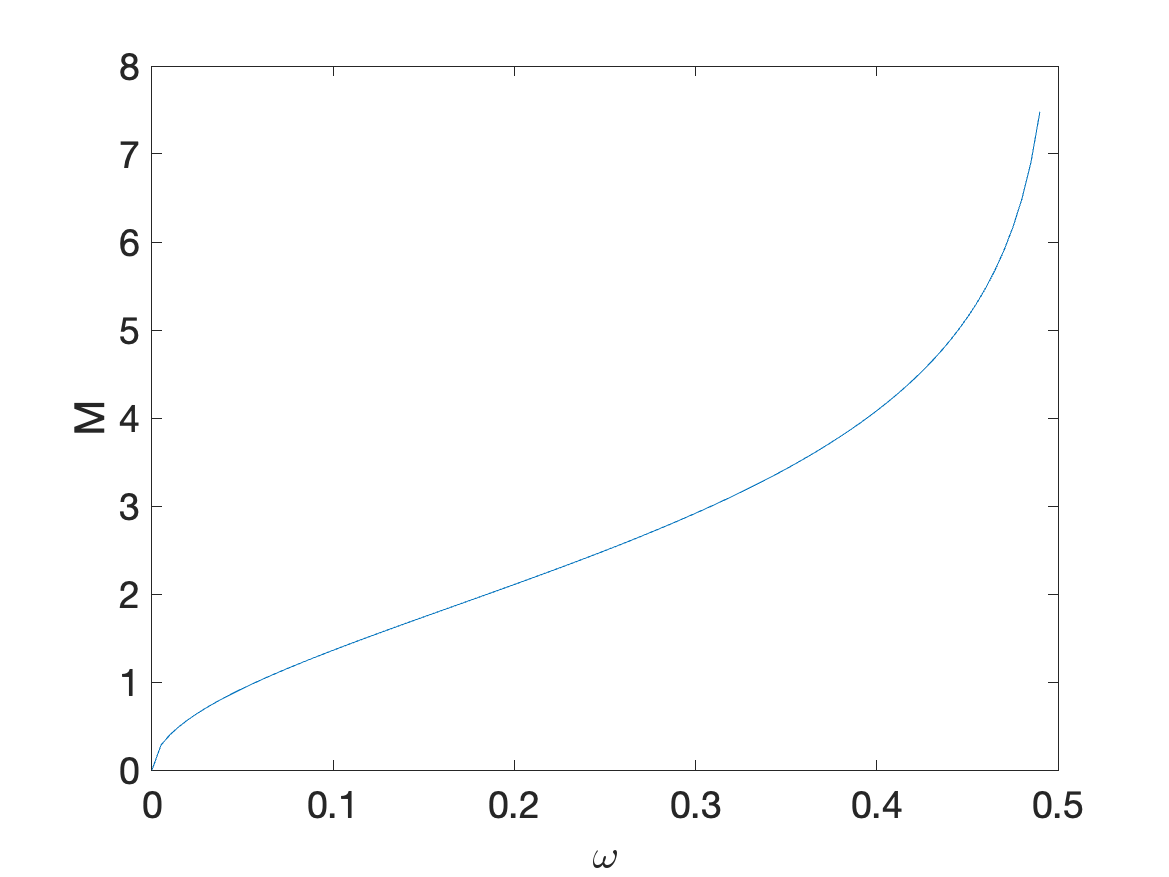}
  \includegraphics[width=0.3\textwidth]{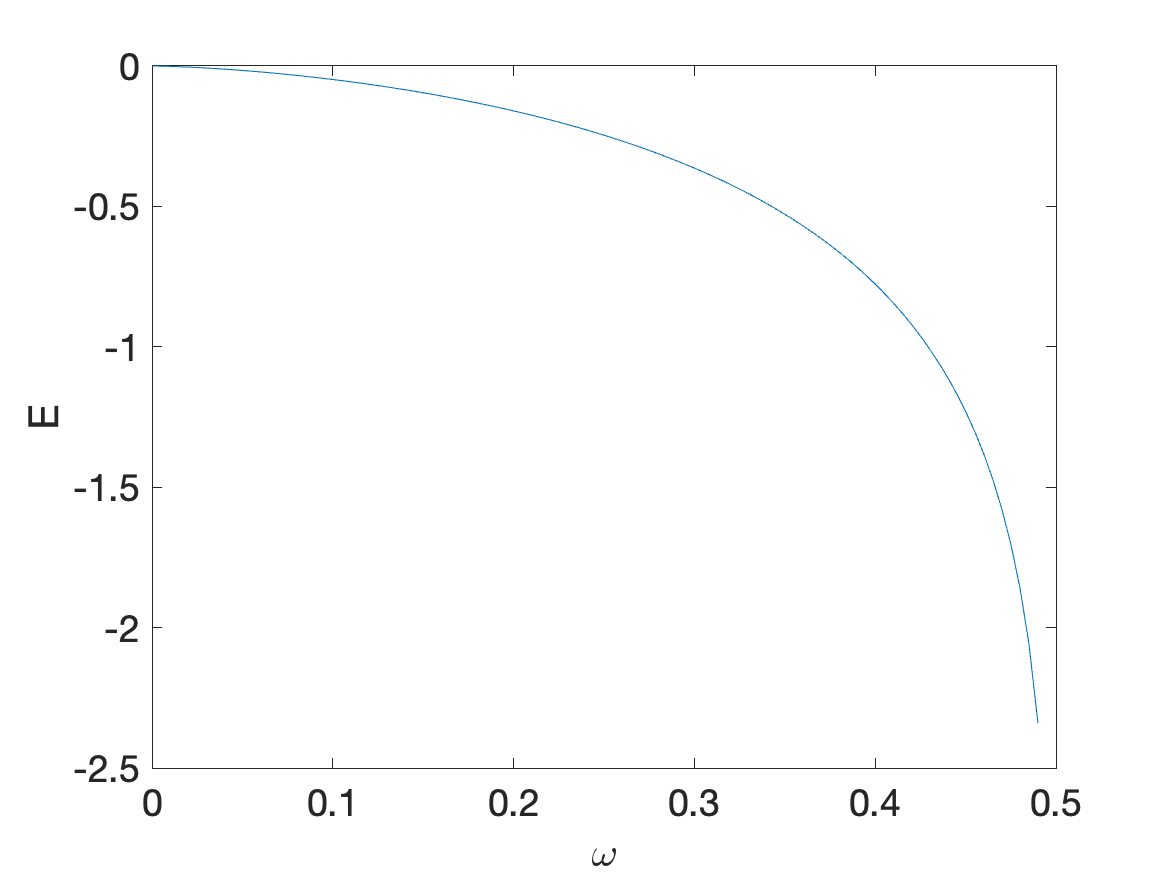}
  \includegraphics[width=0.3\textwidth]{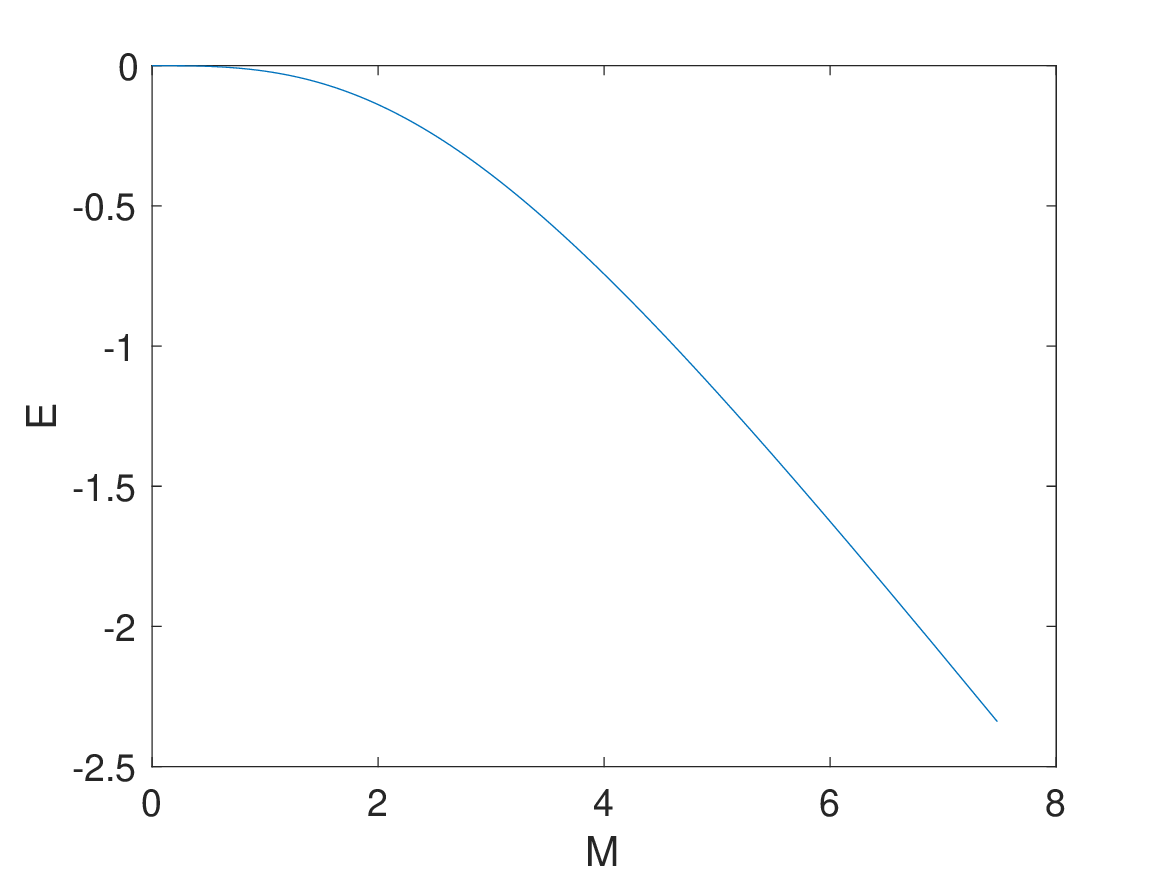}
 \caption{The mass (left) and the energy (middle) for the ground state 
 solutions  (\ref{gsexplicit}) for $\alpha=1$ in dependence of 
 $\omega$, and the energy in dependence of $M$ on the right.}
 \label{figMEalpha1}
\end{figure}

The case $\alpha=1$ is subcritical, the mass and energy are shown in 
Fig.~\ref{figMEalpha1}. Both the mass and the energy are monotone
and diverge at the critical frequency $\omega^{*}=1/(\alpha+1)$. 
Therefore, the energy is also a monotone function of $M$.

\begin{figure}[htb!]
  \includegraphics[width=0.3\textwidth]{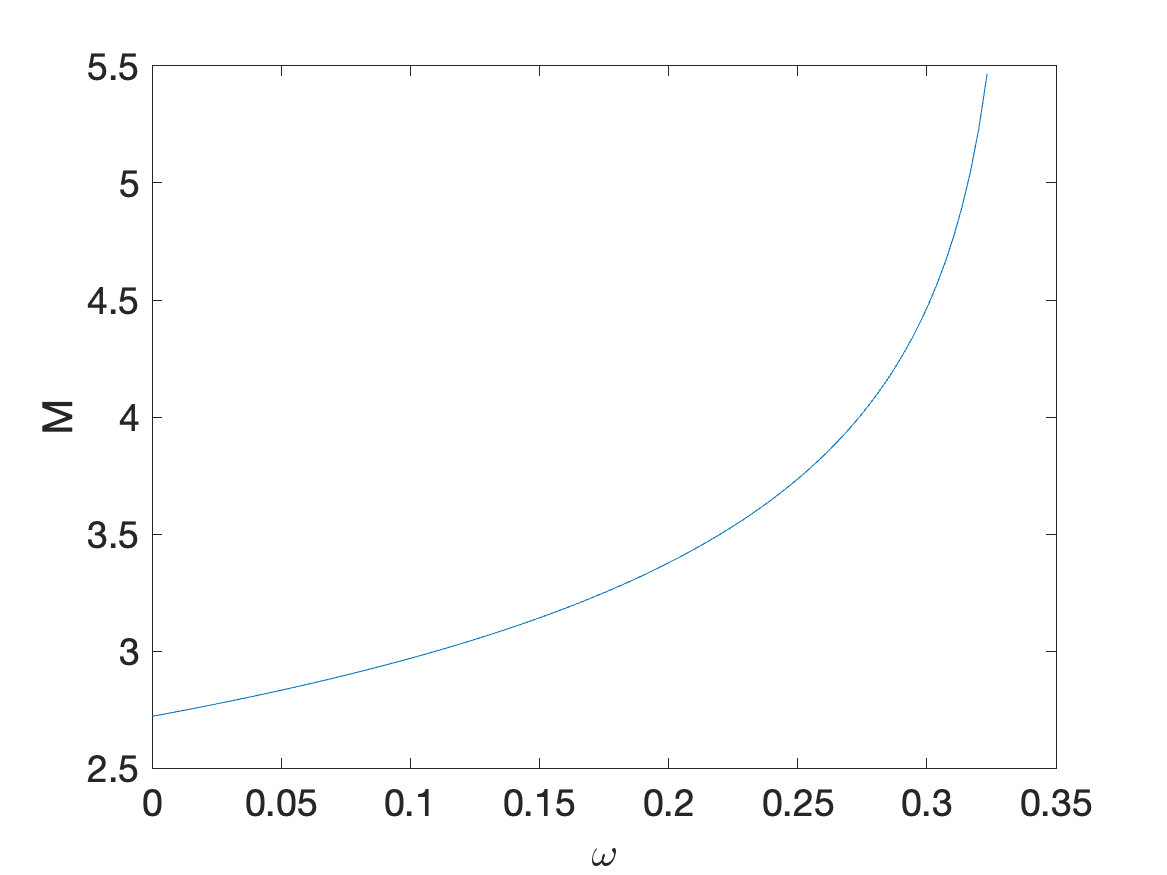}
  \includegraphics[width=0.3\textwidth]{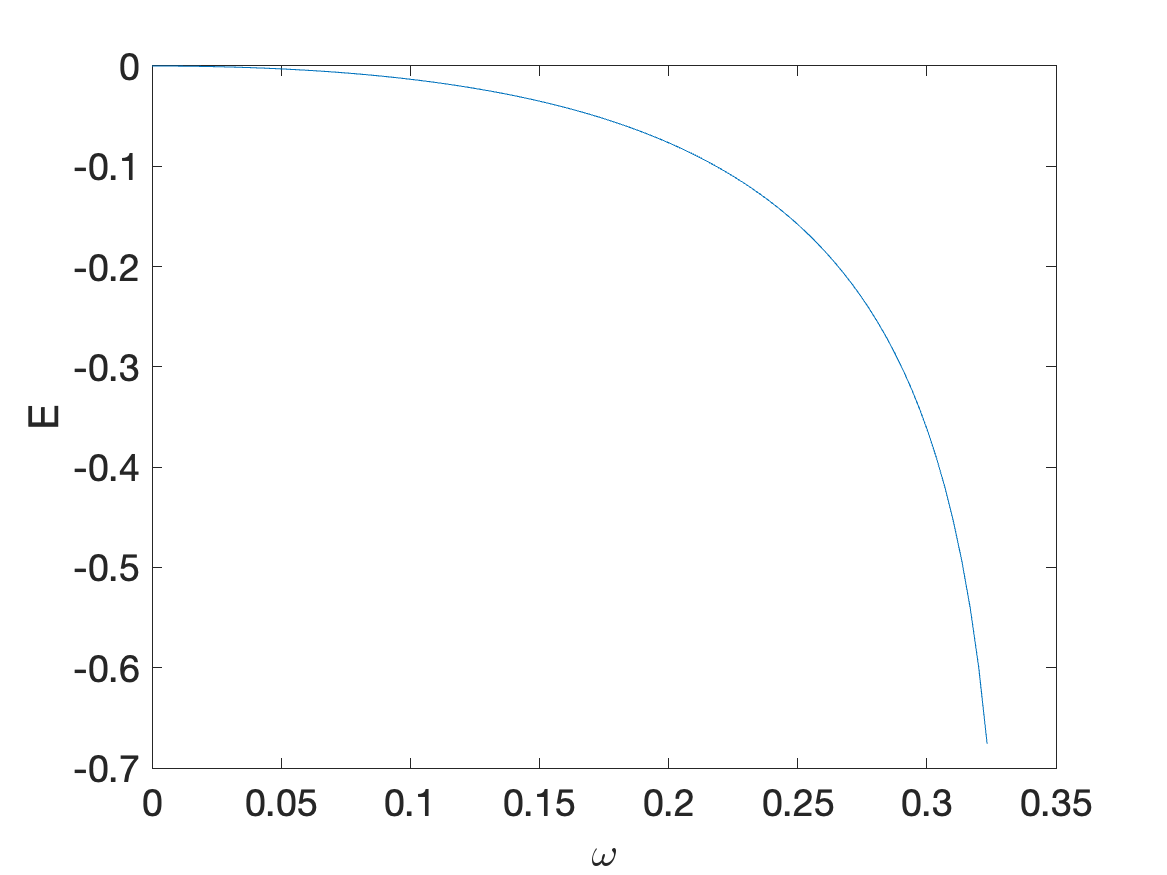}
  \includegraphics[width=0.3\textwidth]{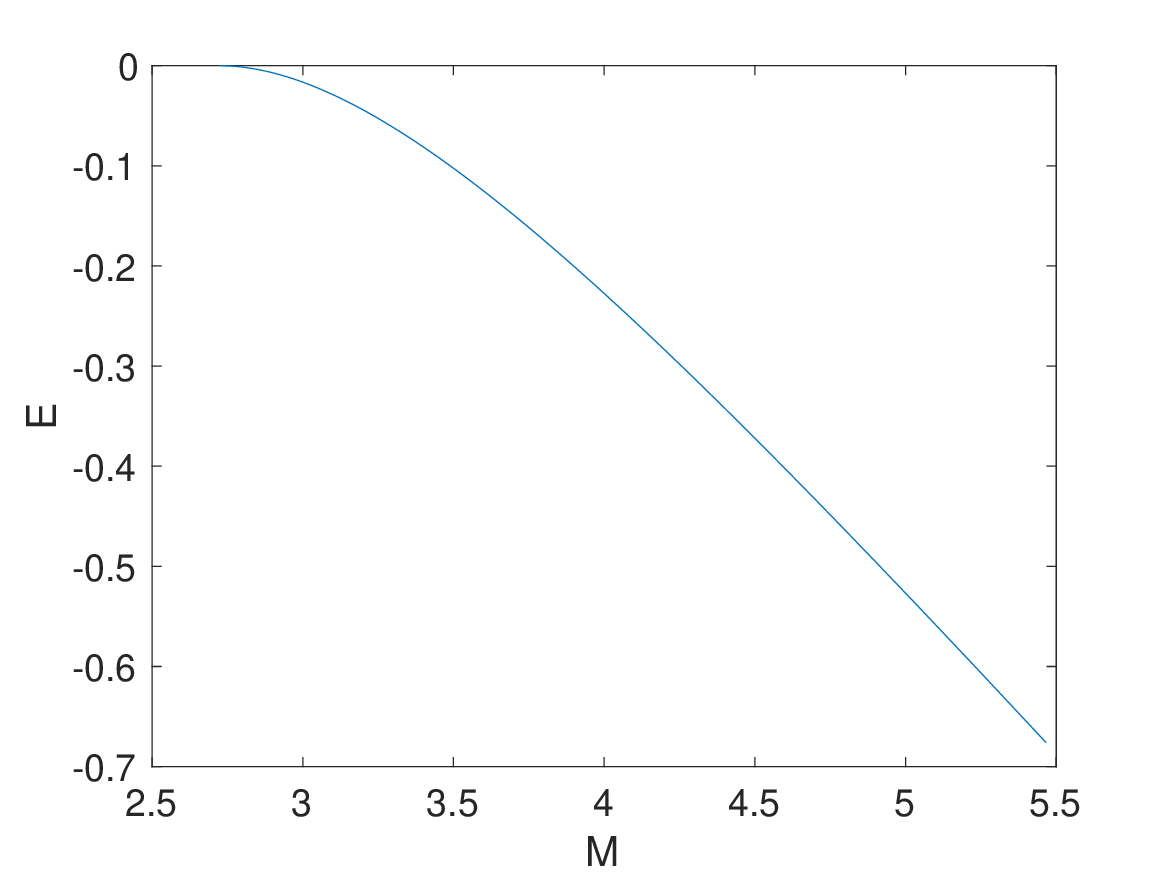}
 \caption{The mass (left) and the energy (middle) for the ground state 
 solutions  (\ref{gsexplicit}) for $\alpha=2$ in dependence of 
 $\omega$, and the energy in dependence of $M$ on the right.}
 \label{figMEalpha2}
\end{figure}

The situation is very similar in the critical case $\alpha=2$ as can 
be seen in Fig.~\ref{figMEalpha2}, but the mass no longer vanishes 
for $\omega\to0$. The energy is
non-positive in both cases.

\begin{figure}[htb!]
  \includegraphics[width=0.3\textwidth]{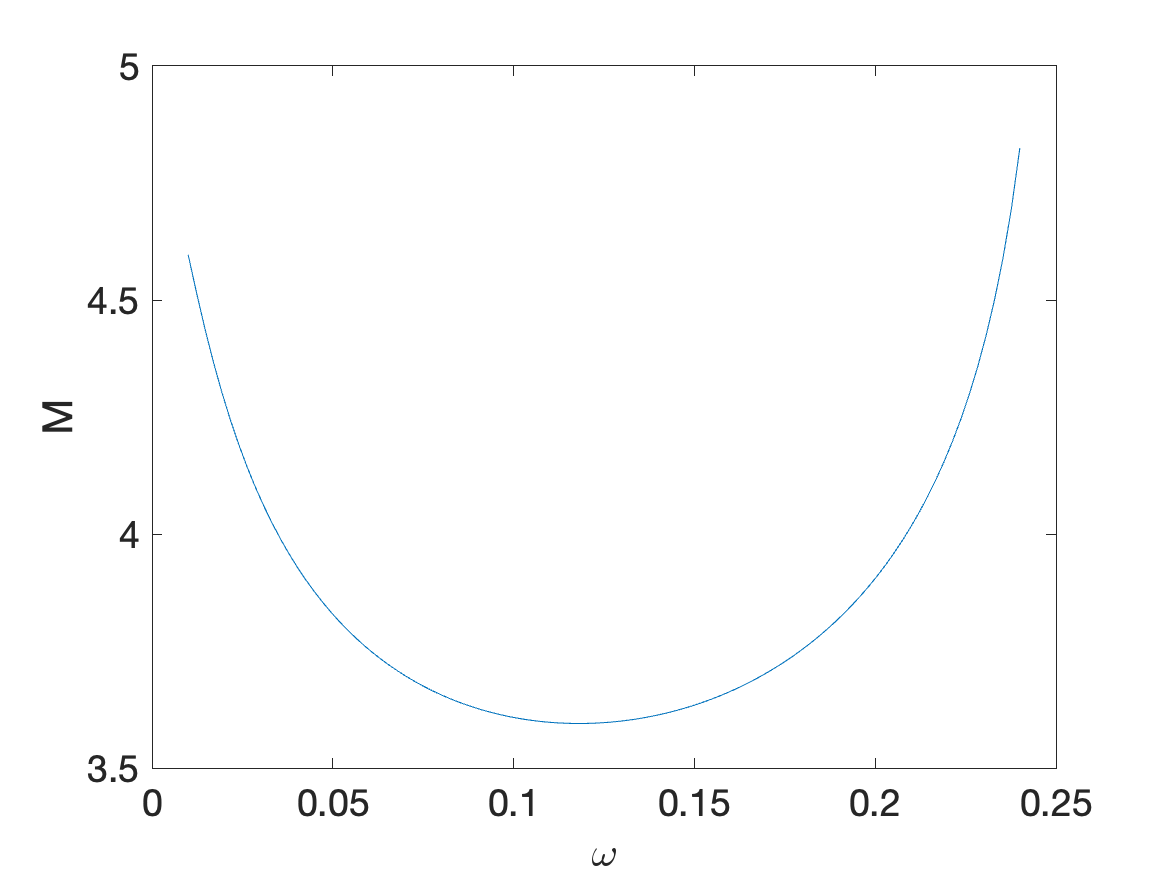}
  \includegraphics[width=0.3\textwidth]{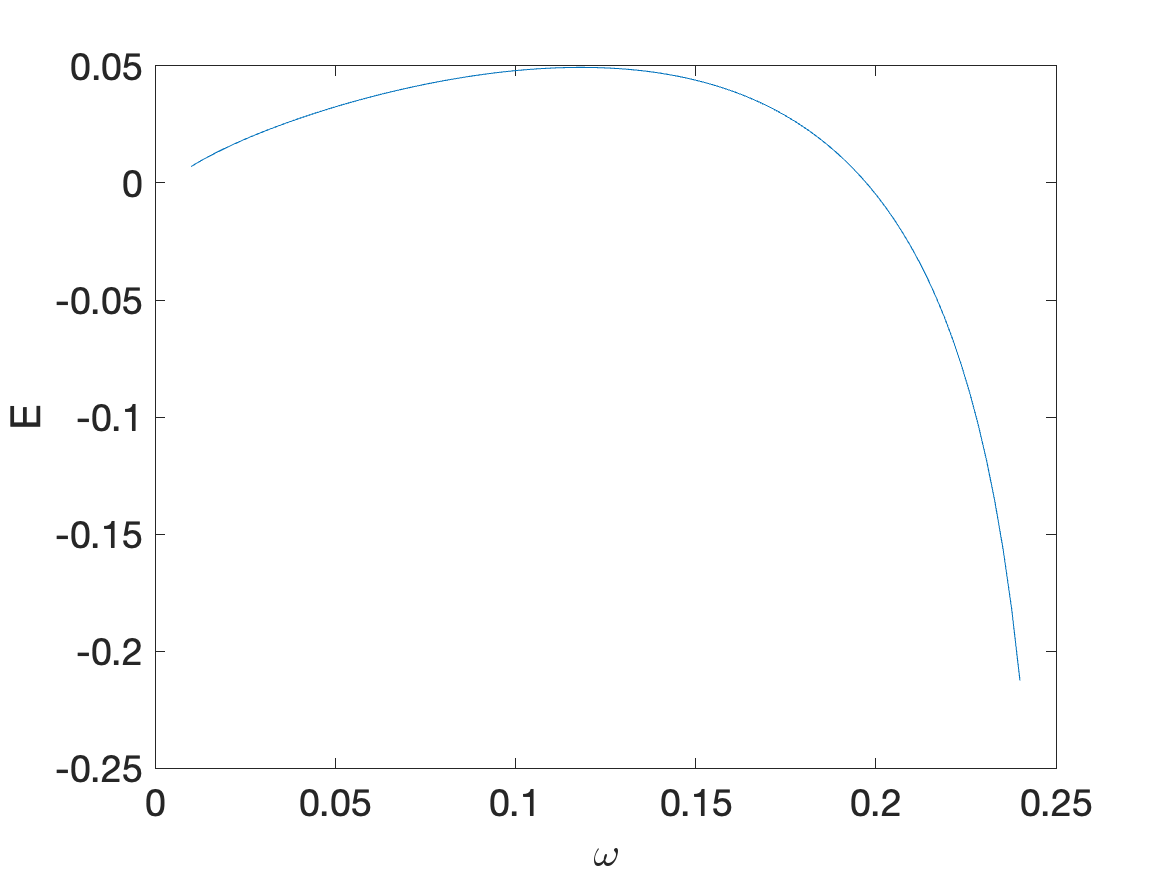}
  \includegraphics[width=0.3\textwidth]{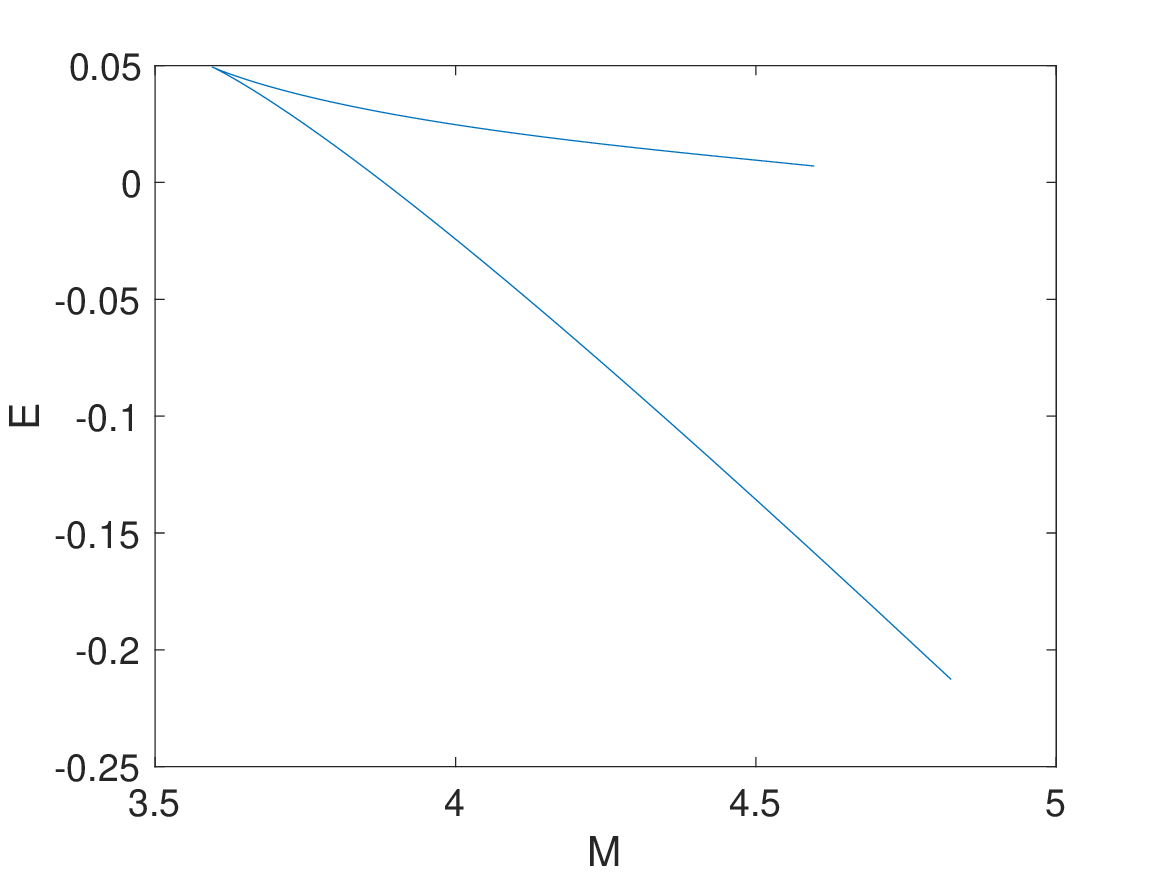}
 \caption{The mass (left) and the energy (middle) for the ground state 
 solutions  (\ref{gsexplicit}) for $\alpha=3$ in dependence of 
 $\omega$, and the energy in dependence of $M$ on the right.}
 \label{figMEalpha3}
\end{figure}

In the supercritical case $\alpha=3$ as shown in 
Fig.~\ref{figMEalpha3}, the situation is considerably different.  The 
mass has a global minimum and diverges also for $\omega\to0$, the 
energy has a maximum and is positive for some values of $\omega$, and 
the energy as a function of $M$ has a cusp. This indicates that there 
are two branches in this case, the one with the lower energy being 
the stable one. 

\subsection{Bifurcation in the higher-dimensional case}

In the higher-dimensional case, we cannot provide an explicit formula for the mass as a function of $\omega$. However, its asymptotic behavior can still be analyzed as $\omega \to 0^+$ and as $\omega\to \omega^*$. In particular, we have the following theorem.

\begin{thm}[\cite{BahhiPhD-25}]\label{thm:asymptoticsD}
Let $d \geq 2$ and $\alpha\in \N^*$. Then the mass $M(\omega)$ satisfies
\begin{enumerate}
  \item {\bf As $\omega\to 0^+$}:
  \begin{itemize}
    \item \textbf{Sub-critical case:} Suppose $d = 2$, or $d \geq 3$ and $\alpha < \frac{2}{d - 2}$. Then the mass $M(\omega)$ satisfies:
    \begin{equation}
\label{assymptotic_mass_phi}
  M(\omega) \underset{\omega \to 0^+}{\sim} \omega^{(\frac{1}{\alpha}-\frac{d}{2})} \|\psi_0\|^2_{L^2(\xR^d)},
\end{equation}
 with $\psi_0$, the unique positive, radial solution of
        \begin{equation} \label{limit_RSNL}
          -\Delta \psi + \psi - |\psi|^{2\alpha} \psi = 0, \notag
        \end{equation}
and
    \begin{itemize}
        \item if $\alpha < \frac{2}{d}$ ($L^2$-sub-critical case), then $\lim_{\omega\to 0^+} M(\omega)=0$,
        \item if $\alpha = \frac{2}{d}$ ($L^2$-critical case), then $\lim_{\omega\to 0^+} M(\omega)=\int_{\R^d} |\psi_0|^2$,
              \item if $\alpha > \frac{2}{d}$, ($L^2$-super-critical case),
        $
        \lim_{\omega\to 0^+} M(\omega)= \infty
        $.
      \end{itemize}
      \item \textbf{Critical case:} Suppose $d \geq 3$ and $\alpha = \frac{2}{d - 2}$. Then
     \begin{equation}
    \lim_{\omega \to 0^+} M(\omega) = \infty. \notag
    \end{equation}
    \item \textbf{Super-critical case:} Suppose $d \geq 3$ and $\alpha > \frac{2}{d - 2}$. Then
    , as $\omega \to 0^+$,
     \begin{equation}
    \lim_{\omega \to 0} M(\omega) =
    \begin{cases}
    \infty, & \text{if } d \in \{3, 4\}, \\
    \|\varphi_0\|_{L^2(\mathbb{R}^d)}, & \text{if } d \geq 5.
    \end{cases} \notag
    \end{equation}
    where $\varphi_0$ is the unique positive, radial solution of
    \begin{equation*}
    -\nabla \cdot \left( \frac{\nabla \varphi}{1 - |\varphi|^{2\alpha}} \right) + \alpha |\varphi|^{2\alpha - 2} \frac{|\nabla \varphi|^2}{(1 - |\varphi|^{2\alpha})^2} \varphi - |\varphi|^{2\alpha} \varphi = 0.
    \end{equation*}
  \end{itemize}

  \item {\bf As $\omega\to \omega^*$}:
  \begin{equation*}
M(\omega) \underset{\omega \to \omega^*}{\sim} \frac{\Gamma}{(\omega^* - \omega)^d}
\end{equation*}
with $\Gamma$ a constant that depends on $d,\alpha$ and $\omega^*$.
\end{enumerate}
\end{thm}

\begin{rem}
\label{rem_assymptotics}
The distinction between the asymptotic regimes comes from the fact that we have different limiting problems, which depend on the exponent $\alpha$. In the $H^1$-subcritical and $H^1$-supercritical settings, the limiting problems can be formulated explicitly, allowing us to derive precise asymptotics for the mass as $\omega \to 0$. In contrast, the critical case involves additional analytical difficulties and remains more delicate to treat.

When $\omega \to \omega^*$, we recover the asymptotic behavior of the mass described above by studying the behavior of $\varphi_\omega$ as $\omega \to \omega^*$. In particular, we proved in \cite{BahhiPhD-25} that $\varphi_\omega$ converges to $1$ locally, causing the denominator of the quasi-linear term in \eqref{SSNL} to blow up. This observation is consistent with the existence theorem, which states 
that no solutions exist for $\omega \geq \omega^*$, and can be 
interpreted as a `saturation phenomenon' for
nuclear matter. As $\omega\to \omega^*$, the solution $\varphi_\omega$ is more and more similar to a smoothed step function (see Fig.~\ref{NLSnuc3Dsol} and Fig.~\ref{NLSnuc3Dsolomlarge} below), with the transition behavior analytically described by a one-dimensional limiting equation (see \cite{BahhiPhD-25}, and \cite{Rota-HDR} for the case $\alpha=1$).
\end{rem}

\begin{rem}
  \label{rem_dimension7}
  In dimension $3$ and $4$, since $M(\omega)\to +\infty$ as $\omega \to 0$, it is natural to conjecture that solitary waves with small $\omega$ are unstable. For $d \ge 5$, the situation is more subtle, as $M(\omega)$ converges to a finite value when $\omega \to 0$. 
  However, following the approach of \cite{LewRot-20, GenRot-24}, one should be able to show that for $d \in \{5,6\}$, $M'(\omega) \to -\infty$ as $\omega \to 0$. Starting from dimension $7$, we expect $\lim_{\omega \to 0} M'(\omega)$ to be finite, but the sign of $M'(\omega)$ for $\omega$ close to $0$ remains unclear (see \cite{LewRot-20, GenRot-24} for a more detailed analysis of related problems). This motivates us to restrict our conjecture to dimensions $d \le 6$.
\end{rem}

The proof of this theorem is given in \cite{BahhiPhD-25} and follows the ideas of \cite{MorMur-14, LewRot-20, GenRot-24}.

\begin{figure}[htb!]
  \includegraphics[width=0.3\textwidth]{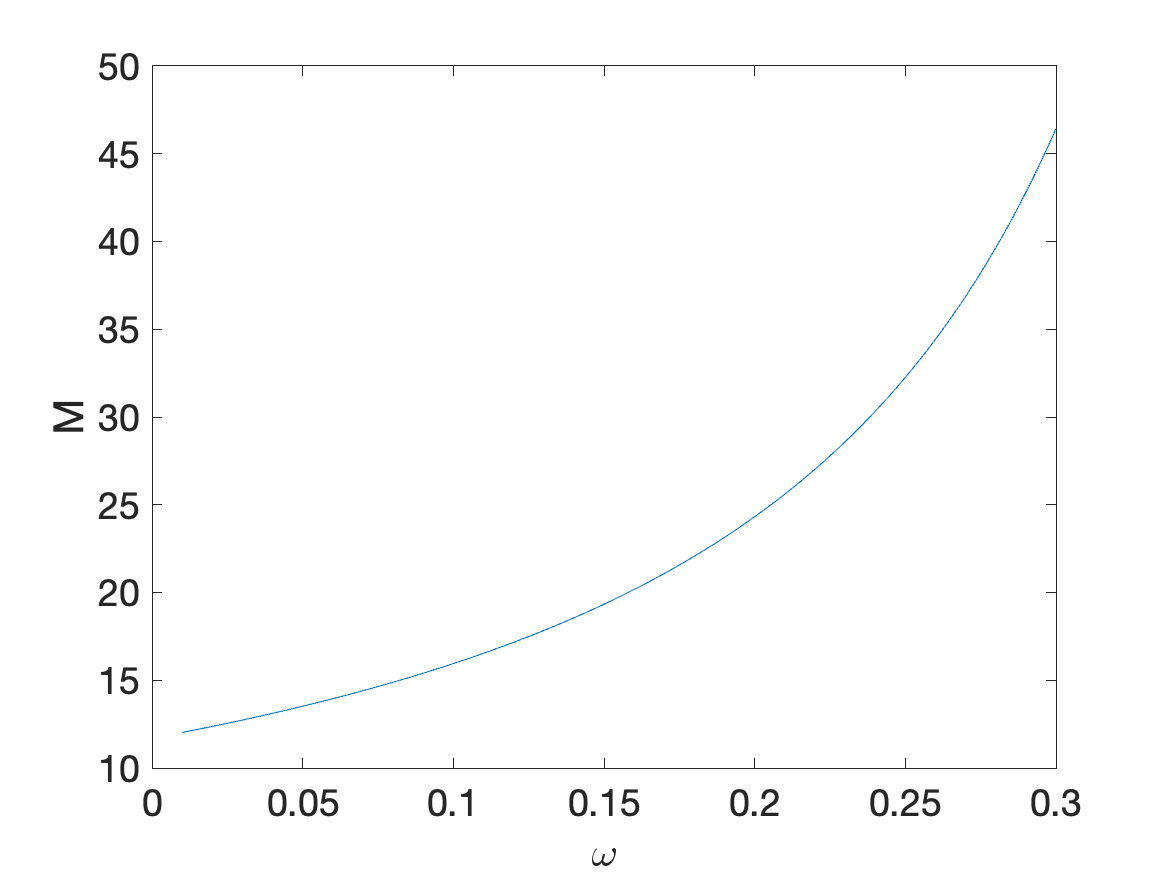}
  \includegraphics[width=0.3\textwidth]{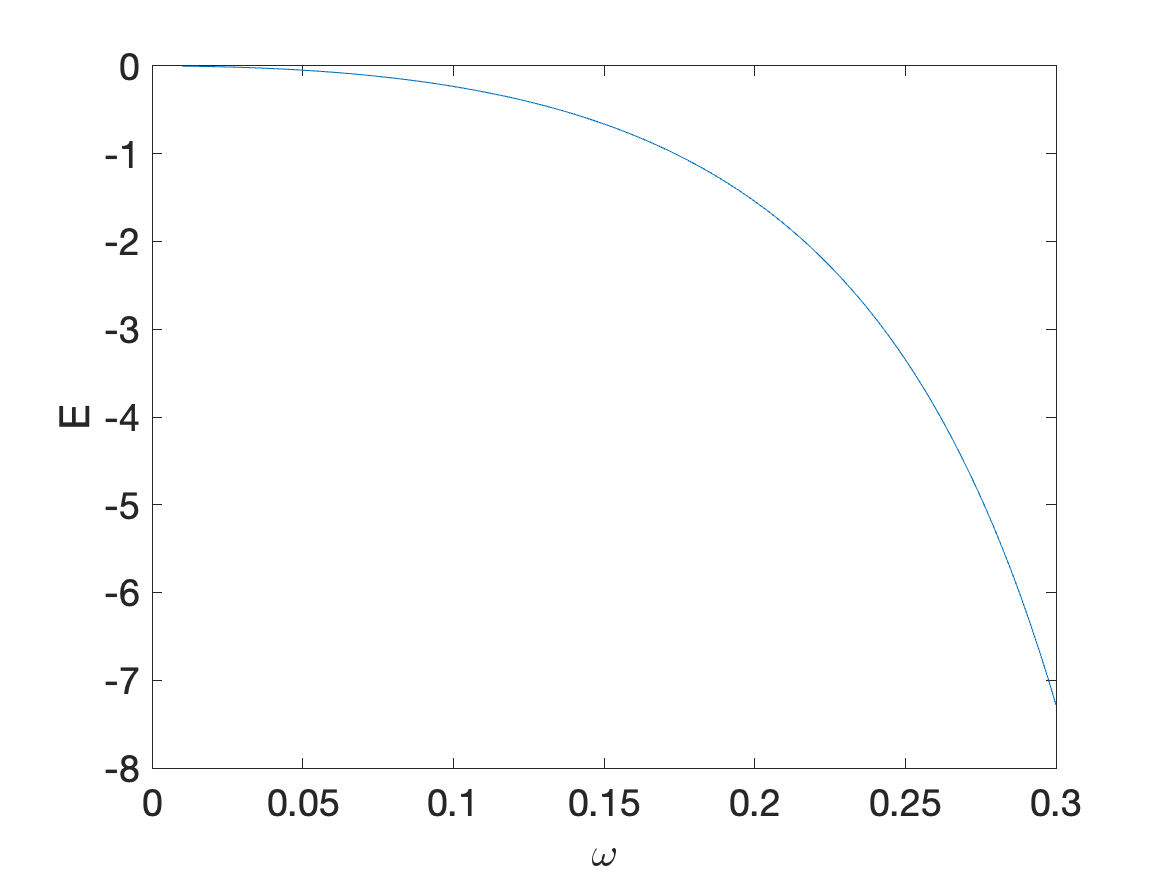}
  \includegraphics[width=0.3\textwidth]{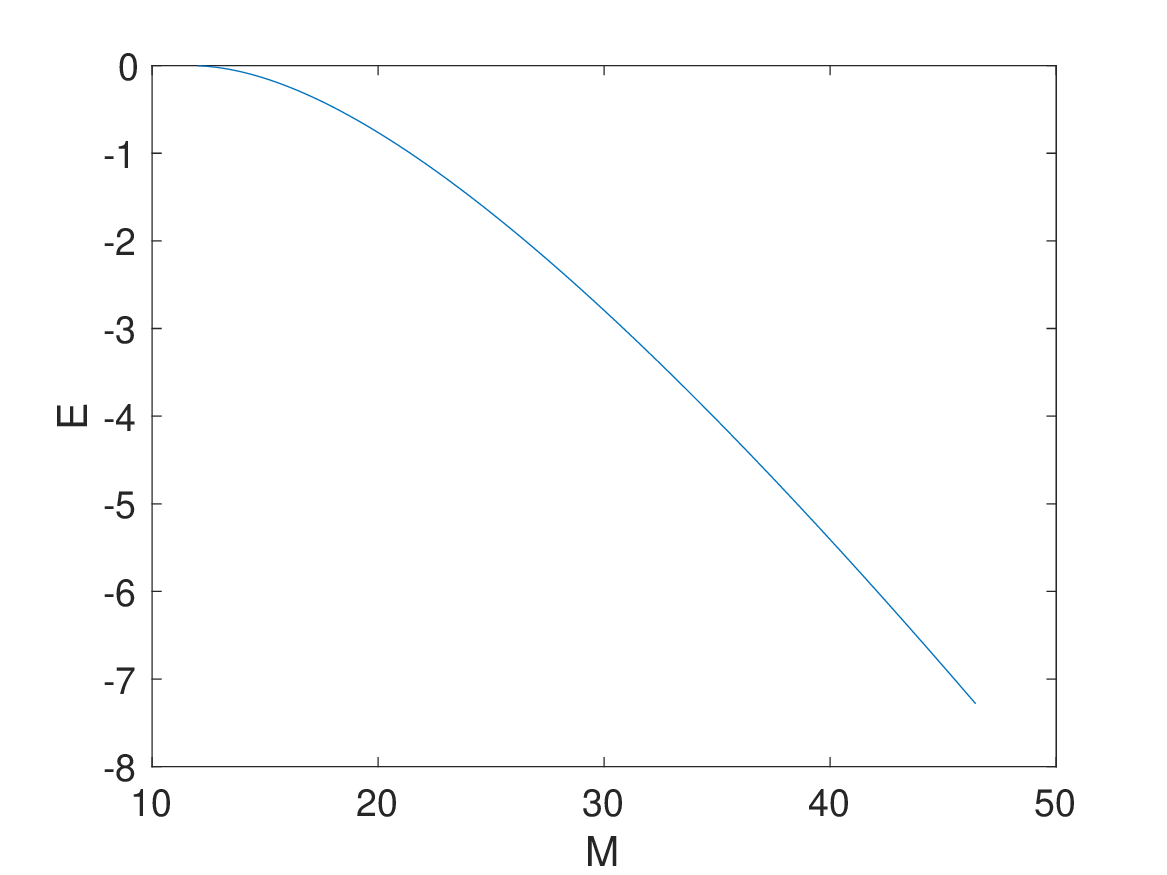}
 \caption{The mass (left) and the energy (middle) for the solitary 
 waves  in $d=2$ for $\alpha=1$ in dependence of
 $\omega$, and the energy in dependence of $M$ on the right.}
 \label{fig2dMEalpha1}
\end{figure}

The case $\alpha=1$ is $L^{2}$ critical for the standard NLS in $d=2$ and
$L^{2}$ supercritical in $d=3$. In the present case this implies that
the solitary waves have a mass tending to a finite value for 
$\omega\to0$ as can be seen in Fig.~\ref{fig2dMEalpha1}. The $\omega$ 
dependence of the energy is also monotone as can be seen in the same
figure. The situation is thus as  for $\alpha=2$ in $d=1$, see
Fig.~\ref{figMEalpha2}.

\begin{figure}[htb!]
  \includegraphics[width=0.3\textwidth]{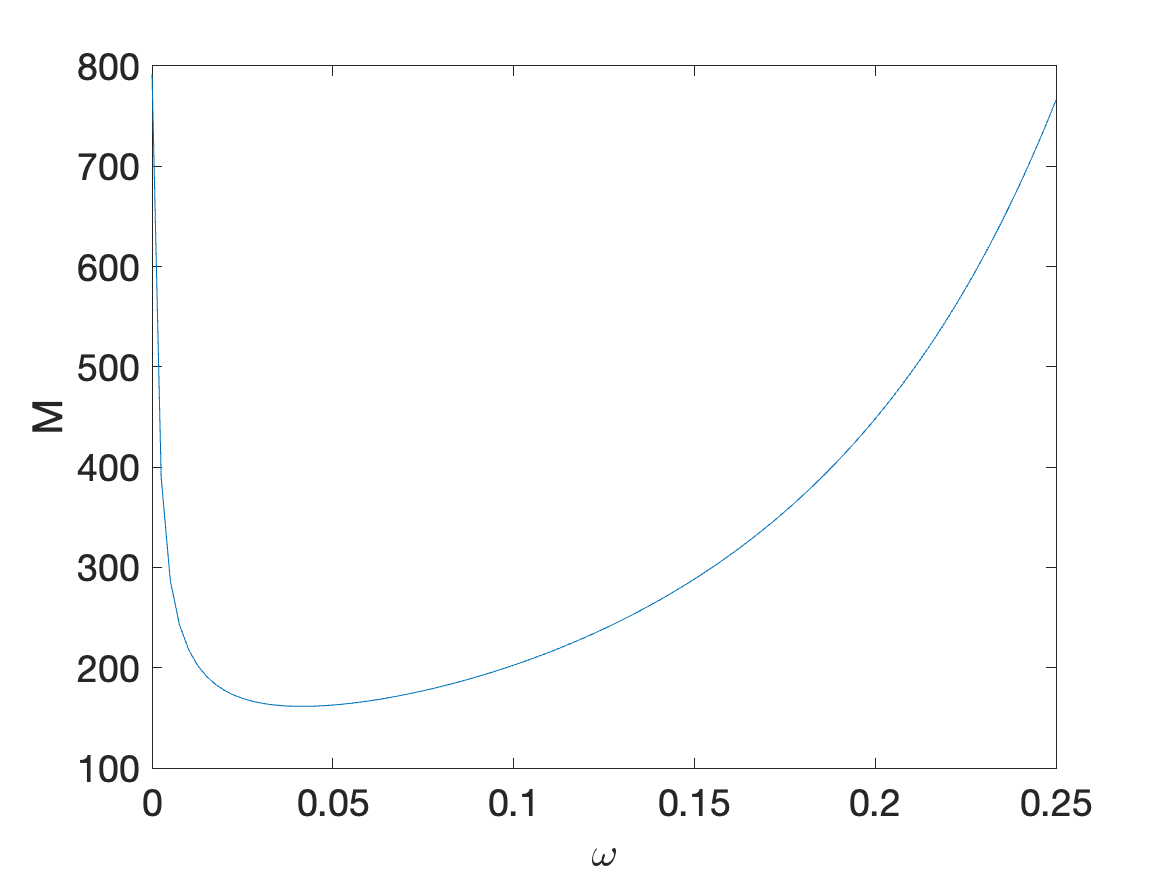}
  \includegraphics[width=0.3\textwidth]{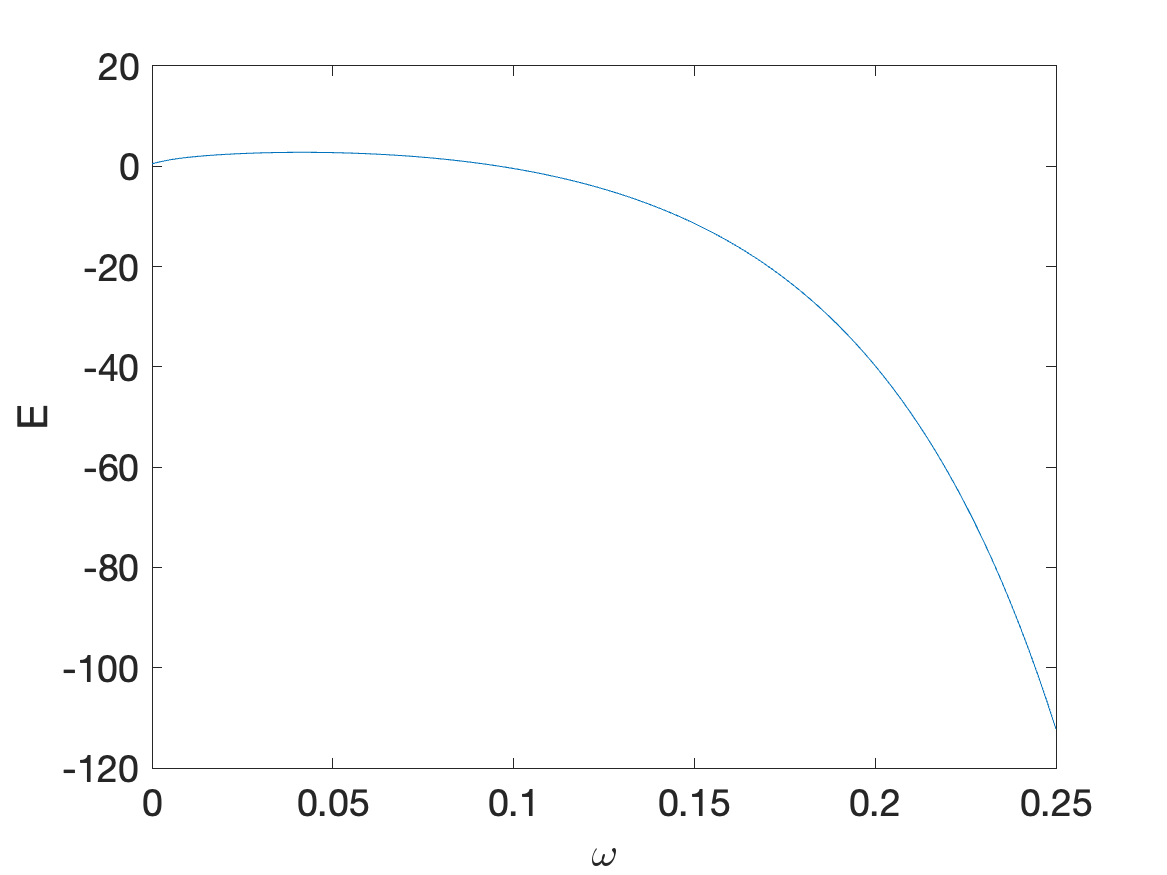}
  \includegraphics[width=0.3\textwidth]{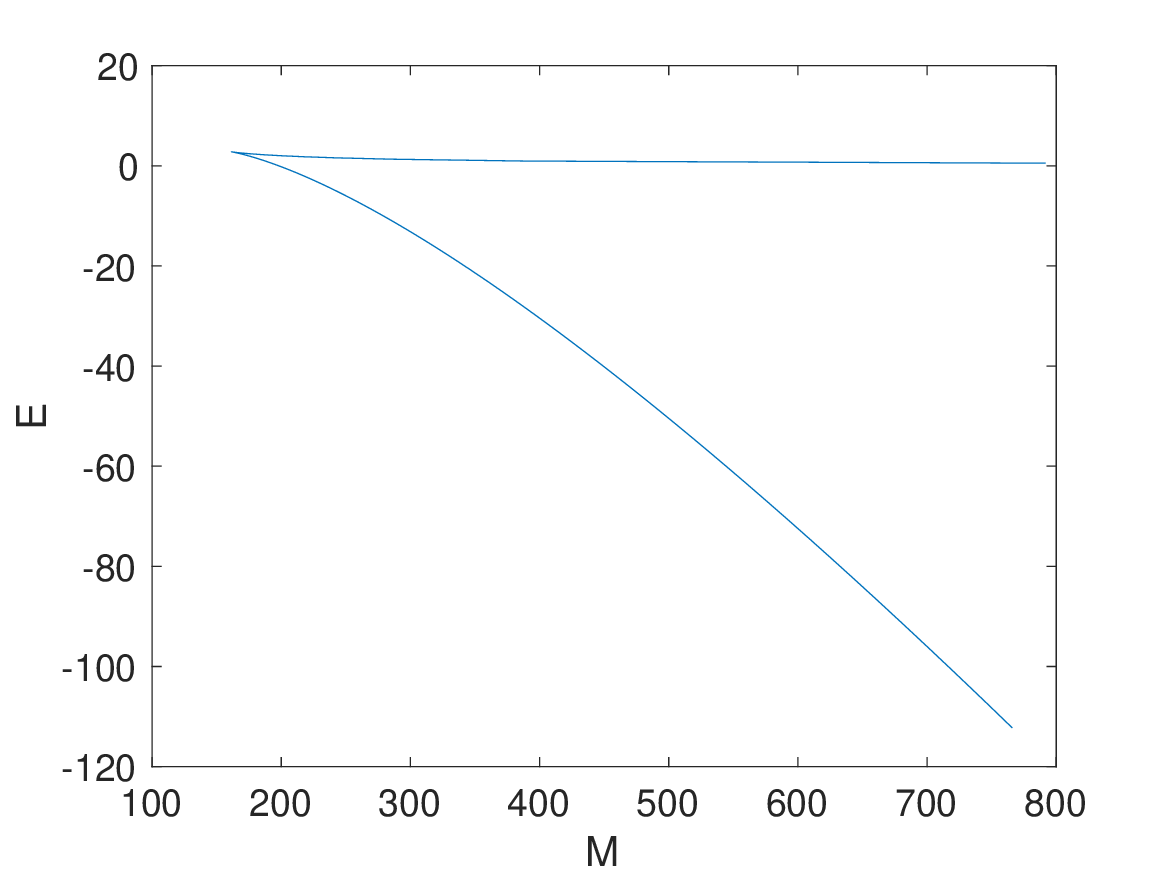}
 \caption{The mass (left) and the energy (middle) for the solitary 
 waves in $d=3$ for $\alpha=1$ in dependence of
 $\omega$, and the energy in dependence of $M$ on the right.}
 \label{fig3dMEalpha1}
\end{figure}

The same nonlinearity in $d=3$ leads to a different behavior similar to
the case $\alpha=3$ in $d=1$, see Fig. \ref{figMEalpha3} which can be
compared to  Fig.~\ref{fig3dMEalpha1}. It is remarkable that the 
unstable branch only appears here for very small values of $\omega$.

\begin{figure}[htb!]
  \includegraphics[width=0.3\textwidth]{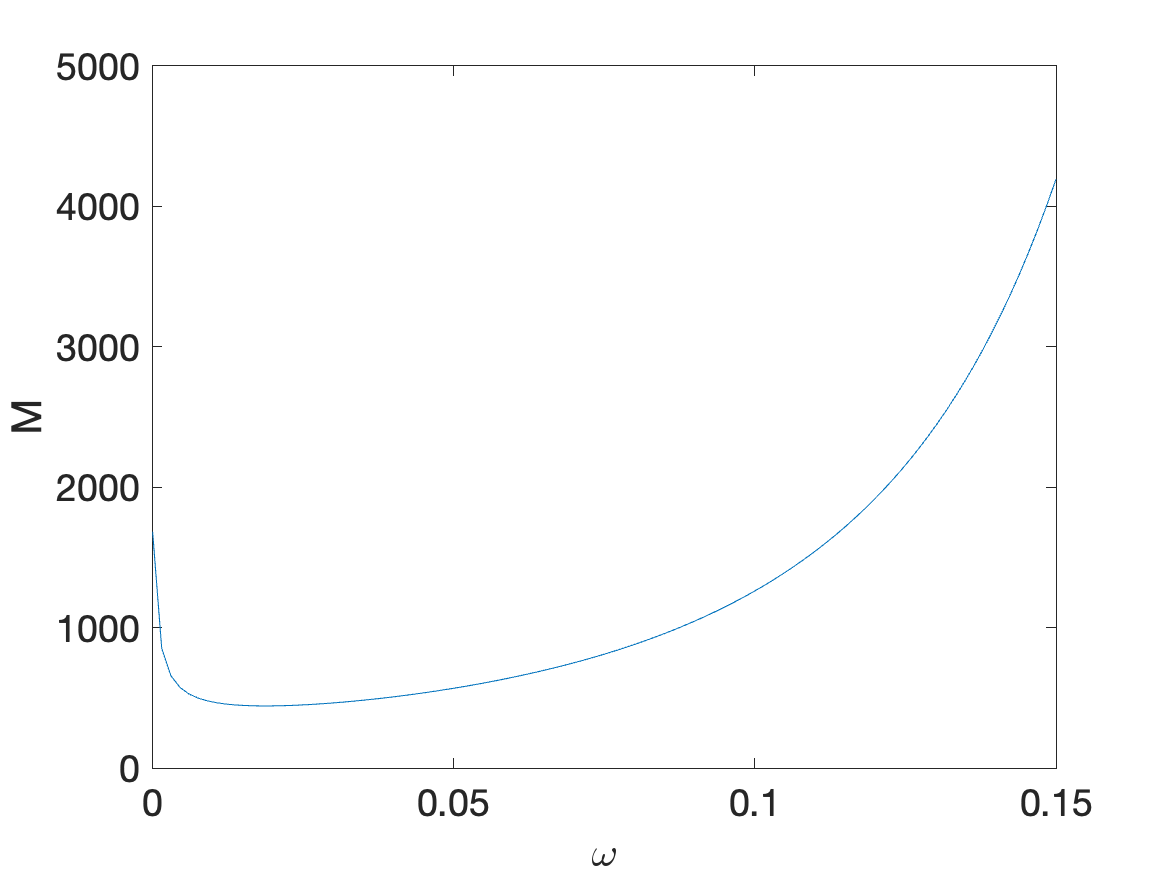}
  \includegraphics[width=0.3\textwidth]{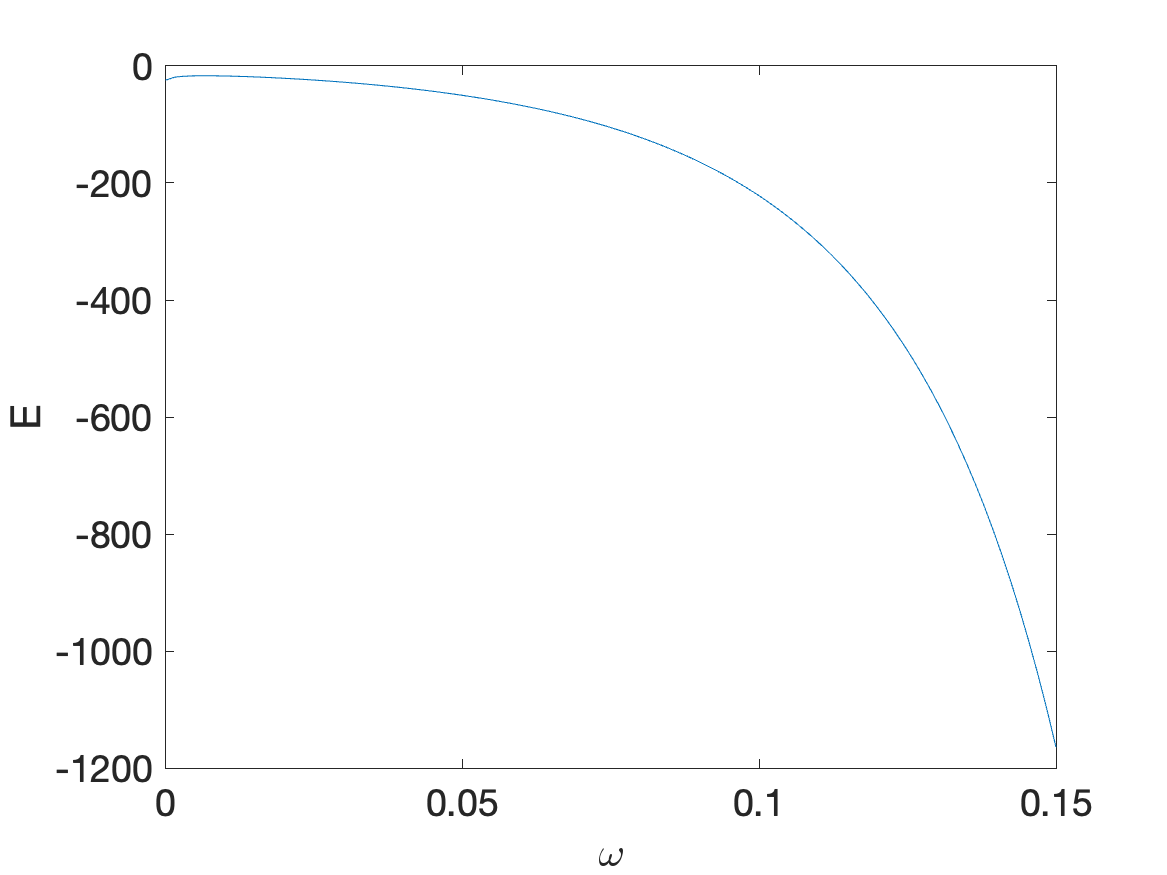}
  \includegraphics[width=0.3\textwidth]{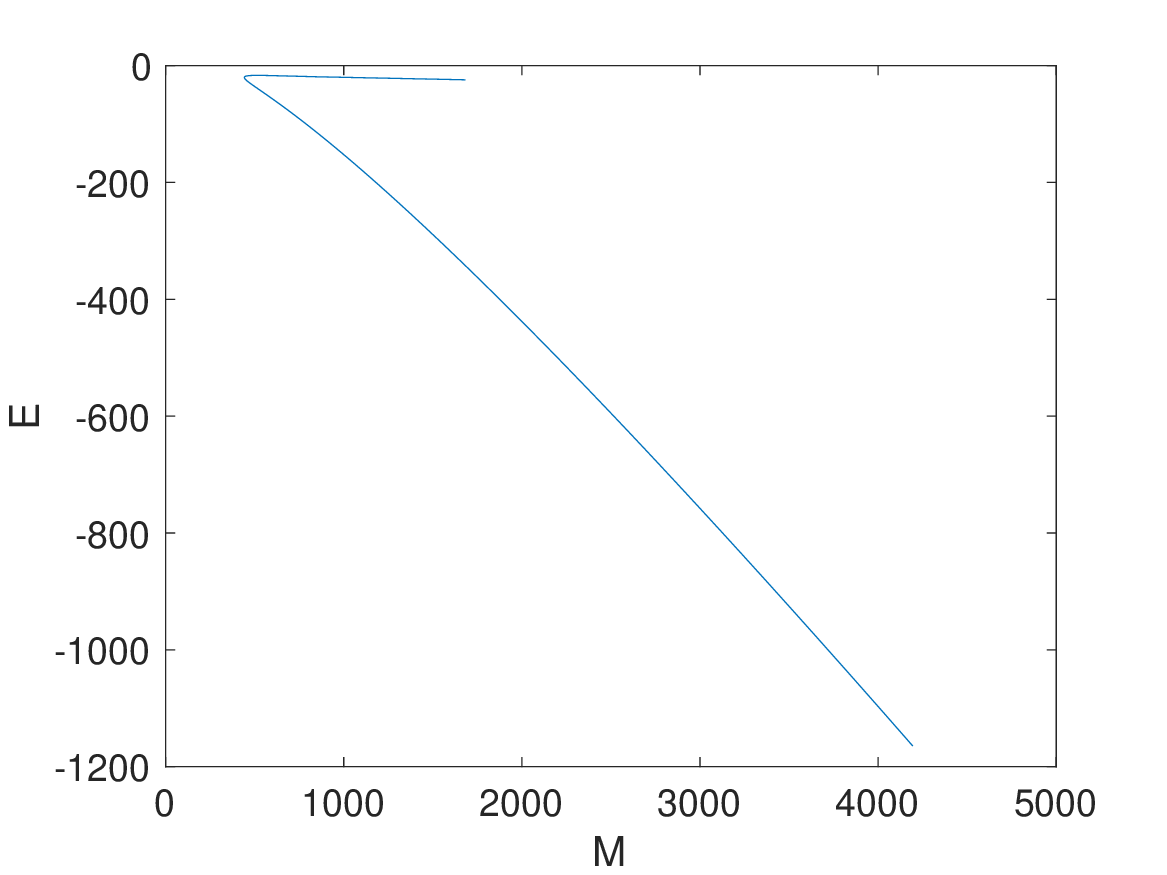}
 \caption{The mass (left) and the energy (middle) for the solitary 
 waves in $d=3$ for $\alpha=3$ in dependence of
 $\omega$, and the energy in dependence of $M$ on the right.}
 \label{fig3dMEalpha3}
\end{figure}

The situation is qualitatively similar for $\alpha=3$ in $d=3$ which
would correspond to an energy supercritical case of the standard NLS 
equation, where no solitary waves exist. Here, there are solitary
waves, and their mass and energy in dependence of $\omega$ can be 
seen in Fig.~\ref{fig3dMEalpha3}.

\section{Numerical study of the time evolution in $d=1$}\label{sect:N1D}

As in \cite{KRN}, we study the stability in the case $\alpha=3$ for 
the two branches with perturbations of the form 
\begin{equation}
	\phi(x,0) = \lambda \varphi_\omega(x),\quad \lambda\approx 1
	\label{pert1d}.
\end{equation}
The numerical approach is the same as in \cite{KRN}, a discrete 
Fourier transform (DFT) in $x$ approximating a situation on $\mathbb{R}$ 
via a periodic setting, $x\in L_{x}[-\pi,\pi]$, $L_{x}>0$, and the 
classical explicit Runge-Kutta method of 4th order in time. The 
accuracy of the numerical results is controlled in $x$ via the decrease  
of the DFT coefficients (we only approximate smooth functions here, 
and for these the highest DFT coefficients indicate the numerical 
error in $x$) and the conservation of the numerically computed mass 
and energy. Both are exactly conserved, but due to unavoidable 
numerical errors, the numerically computed quantities will be time 
dependent. In all numerical experiments the relative error in the conservation of
energy is smaller than $10^{-3}$.  The reader is referred to \cite{KRN} for details. 

\begin{figure}[htb!]
  \includegraphics[width=0.7\textwidth]{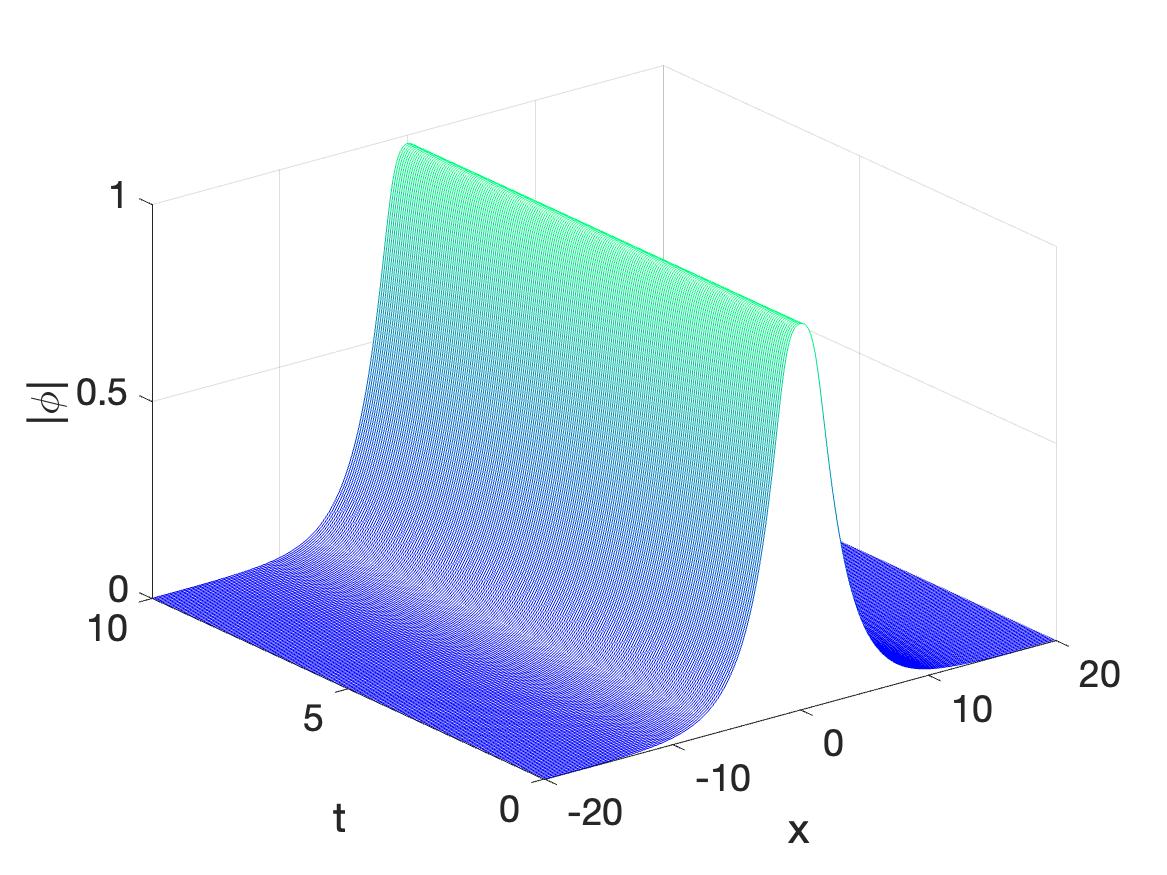}
 \caption{The solution to equation (\ref{SSNLt}) in $d=1$ for the initial data
 $\phi(x,0)=1.001\varphi(x)$ for $\omega=0.22$ and $\alpha=3$.}
 \label{figalpha3_b022_1001}
\end{figure}

First, we consider the case $\omega=0.22$ (recall that
$\omega^{*}=0.25$) which is expected to be stable. We apply 
$N_{x}=2^{12}$ Fourier modes for $x\in 
L_{x}[-\pi,\pi]$ with $L_{x}=30$ and $N_{t}=10^{6}$ time steps for 
$t\leq 10$. The resulting solution for $\lambda=1.001$ can be seen in 
Fig.~\ref{figalpha3_b022_1001}. 

\begin{figure}[htb!]
  \includegraphics[width=0.49\textwidth]{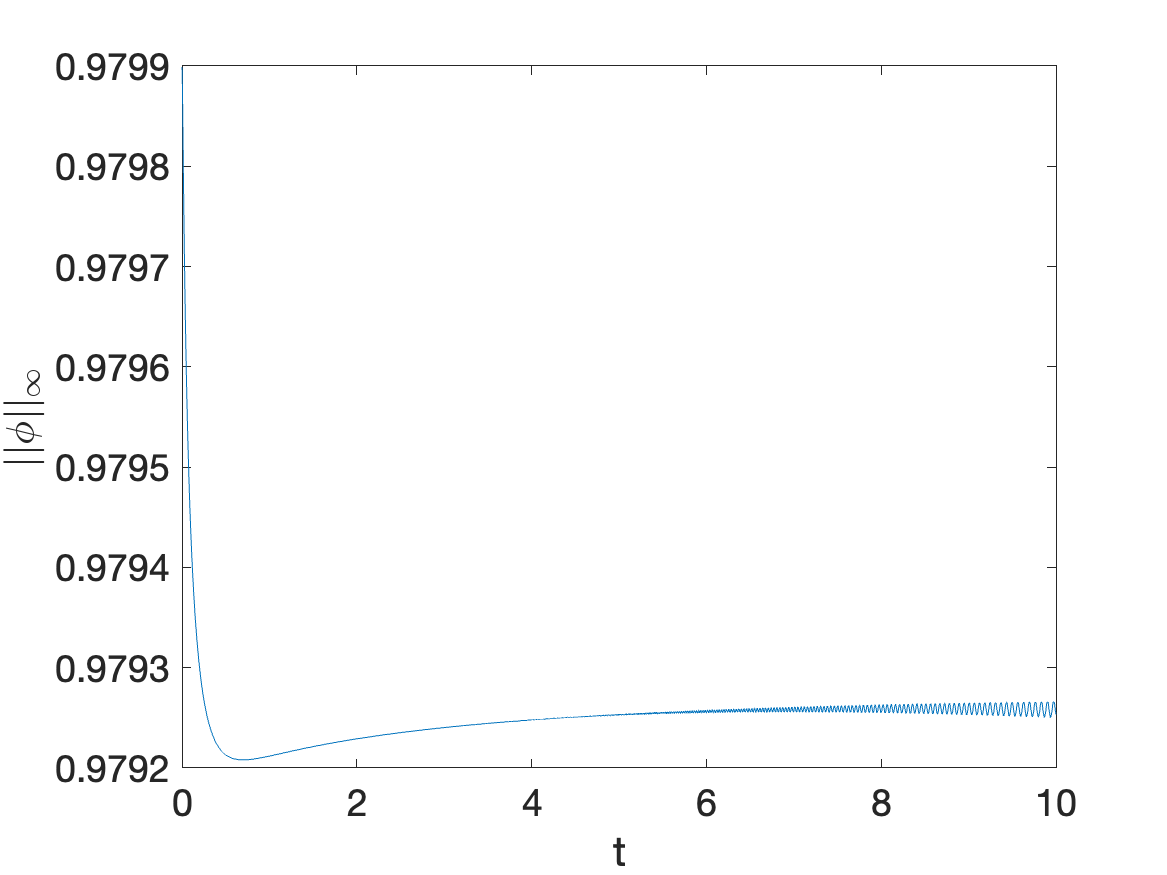}
  \includegraphics[width=0.49\textwidth]{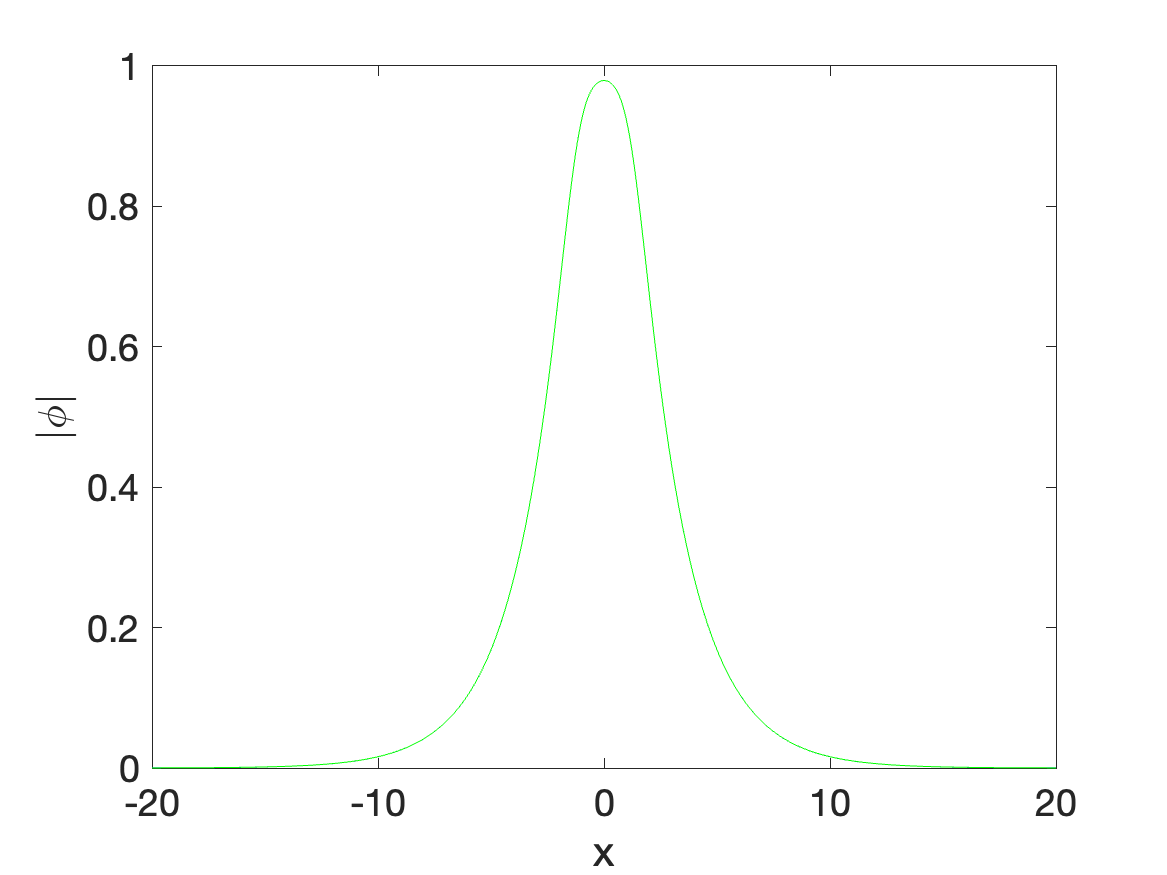}
 \caption{The $L^{\infty}$ norm of the solution to equation (\ref{SSNLt}) in $d=1$ for the initial data
 $\phi(x,0)=1.001\varphi(x)$ for $\omega=0.22$ and $\alpha=3$ on the 
 left, and the solution at the final time in blue with 
 a fitted ground state in green.}
 \label{figalpha3_b022_1001inf}
\end{figure}

There is some radiation emitted at 
the beginning, but then the solution appears to approach a ground 
state with a slightly higher mass than the unperturbed ground state. 
This interpretation is confirmed 
by the $L^{\infty}$ norm of the solution shown on the left in 
Fig.~\ref{figalpha3_b022_1001inf}. The $L^{\infty}$ norm shows some 
oscillations for larger times due to radiation reentering the 
computational domain because of the imposed periodicity, but appears 
to reach a  ground state. 
The oscillations are around an estimated final value which leads via 
(\ref{gsexplicit}) to a fitted value of $\omega\approx0.2205$. The
fitting is based on the maximum of $\varphi_{\omega}$ being given by 
$(\omega/\omega^{*})^{1/(2\alpha)}$, which gives $\omega$ for an 
estimated maximum of the final state. The 
corresponding ground state is shown in green on the left of 
Fig.~\ref{figalpha3_b022_1001inf} together with the solution $\phi$ 
for $t=20$ in blue. The agreement is to the order of $10^{-4}$ which is also 
the level of the radiation in the figure. 

\begin{figure}[htb!]
  \includegraphics[width=0.7\textwidth]{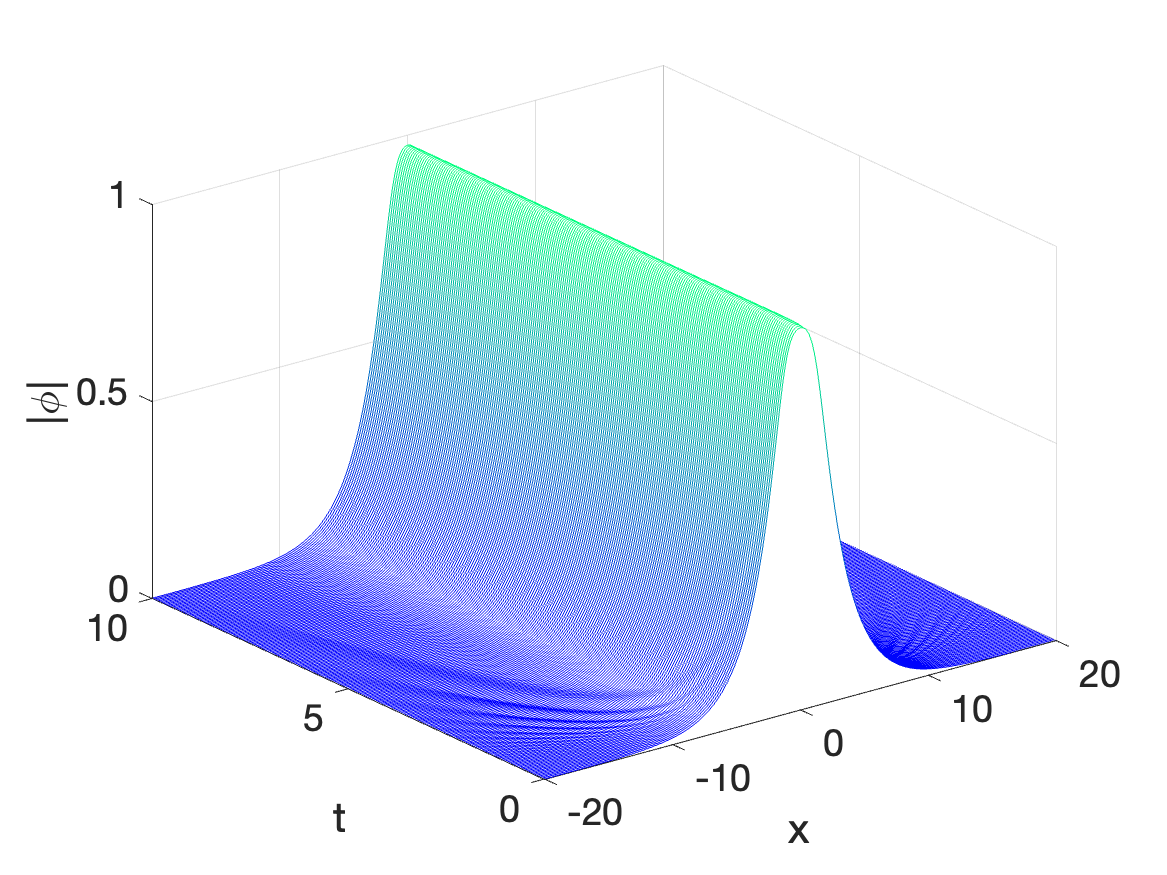}
 \caption{The solution to equation (\ref{SSNLt}) in $d=1$ for the initial data
 $\phi(x,0)=0.99\varphi(x)$ for $\omega=0.22$ and $\alpha=3$.}
 \label{figalpha3_b022_099}
\end{figure}

Thus the solution is clearly stable for these initial data. The 
situation is very similar for $\lambda=0.99$ for the initial data in 
(\ref{pert1d}), i.e., a perturbation with smaller mass than the 
unperturbed ground state. The solution is shown in 
Fig.~\ref{figalpha3_b022_099}.
Once more, radiation is emitted at the beginning, but then the
solution appears to approach a ground state with slightly smaller 
mass than the unperturbed ground state.
The fitted value from the
$L^{\infty}$ norm is $\omega\approx0.2149$, with excellent  agreement (of the order of $10^{-4}$).

\begin{figure}[htb!]
  \includegraphics[width=0.7\textwidth]{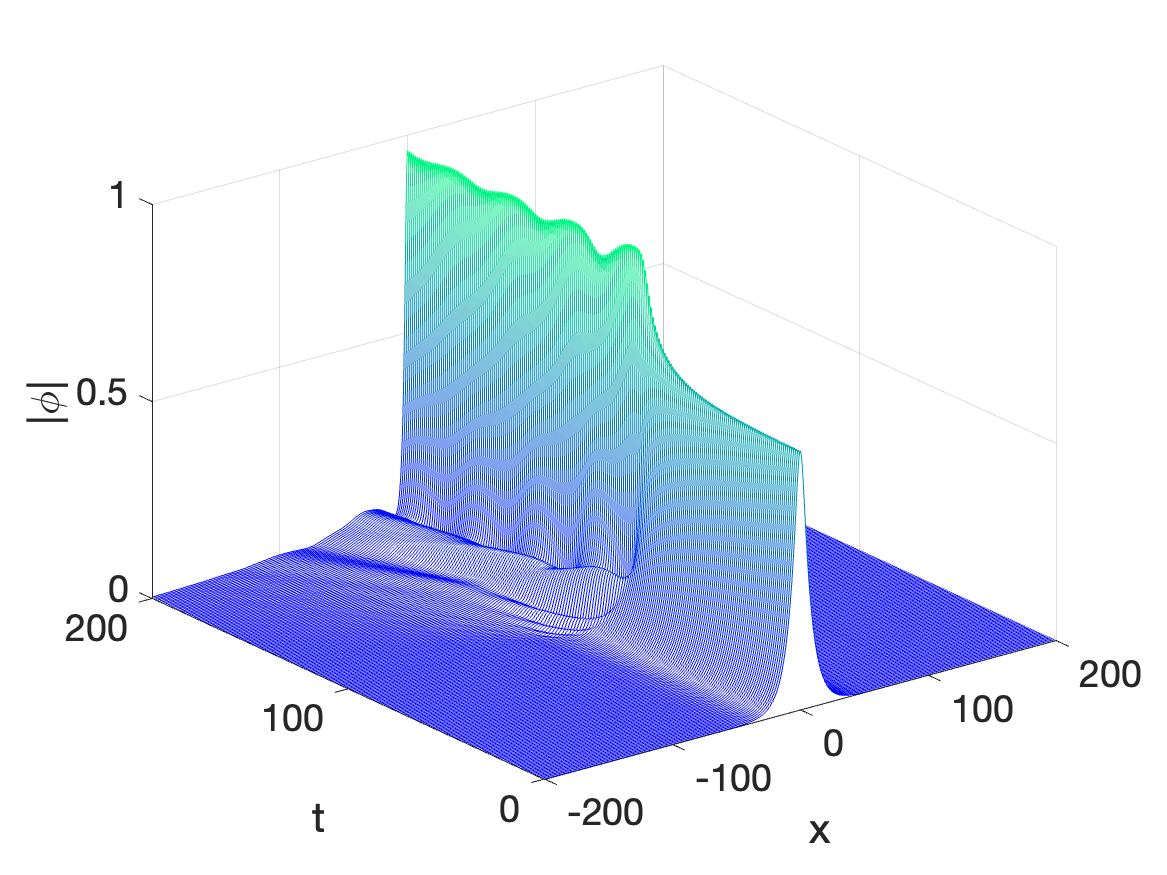}
 \caption{The solution to equation (\ref{SSNLt}) in $d=1$ for the initial data
 $\phi(x,0)=1.001\varphi(x)$ for $\omega=0.02$ and $\alpha=3$.}
 \label{figalpha3_b002_1001}
\end{figure}

In a similar way we study the case $\omega=0.02$ which is expected 
to be unstable. We use 
$N_{x}=2^{12}$ Fourier modes for $x\in 
L_{x}[-\pi,\pi]$ with $L_{x}=100$ and $N_{t}=10^{6}$ time steps for 
$t\leq 200$. Obviously we have to allow for much larger times in the 
case of small $\omega$ and larger domains than above. 
The resulting solution for $\lambda=1.001$ can be seen in 
Fig.~\ref{figalpha3_b002_1001}. For some time the initial peak is 
compressed laterally and growing in height, but then a jump to a 
different form happens. The latter appears to correspond to a ground 
state plus strong radiation.

\begin{figure}[htb!]
  \includegraphics[width=0.49\textwidth]{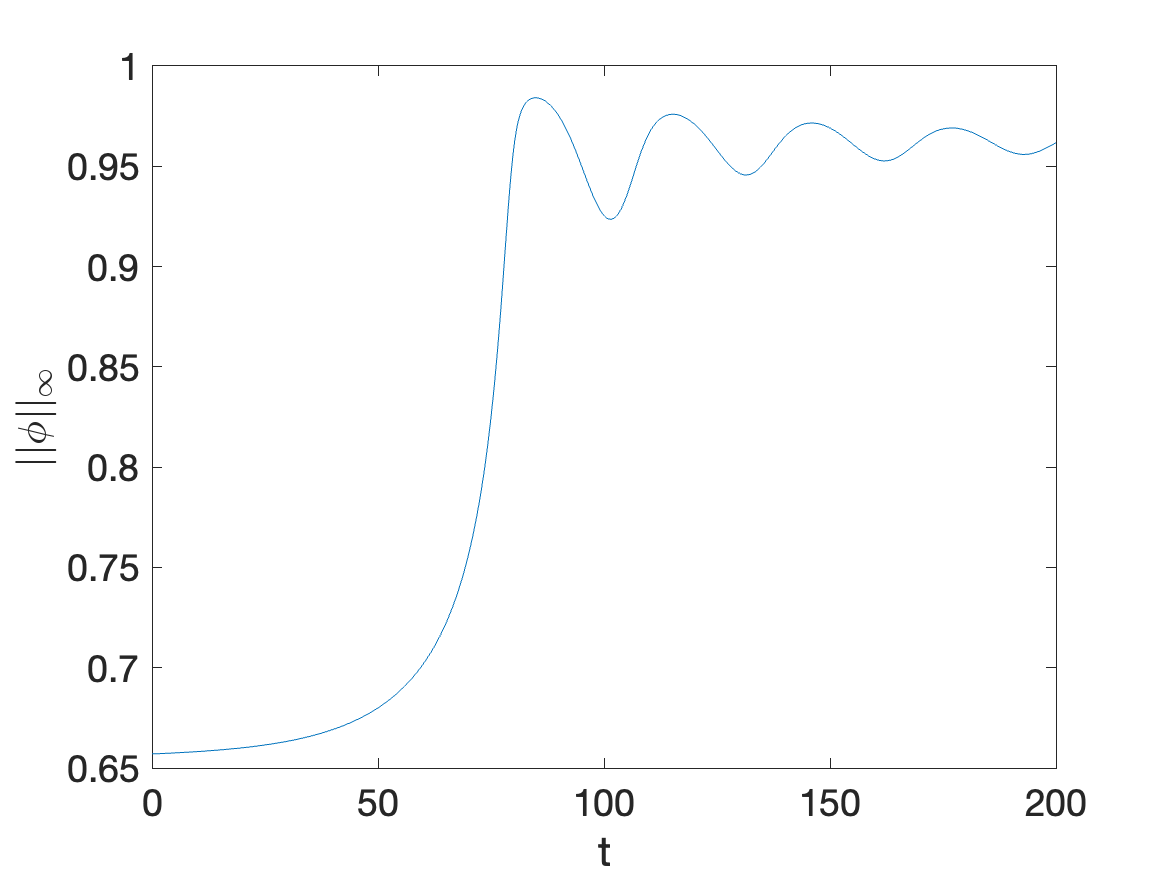}
  \includegraphics[width=0.49\textwidth]{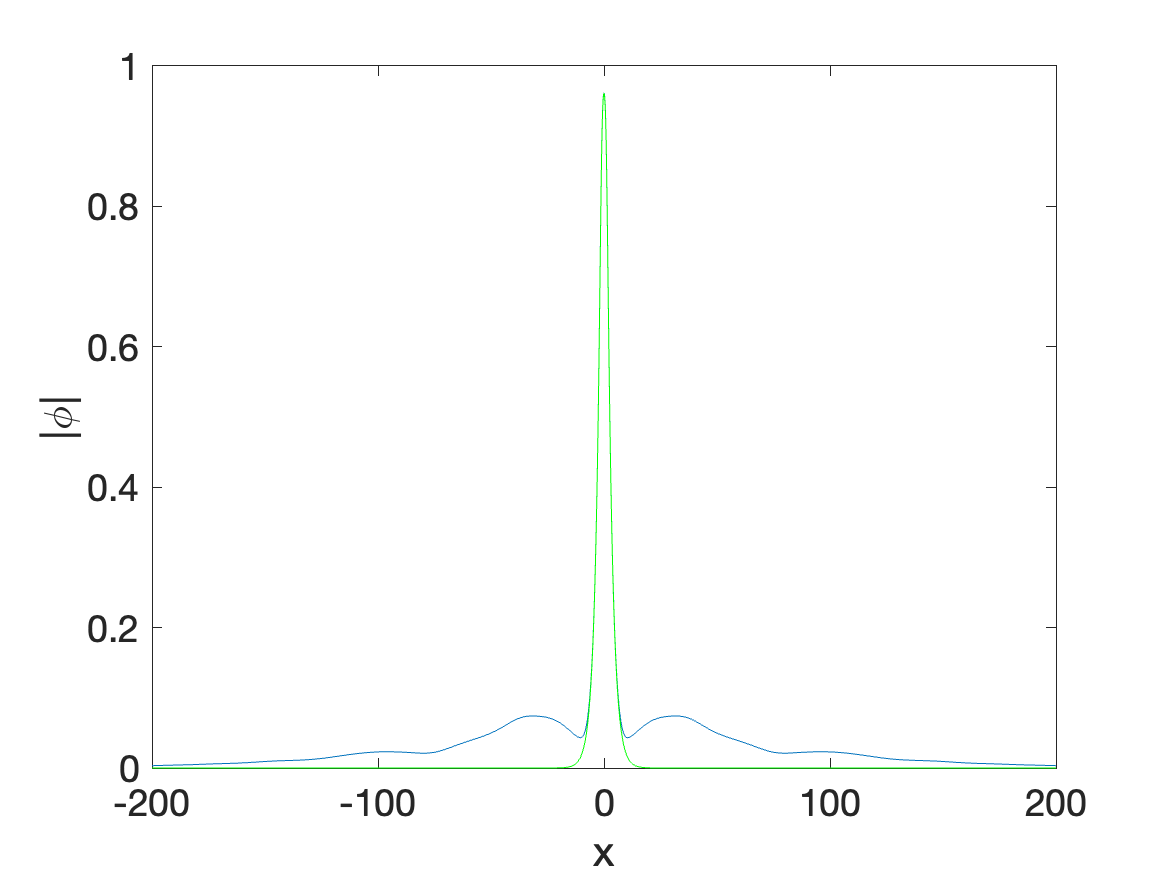}
 \caption{The $L^{\infty}$ norm of the solution to equation (\ref{SSNLt}) in $d=1$ for the initial data
 $\phi(x,0)=1.001\varphi(x)$ for $\omega=0.02$ and $\alpha=3$ on the 
 left, and the solution at the final time in blue with 
 a fitted ground state in green.}
 \label{figalpha3_b002_1001inf}
\end{figure}

The interpretation of the final state of the solution to correspond 
to a ground state is confirmed 
by the $L^{\infty}$ norm of the solution shown on the left in 
Fig.~\ref{figalpha3_b022_1001inf}. It appears to oscillate around an 
asymptotically constant level. The estimated final value of the 
$L^{\infty}$ norm leads via 
(\ref{gsexplicit}) to a fitted value of $\omega\approx0.1972$. The
corresponding ground state is shown in green on the left of 
Fig.~\ref{figalpha3_b002_1001inf} together with the solution $\phi$ 
for $t=200$ in blue. It can be seen that the agreement near the peak 
is excellent, but there is still a considerable amount of radiation 
near the peak at this time. 

\begin{figure}[htb!]
  \includegraphics[width=0.7\textwidth]{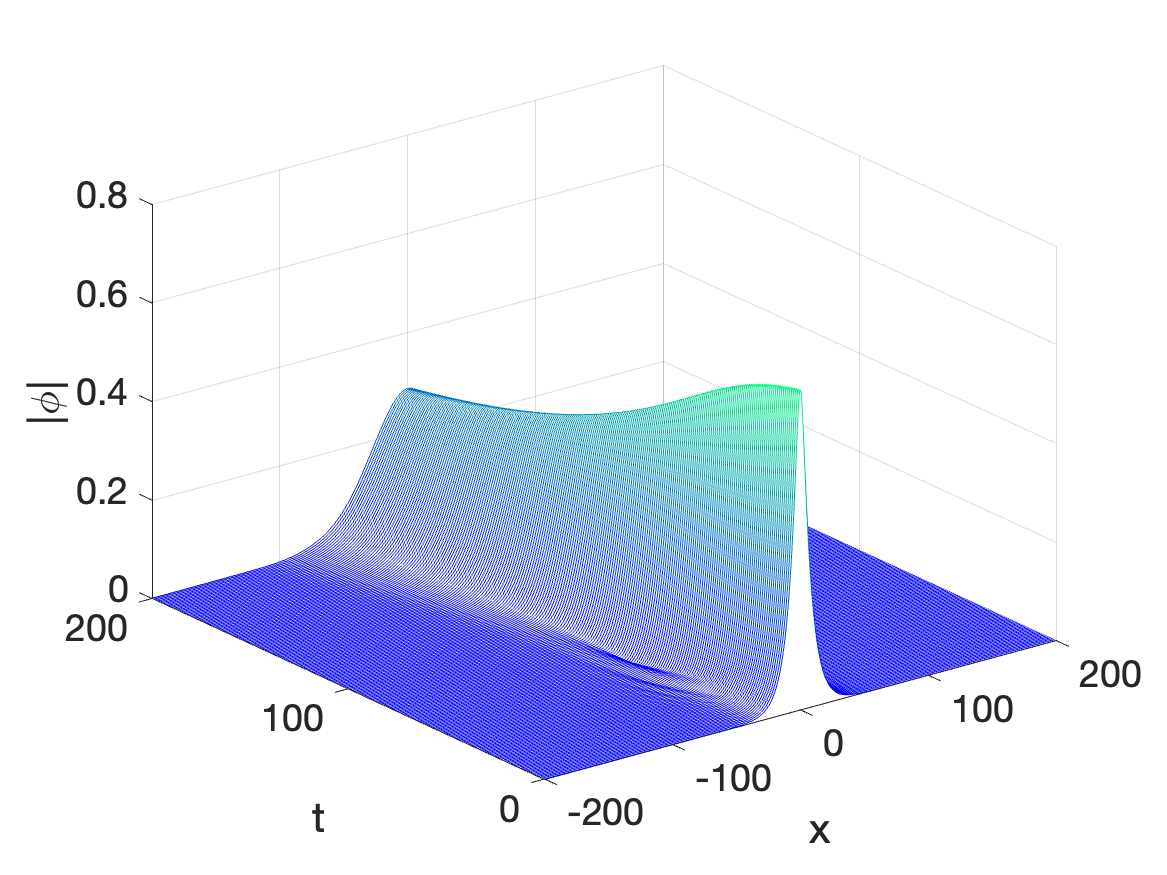}
 \caption{The solution to equation (\ref{SSNLt}) in $d=1$ for the initial data
 $\phi(x,0)=0.99\varphi(x)$ for $\omega=0.02$ and $\alpha=3$.}
 \label{figalpha3_b002_099}
\end{figure}

Thus a perturbation of a stationary solution (\ref{gsexplicit}) on the
unstable branch with larger mass than the unperturbed state is 
unstable against the formation of a ground state on the stable branch. 
The 
situation is completely different for perturbations with smaller mass 
than the unperturbed solution (\ref{gsexplicit}) on the unstable 
branch, say (\ref{pert1d}) for $\lambda=0.99$. The solution is shown in 
Fig.~\ref{figalpha3_b002_099}. 

\begin{figure}[htb!]
  \includegraphics[width=0.49\textwidth]{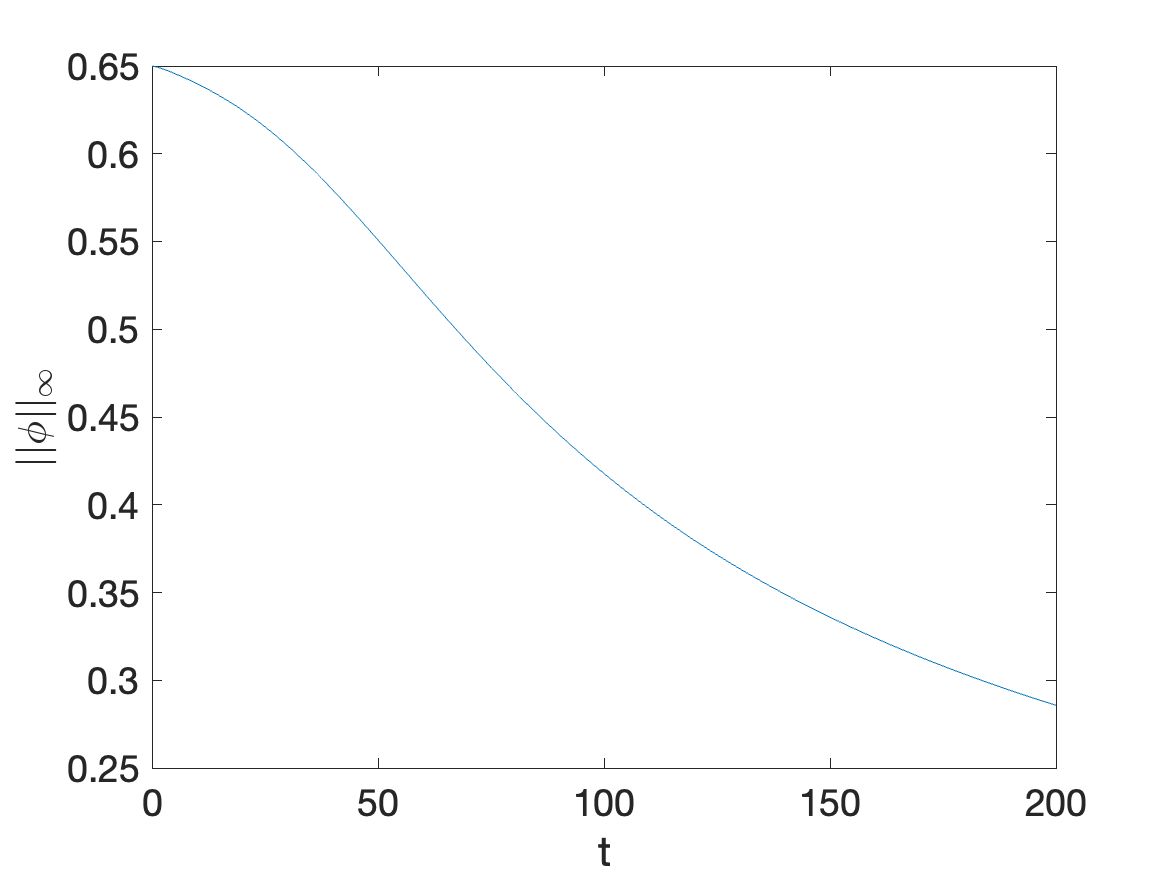}
  \includegraphics[width=0.49\textwidth]{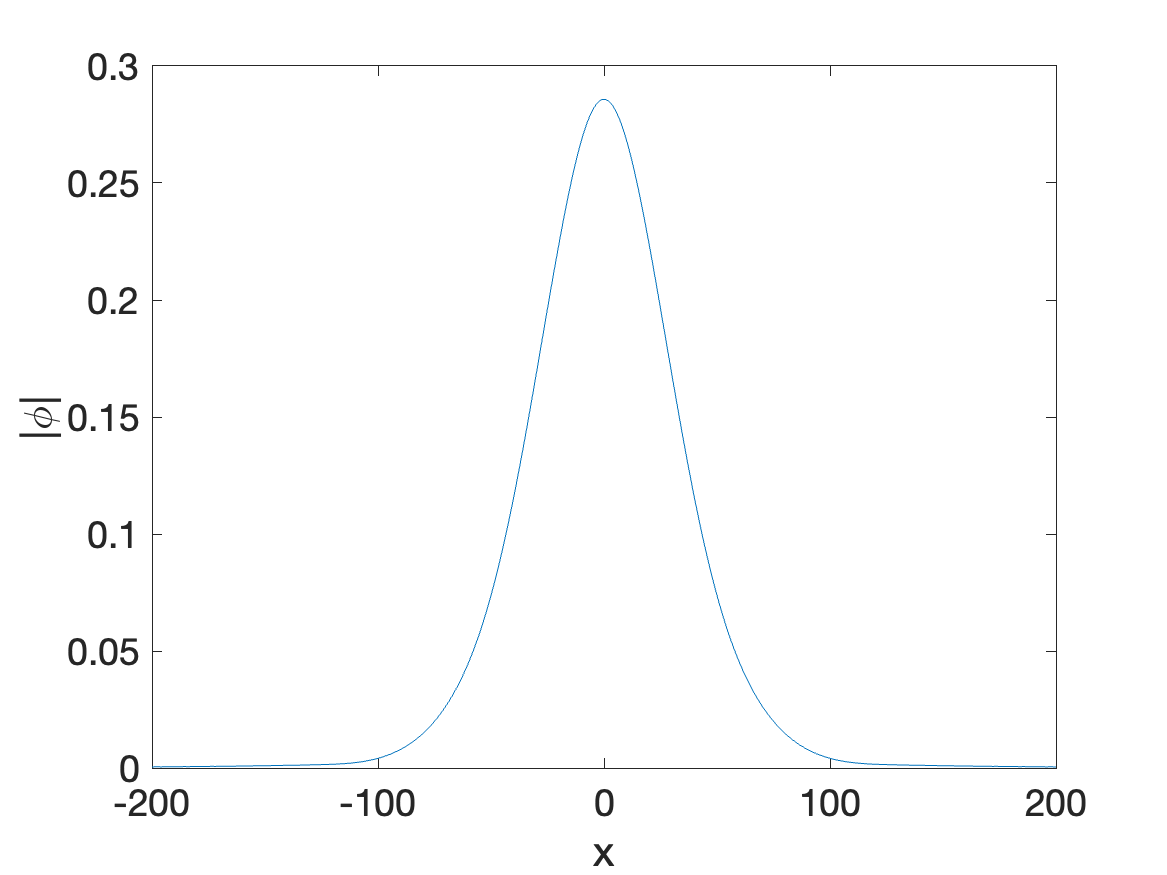}
 \caption{The $L^{\infty}$ norm of the solution to equation (\ref{SSNLt}) in $d=1$ for the initial data
 $\phi(x,0)=0.99\varphi(x)$ for $\omega=0.02$ and $\alpha=3$ on the 
 left, and the solution at the final time.}
 \label{figalpha3_b002_099inf}
\end{figure}

It appears the initial data are simply dispersed in this case, no 
stable ground state seems to appear for long times. 
This interpretation is 
confirmed by the $L^{\infty}$ norm of the solution on the left of 
Fig.~\ref{figalpha3_b002_099inf}. The solution at the final time is 
shown on the right of the same figure. Even for larger times it 
appears to just get broader and flatter. 

\section{Numerical study of solitary waves in higher dimensions}\label{sect:NdD}
In this section we present a numerical approach to equation 
(\ref{SSNLt}) in a radially symmetric situation in higher dimensions. 
First we construct ground state solutions which will illustrate some 
of the results in section 2. Then we will consider perturbations of 
these ground states and the time evolution of 
localized initial data. This provides numerical evidence for the
conjectures. 

\subsection{Construction of solitary waves}

In this section, we numerically construct localized solutions of
equation (\ref{SSNL}) for given $\omega<1/(1+\alpha)$ in dimensions 
greater than 1. Since the 
wanted solution is 
known to be real and radially symmetric with a single maximum at the origin, 
we construct solutions  of the equation
\begin{equation}
	\frac{\Delta\varphi}{1-\varphi^{2\alpha}}-\omega \varphi
	-\alpha\varphi^{2\alpha-2}\frac{(\partial_{r}\varphi)^{2}\varphi}{(1-\varphi^{2\alpha})^{2}} +\phi^{2\alpha+1}=0
	\label{phirad}
\end{equation}
where $\Delta=\partial_{rr}+[(d-1)/r]\partial_{r}$. This will be done as 
in \cite{CKS} with a Newton iteration. 

 As in \cite{CKS} we use the independent variable 
\begin{equation}
	s=r^{2},
	\label{s}
\end{equation}
in order to get a less singular equation. With this coordinate, we 
get for (\ref{phirad})
\begin{equation}\label{Qeqs}
	4s \varphi'' +2d \varphi' 
	+(\varphi^{2\alpha+1}-\omega\varphi)(1-\varphi^{2\alpha})+
	\frac{4\alpha s 
	(\varphi')^{2}\varphi^{2\alpha-1}}{(1-\varphi^{2\alpha})}=0,
\end{equation}
where the prime denotes derivative with respect to $s$. 

Since it is known that the ground states are exponentially 
decreasing, and since we work with finite (double) precision (roughly 
of the order of  $10^{-16}$), 
we solve equation (\ref{Qeqs}) on a finite interval $[0,s_{0}]$, where $s_{0}>0$ is chosen 
such that the solution $\varphi$ vanishes for $s>s_{0}$ with 
numerical precision. 

To approximate the derivatives in (\ref{Qeqs}), we apply \emph{Chebyshev 
differentiation matrices}, see \cite{trefethen}. 
The interval $[0,s_{0}]$ is mapped via $s = \frac{s_{0}}{2}(1+l)$, $l\in[-1,1]$ to the 
interval $[-1,1]$. On the latter we introduce standard \emph{Chebyshev 
collocation points} $l_{n}=\cos(n\pi/N$), $n=0,\ldots,N$, $N\in 
\mathbb{N}$ to discretize 
the problem. For any given $\omega>0$ in the admissible range, the function $\varphi$ is approximated via the {\it Lagrange 
interpolation polynomial} $P_{N}(\ell)$ of degree $N$, coinciding with $\varphi$ at the collocation 
points, 
\[
P_{N}(l_{n}) = \varphi(l_{n}),\quad n = 0,\ldots,N.
\] 
The derivative of $\varphi$ is approximated via the derivative of the Lagrange 
polynomial, i.e.
\[
\frac{\partial}{\partial s} \varphi(s(l_n))\approx P_{N}'(l_n).
\] 
This implies that the derivative is approximated via a 
differentiation matrix $ D $,  see \cite{trefethen,WR} for a 
numerical implementation.

With the above discretization, equation (\ref{Qeqs}) is approximated 
by a system of nonlinear equations for the vector $\varphi$ 
(in an abuse of notation we use the same symbol for 
the function $\varphi$ and its discretisation) which 
can be formally written in the form $\mathbb{F}(\varphi)=0$. The 
vanishing condition for $s=s_{0}$ is implemented as in 
\cite{trefethen} by 
eliminating the column and the line corresponding to $s_{0}$. Thus
$\mathbb{F}(\varphi)=0$ is an $N$-dimensional system of 
nonlinear equations for the $N$ components $\varphi(l_{n})$, 
$n=1,\ldots,N$. 
This system will be solved via a Newton iteration.

Note, however, that equation (\ref{Qeqs}) has the trivial solution 
$\varphi=0$. Therefore, we apply a continuation technique, which means we
construct a non-trivial solution for $\omega\ll 1/(1+\alpha)$ as 
detailed below, and use the result as the initial iterate for a 
slightly larger $\omega$ and so on. 

Below we concentrate on the 
case of a cubic nonlinearity $\alpha=1$ and $d=3$ 
though the code is set up to handle general values (see also the 
examples in the previous subsection for mass and energy plots). We use $N=200$ 
collocation points  and $s_{0}=10^{3}$. For $\omega=0.1$, we take 
$\varphi^{(0)}=0.9\exp(-s/50)$ as the initial iterate and apply some 
relaxation; this means that instead of accepting the new iterate 
$\varphi^{(n+1)}$ of a standard Newton iteration for given iterate 
$\varphi^{(n)}$, we use $\mu\varphi^{(n+1)}+(1-\mu)\varphi^{(n)}$ 
with $0<\mu<1$ as 
the new iterate. The relaxation is needed if the maximum of $\varphi$ 
is close to 1 in which case the iteration can produce values of 
$\varphi$ larger than 1 which leads the iteration to break down. To 
avoid this, we generally use a value of $\mu=0.1$. 
The iteration is stopped once the residual of 
$\mathbb{F}$ is smaller than some threshold, typically $10^{-10}$. 

\begin{figure}[htb!]
\includegraphics[width=0.7\textwidth]{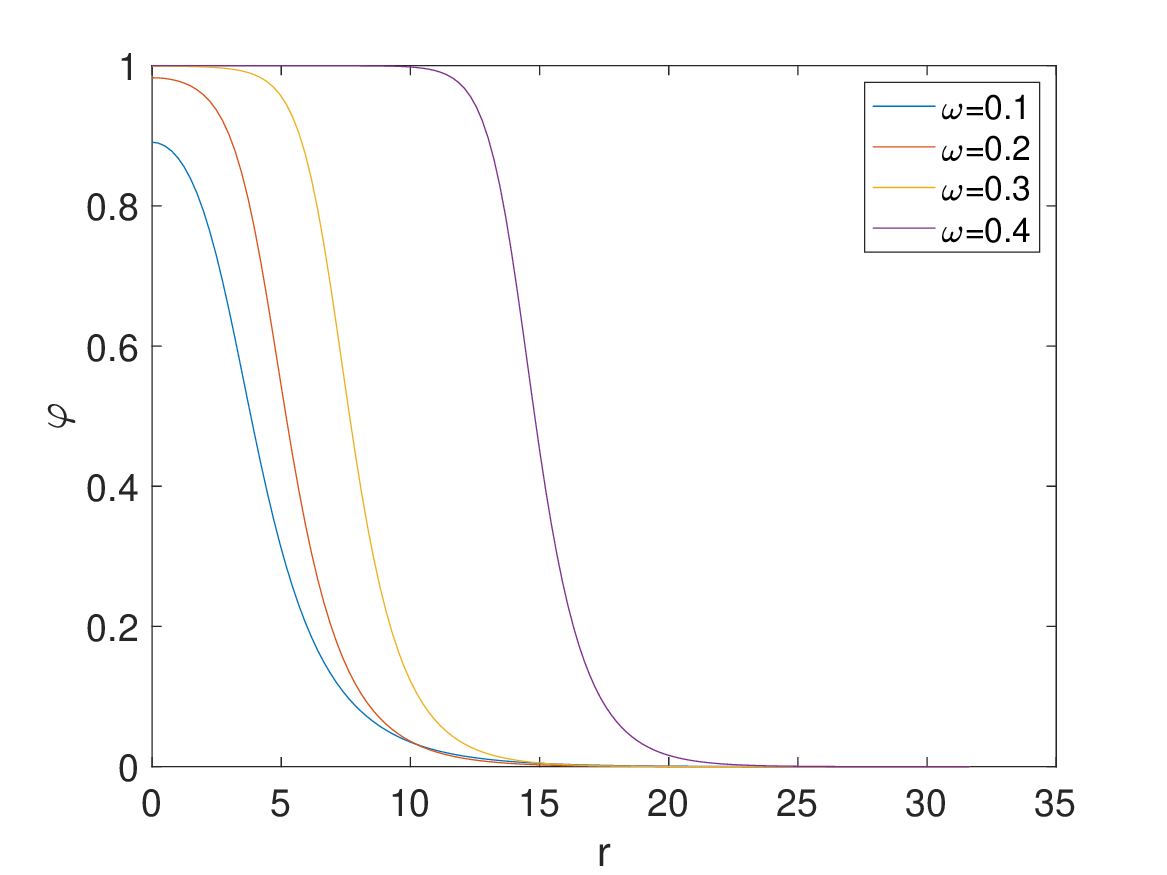}
\caption{Ground state solutions of (\ref{Qeqs}) for several values of 
$\omega$.}
 \label{NLSnuc3Dsol}
\end{figure}

We show the solution for several values of $\omega$ in 
Fig.~\ref{NLSnuc3Dsol}. It can be seen that the maximum of the
solution becomes closer and closer to 1 for larger values of $\omega$ 
(for $\omega=0.4$ the difference between the maximum and 1 is of the 
order of $10^{-7}$). In addition the solution broadens. 

\begin{figure}[htb!]
\includegraphics[width=0.7\textwidth]{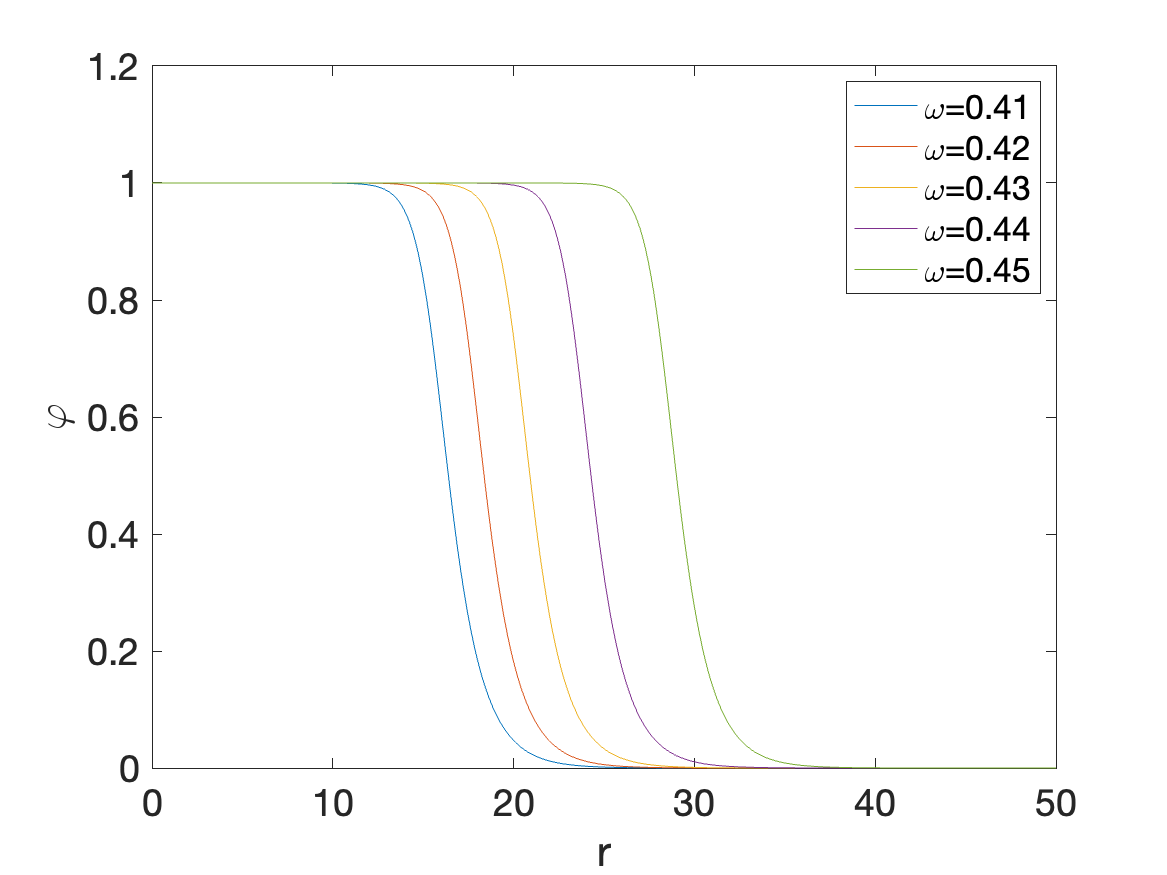}
\caption{Ground state solutions of (\ref{Qeqs}) for several values of 
$\omega>0.4$.}
 \label{NLSnuc3Dsolomlarge}
\end{figure}

To reach even higher values of $\omega$, one has to consider a larger 
interval, $s_{0}=10^{4}$, a higher number of collocation points, 
$N=10^{3}$, a smaller value of the relaxation parameter, 
$\mu=10^{-2}$ and smaller steps in the continuation process. We show 
the solutions for several values of $\omega>0.4$ in 
Fig.~\ref{NLSnuc3Dsolomlarge}. It can be seen that the solutions are
like smoothed steps where the inflection point appears for larger and 
larger values of $r$. For $\omega=0.45$, the difference between the 
maximum and 1 is of the order of $10^{-15}$. Thus it is not possible 
in double precision to reach even larger values of $\omega$ (see Remark \ref{rem_assymptotics}). 

An advantage of the use of Chebyshev collocation points in the 
solution of equation (\ref{Qeqs}) is that the approximation by a 
Lagrange polynomial is related to an expansion of the function 
$\varphi$ in terms of Chebyshev polynomials $T_{n}(l)=\cos(n 
\arccos(l))$, $n=0,1,\ldots$,
\begin{equation}
	\varphi(s(l))\approx \sum_{m=0}^{N}a_{m}T_{m}(l)
	\label{cheb},
\end{equation}
where the \emph{spectral coefficients} $a_{n}$ are determined by a 
collocation method, i.e., as a solution of system
\begin{equation}
	\varphi(s(l_{n}))= \sum_{m=0}^{N}a_{m}T_{m}(l_{n}),\quad 
	n=0,1,\ldots,N
	\label{col}.
\end{equation}
These coefficients can be efficiently computed with a \emph{Fast 
Cosine Transform} that is related to the Fast Fourier transform, see 
the discussion in \cite{trefethen}. As is known for the Fourier 
coefficients of an analytic function, the Chebyshev coefficients 
$a_{n}$ of an analytic function decrease exponentially with $N$. Thus 
the numerical error in truncating a Chebyshev series is of the order 
of the coefficient $a_{N}$. This can be used to estimate the 
numerical error in the absence of a priori estimates. We will always 
control the spatial resolution via the spectral coefficients in this 
paper. 

It is also straight forward to numerically approximate the integrals 
of a function sampled on Chebyshev collocation points with the 
\emph{Clenshaw-Curtis algorithm} \cite{CC}. The integral of a 
function $f(l)$ is approximated via $\int_{-1}^{1}f(l)dl \approx 
\sum_{n=0}^{N}f(l_{n})w_{n}$, where the $l_{n}$, $n=0,\ldots,N$ are 
the Chebyshev collocation points and where the $w_{n}$, 
$n=0,\ldots,N$ are known weights, see \cite{trefethen}. The 
Clenshaw-Curtis method is again a \emph{spectral method}, i.e., the 
error in the numerical approximation of an analytic function 
decreases exponentially with $N$. With this algorithm one can 
compute the 
$L^{2}$-mass  and energy (\ref{ener}) of the ground states that
are shown in Figs.~\ref{fig2dMEalpha1} to \ref{fig3dMEalpha3}.

\subsection{Numerical approach for the time evolution in the radially
symmetric case}
In this section we detail the numerical approach for the time 
integration of  equation (\ref{SSNLt}) in the radially symmetric 
case. 

We apply the same approach for the spatial dependence as in the 
previous subsection, i.e., a Chebyshev grid in the variable $s=r^{2}$. 
The time integration will be done with a Crank-Nicolson (CN) method, a 
second order implicit scheme which is unconditionally stable (this 
means the time step can be chosen in dependence of the aimed at 
accuracy only, without stability considerations). We
discretize time in the form $t_{n}=nh$, $n=0,1,\ldots$, $h>0$ the 
time step (of course it is possible to use different values for $h$ 
in each time step, but for simplicity we choose a constant $h$ for each 
run). For an 
equation of the form $\partial_{t}\phi = \mathcal{F}(\phi)$ the CN 
scheme then  reads
\begin{equation}
	\phi(t_{n+1})-\phi(t_{n}) = \frac{h}{2}[\mathcal{F}(\phi(t_{n}))
+\mathcal{F}(\phi(t_{n+1})]
	\label{CN}.
\end{equation}
For equation (\ref{SSNLt}), we get after the spatial discretisation 
of the previous subsection an $N+1$ dimensional system of ordinary 
differential equations (ODEs) that will be integrated with 
(\ref{CN}) which is equivalent to a nonlinear equation system for 
$\phi(t_{n+1})$ (we again use the same symbol for the function $\phi$ 
and the vector $\phi$,  thus 
$\phi(t_{n+1})\in\mathbb{C}^{N+1}$ here). 

In \cite{CKS} we have solved the resulting system with a simplified 
Newton iteration, i.e., an inversion of the linear part of the 
equation containing the Laplacian. This is not possible here since 
the Laplacian appears in a nonlinear term. Since the term 
$1-|\phi|^{2\alpha}$ in the denominator of 
(\ref{SSNLt}) can be become small, an efficient iteration that is 
rapidly converging has to be applied. Therefore, we use as in the
previous subsection a Newton iteration. The problem is that in contrast 
to the stationary solutions $\varphi$, the general solution $\phi$ of 
(\ref{SSNLt}) will be complex whereas equation (\ref{SSNLt}) is not 
complex linear in $\phi$. Therefore, we have to consider also its
complex conjugate,
\begin{equation}
\label{SSNLtcc}
-i\partial_{t}\bar{\phi} = -\nabla \left(\frac{\nabla\bar{\phi}}{1-|\phi|^{2\alpha}} \right) +\alpha |\phi|^{2\alpha-2}\frac{|\nabla\phi|^2}{(1-|\phi|^{2\alpha})^2}
\bar{\phi} - |\phi|^{2\alpha}\bar{\phi} 
\end{equation}
with the same spatial discretisation. The CN approach (\ref{CN}) 
applied to both these equations leads to a nonlinear system $F(\Phi)=0$ 
where $\Phi$ is the $2N+2$ dimensional vector formed by 
$\phi(t_{n+1})$ and its complex conjugate. Equation $F(\Phi)=0$ is solved by a 
standard Newton iteration which is stopped once the residual 
$||F||_{\infty}<10^{-10}$. The initial iterate for time step 
$t_{n+1}$ is always given by $\phi(t_{n})$.

The accuracy of  the spatial discretisation is controlled as in the 
previous subsection via the Chebyshev coefficients (\ref{col}). The 
resolution in time is as in the 1D case via the conservation 
of the quantity $\delta:=|E(t)/E(0)-1|$. 
Since it was shown  in \cite{etna} that this quantity overestimates 
the numerical accuracy by 1-2 orders of magnitude, we  
always aim at values of $\delta<10^{-3}$. 

\begin{figure}[htb!]
  \includegraphics[width=0.49\textwidth]{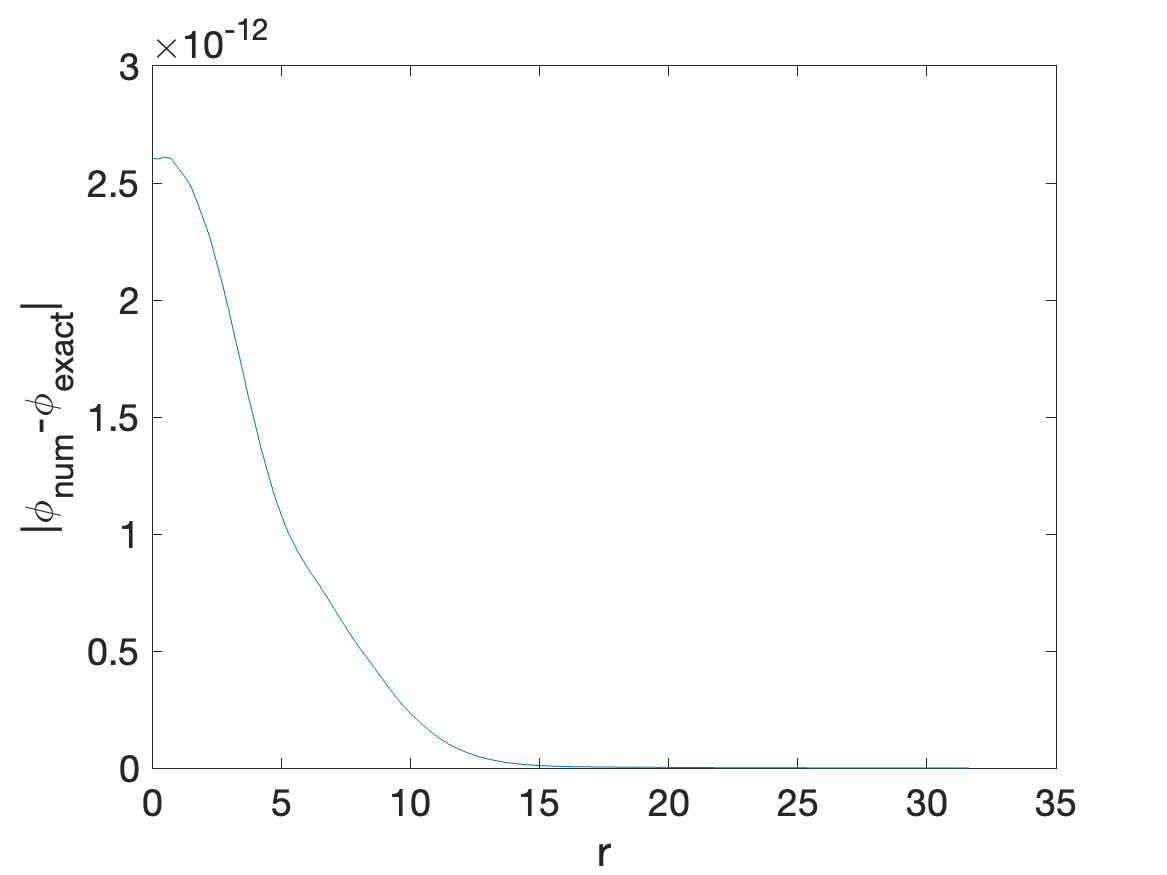}
  \includegraphics[width=0.49\textwidth]{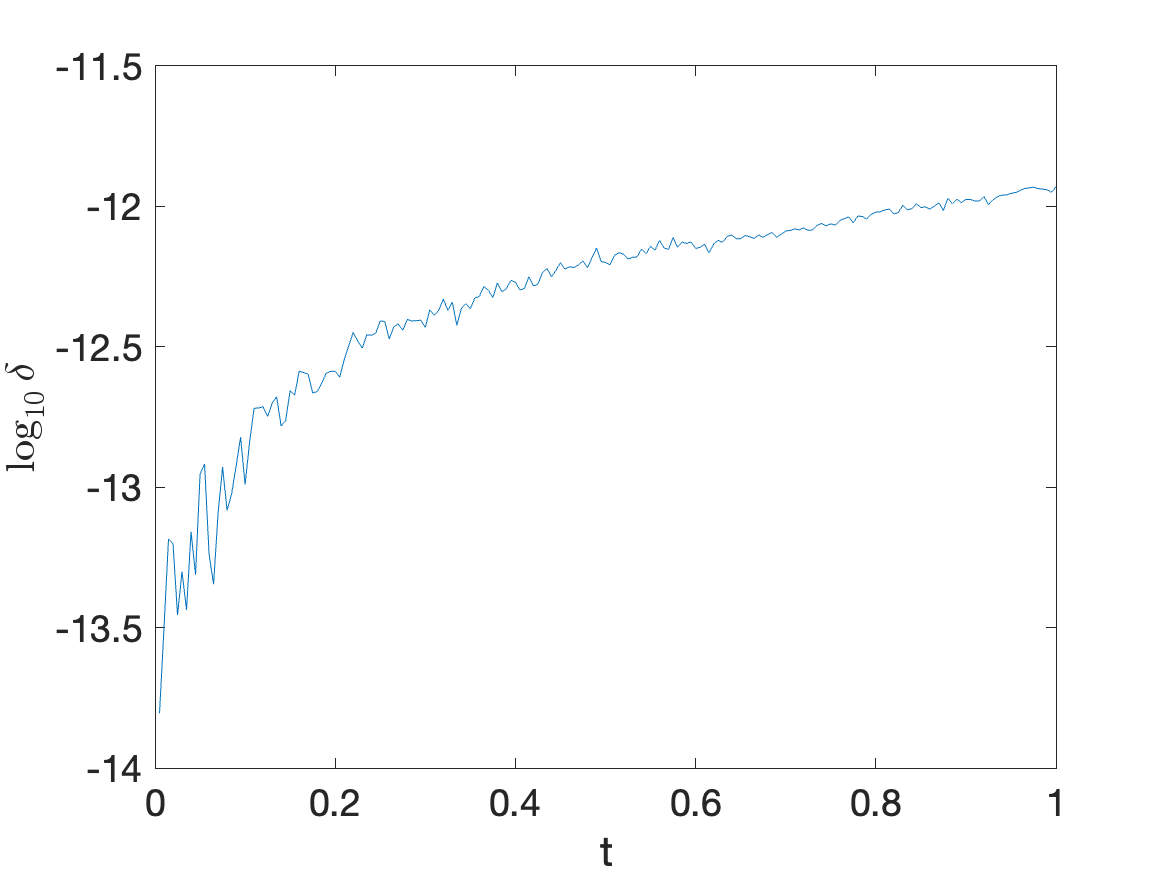}
 \caption{Time evolution of the stationary solution $\varphi$ for 
 $\omega=0.1$ constructed in the previous subsection: on the left the 
 difference between numerical and exact solution for $t=1$, on the 
 right the relative energy $\delta$ in dependence of time.}
 \label{soltest}
\end{figure}

To test the code, we use the stationary solution of the previous 
subsection for $\omega=0.1$ as initial data with $N_{t}=4000$ time steps 
for $t\leq 1$. The difference between the numerical and the exact 
solution $\phi_\omega=\varphi_\omega e^{i\omega t}$ at the final time is seen in
Fig.~\ref{soltest} on the left. It is of the order of $10^{-12}$. 
Note that this tests also the solution $\varphi$ constructed in the 
previous subsection since an error there would also appear in the time 
evolution.  On the right side of the same figure we show the relative 
energy conservation. It can be seen that it is of the order of 
$10^{-12}$ during the whole computation. 

\subsection{Perturbed ground states}
In this subsection we study numerically the evolution of radial 
perturbations of the solutions to the stationary 
equation~\eqref{SSNL}.  As in 1D, we consider initial data of the form
\begin{equation}
	\phi(r,0) = \lambda \varphi_\omega(r),
	\label{solpert}
\end{equation}
where $\lambda$ is a real constant close to 1. Its value is chosen 
such that $||\phi(r,0)||_{\infty}<1$. 

As an example we consider the solution for $\omega=0.1$, which is
expected to be stable,  with 
$\lambda=1.01$ and $\lambda=0.99$, i.e., a perturbation of the order 
of 1\%. 
Numerically one has to apply a finite perturbation in order to see an 
effect in finite time. This will lead to some radiation which will 
propagate towards infinity. To reduce the effect of this radiation being 
reflected at the boundary of the numerical domain, we choose a larger 
domain, $s\leq 10^{4}$. In the examples studied below, the effect of the 
radiation will be of the order of $10^{-6}$. If longer times are to 
be studied, one either has to use larger domains or several domains, 
possibly with a compactification of the real line as in \cite{birem}. 
Here we apply $N=1000$ Chebyshev collocation points for 
$s\in[0,10^{4}]$ and $N_{t}=10^{4}$ time steps for $t\leq 20$. 

\begin{figure}[htb!]
  \includegraphics[width=0.49\textwidth]{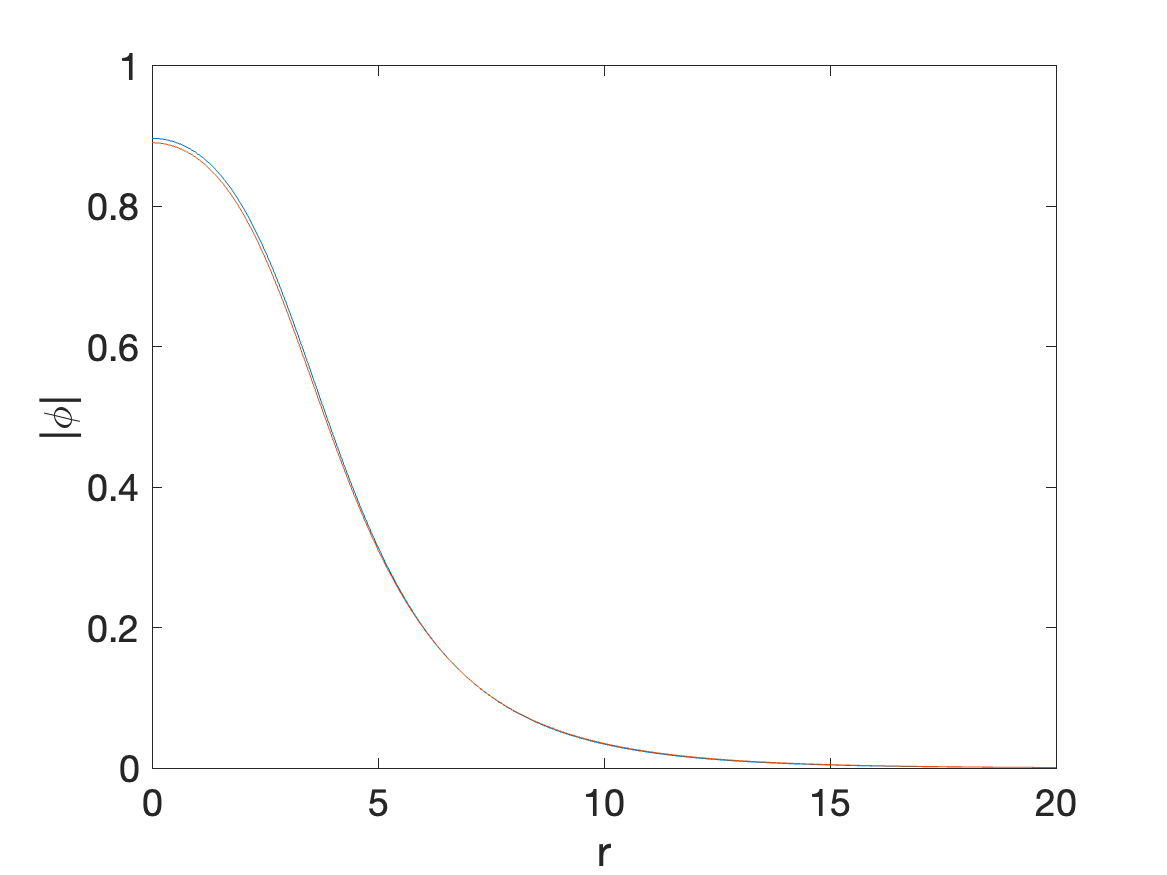}
  \includegraphics[width=0.49\textwidth]{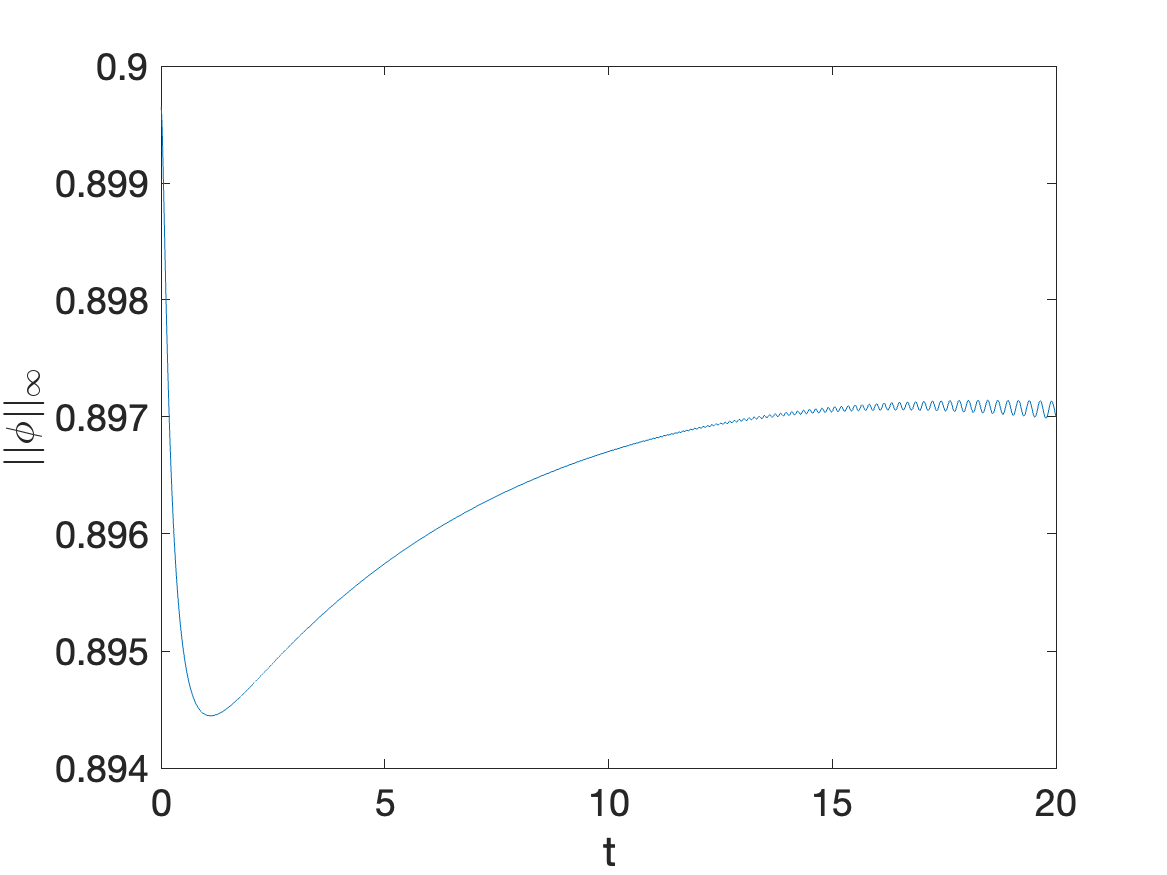}
 \caption{Solution to equation (\ref{SSNLt}) for initial data 
 (\ref{solpert}) with $\lambda=1.01$: on the left the 
  solution for $t=20$ in blue and the initial condition in red, on the 
 right the $L^{\infty}$ norm of the solution in dependence of time.}
 \label{NLS3Dsolom01_101_t20}
\end{figure}

For $\lambda=1.01$ in (\ref{solpert}), we get the solution at the 
final time shown in Fig.~\ref{NLS3Dsolom01_101_t20} on the left. The
initial condition is shown in red in the same figure. It appears that 
the final state is a ground state with a slightly larger mass and 
maximum than the unperturbed ground state. This interpretation is 
backed by the $L^{\infty}$ norm of the solution on the right of the 
same figure. It appears to saturate at some point which indicates 
that a stationary solution is the final state. The small oscillations 
in the $L^{\infty}$ norm are due to some back scattering of 
radiation. 

\begin{figure}[htb!]
  \includegraphics[width=0.7\textwidth]{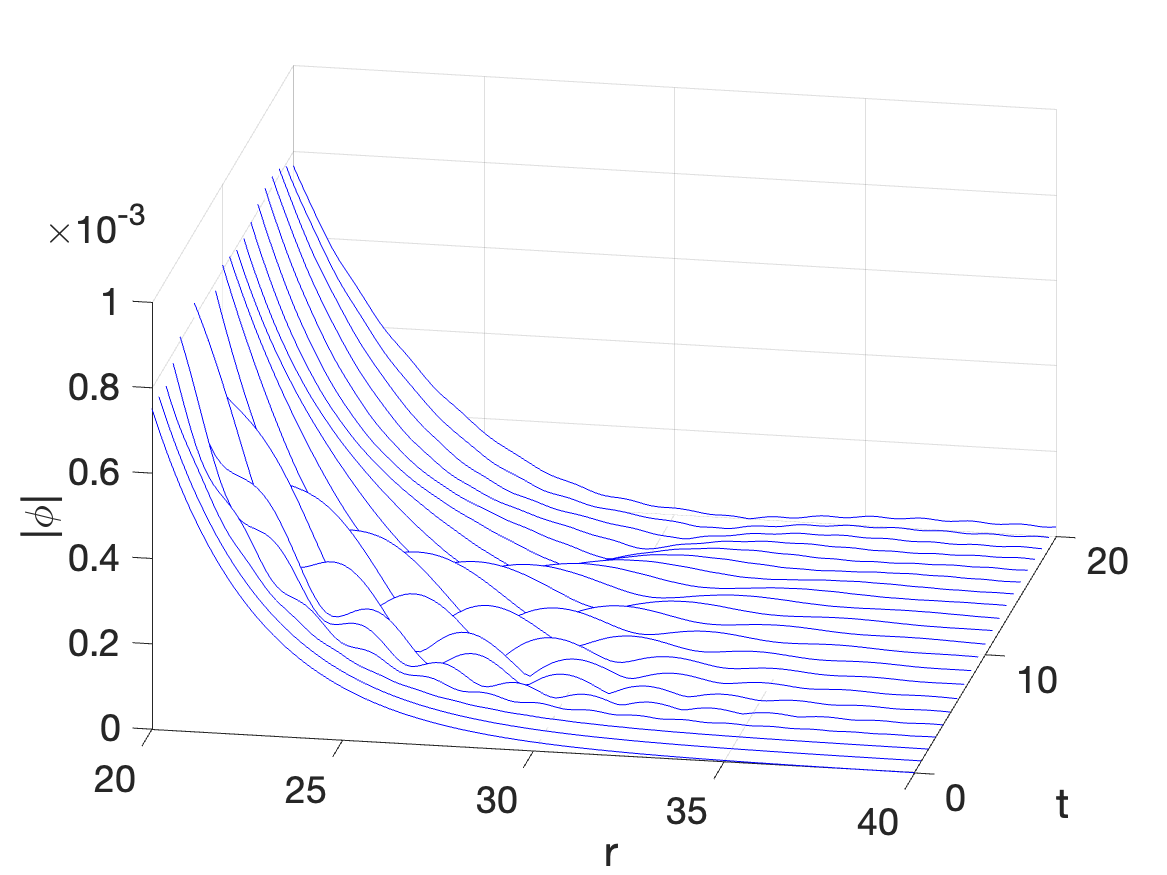}
 \caption{Close-up of the solution to equation (\ref{SSNLt}) for initial data 
 (\ref{solpert}) with $\lambda=1.01$.}
 \label{NLS3Dsolom01_101watercu}
\end{figure}

This radiation is hardly visible in Fig.~\ref{NLS3Dsolom01_101_t20}.
Therefore, we show a close-up of the solution in
Fig.~\ref{NLS3Dsolom01_101watercu}. It can be seen that radiation is
emitted towards infinity, but that the final state appears to be 
another ground state. Since the ground states are not explicitly known as in $d=1$,
it is less straight forward to identify the value of $\omega$ corresponding to 
the final state than in 1D (one would have to compute the 
$L^{\infty}$ norm for a sufficient number of values of $\omega$ in 
order to identify as in 1D the $\omega$ corresponding to a certain 
final state). 

\begin{figure}[htb!]
  \includegraphics[width=0.49\textwidth]{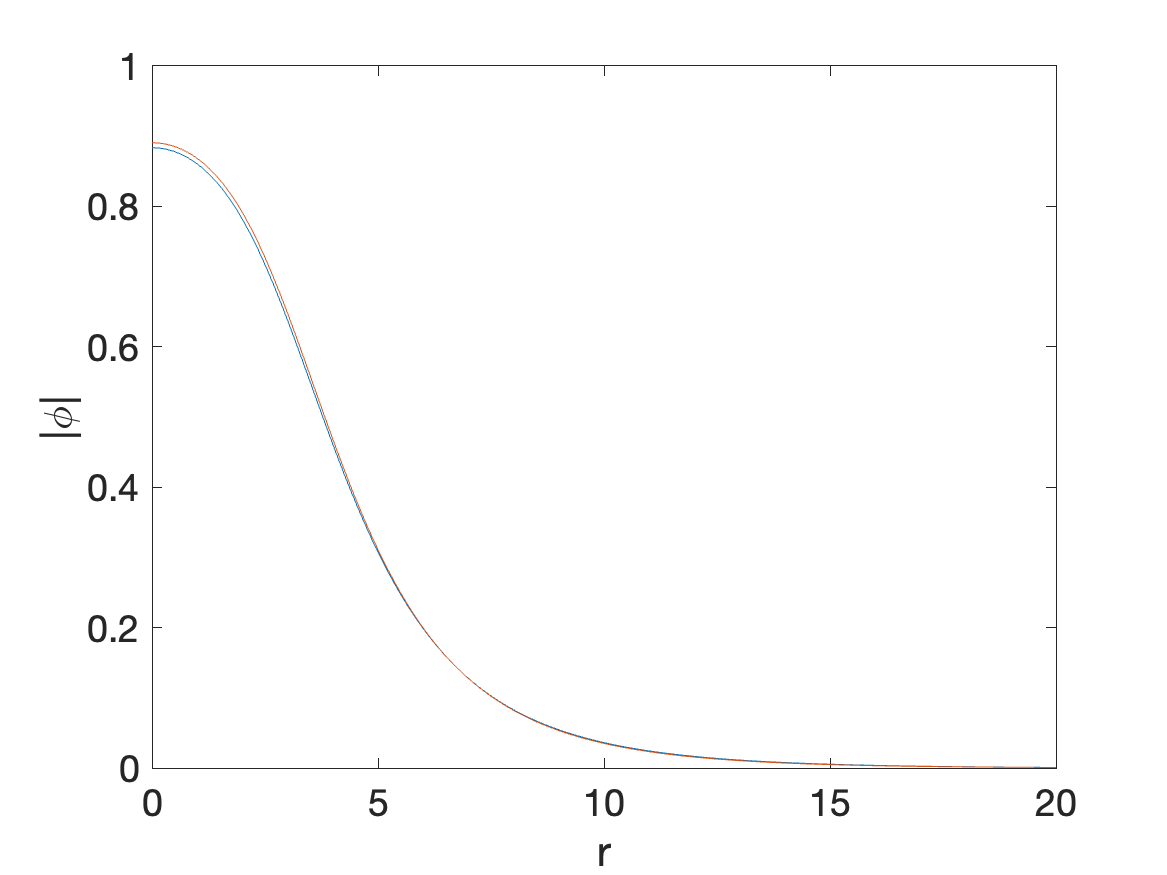}
  \includegraphics[width=0.49\textwidth]{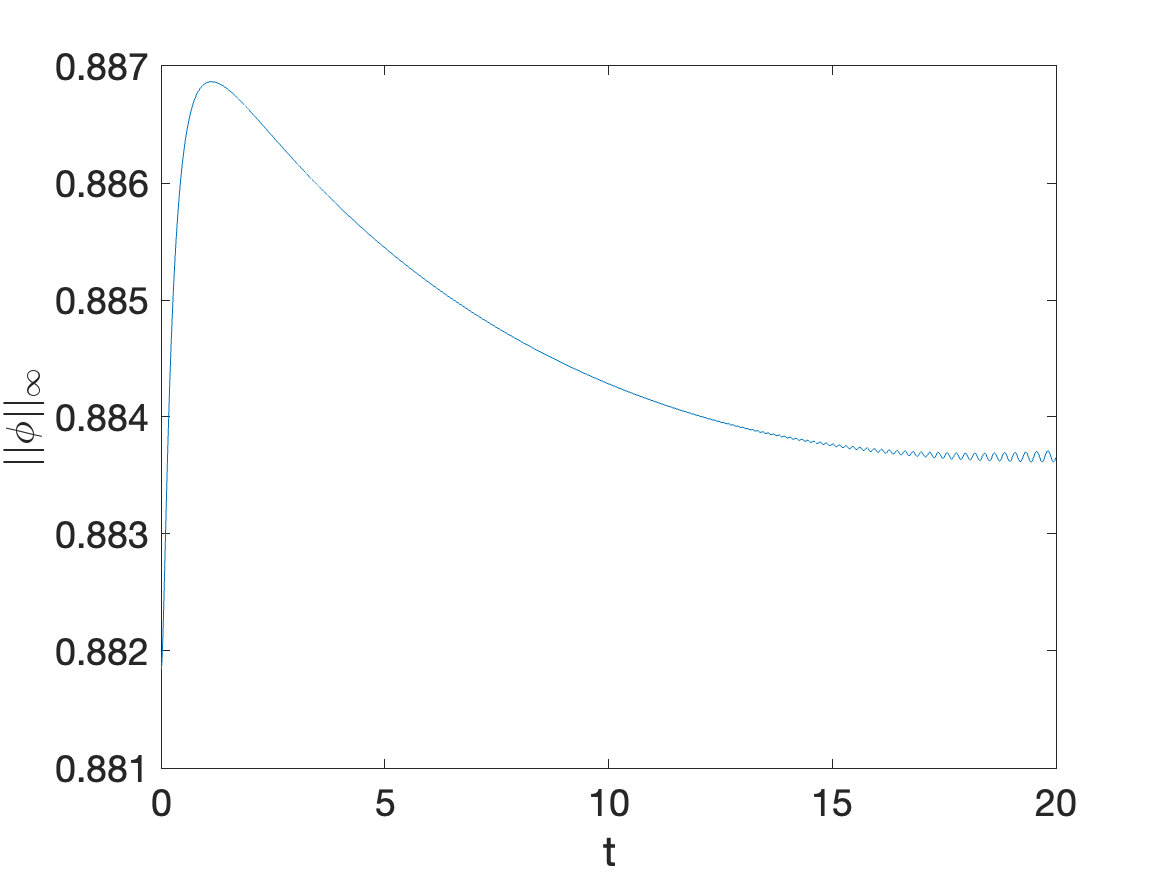}
 \caption{Solution to equation (\ref{SSNLt}) for initial data 
 (\ref{solpert}) with $\lambda=0.99$: on the left the 
  solution for $t=20$ in blue and the initial condition in red, on the 
 right the $L^{\infty}$ norm of the solution in dependence of time.}
 \label{NLS3Dsolom01_099_t20}
\end{figure}

The situation is similar for initial data of the form (\ref{solpert}) 
with $\lambda=0.99$, i.e., a perturbation with smaller mass than the 
unperturbed ground state. The solution at the final time can be seen 
in Fig.~\ref{NLS3Dsolom01_099_t20} on the left together with the
initial data in red. The final state appears to be a ground state of 
slightly lower mass than the unperturbed ground state. This is in 
accordance with the $L^{\infty}$ norm of the solution on the right of 
the same figure. Thus the ground state appears to be asymptotically 
stable. 

\begin{figure}[htb!]
  \includegraphics[width=0.7\textwidth]{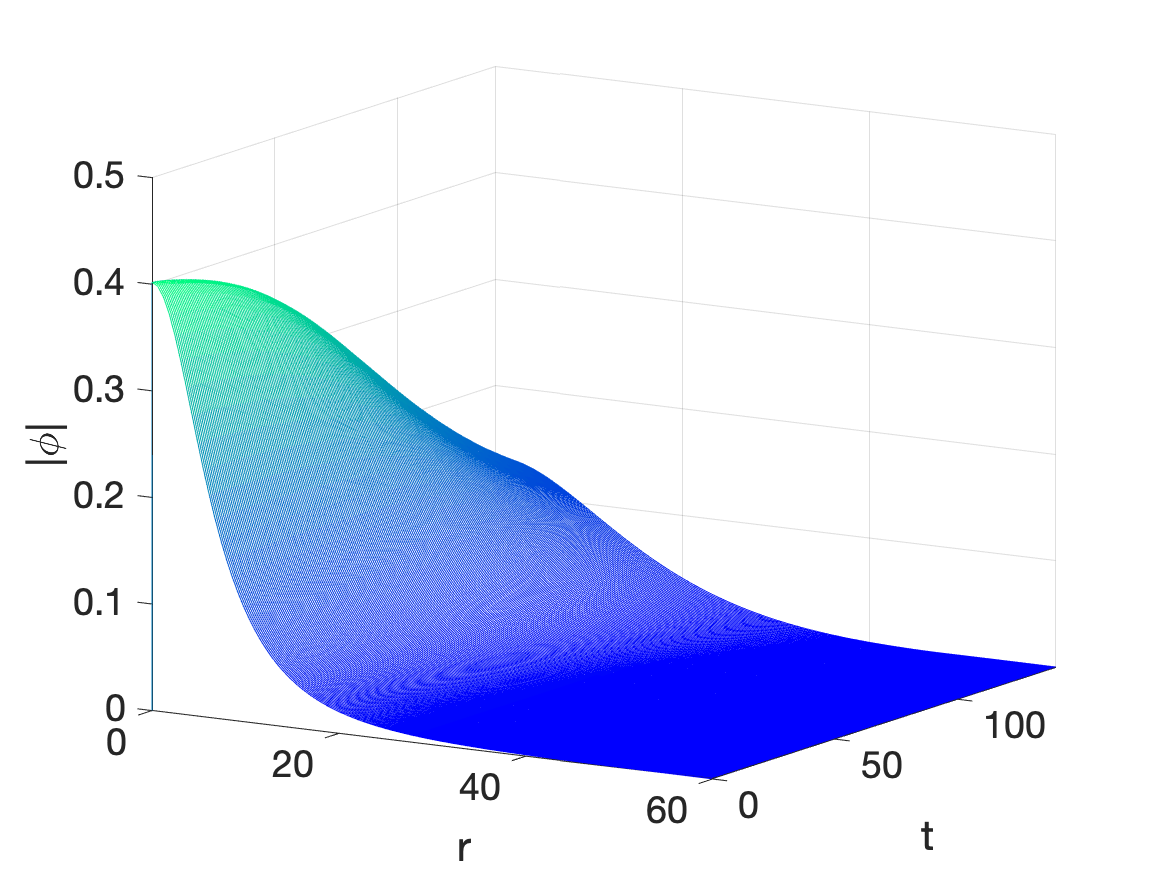}
 \caption{Solution to equation (\ref{SSNLt}) for initial data 
 (\ref{solpert}) with $\lambda=0.99$.}
 \label{NLS3Dsolom001_099}
\end{figure}

If we perturb solitary waves on the unstable branch, say for 
$\omega=0.01$, the behavior changes considerably. Since  $\omega$
is quite small in comparison to the  case $d=1$ studied before, we need
to compute for much longer times. In 
Fig.~\ref{NLS3Dsolom001_099} we show the solution for initial data of
the form $\phi(x,0)=0.99 \varphi_\omega(r)$. The solution appears to be
simply dispersed as in the $d=1$ case.

\begin{figure}[htb!]
  \includegraphics[width=0.49\textwidth]{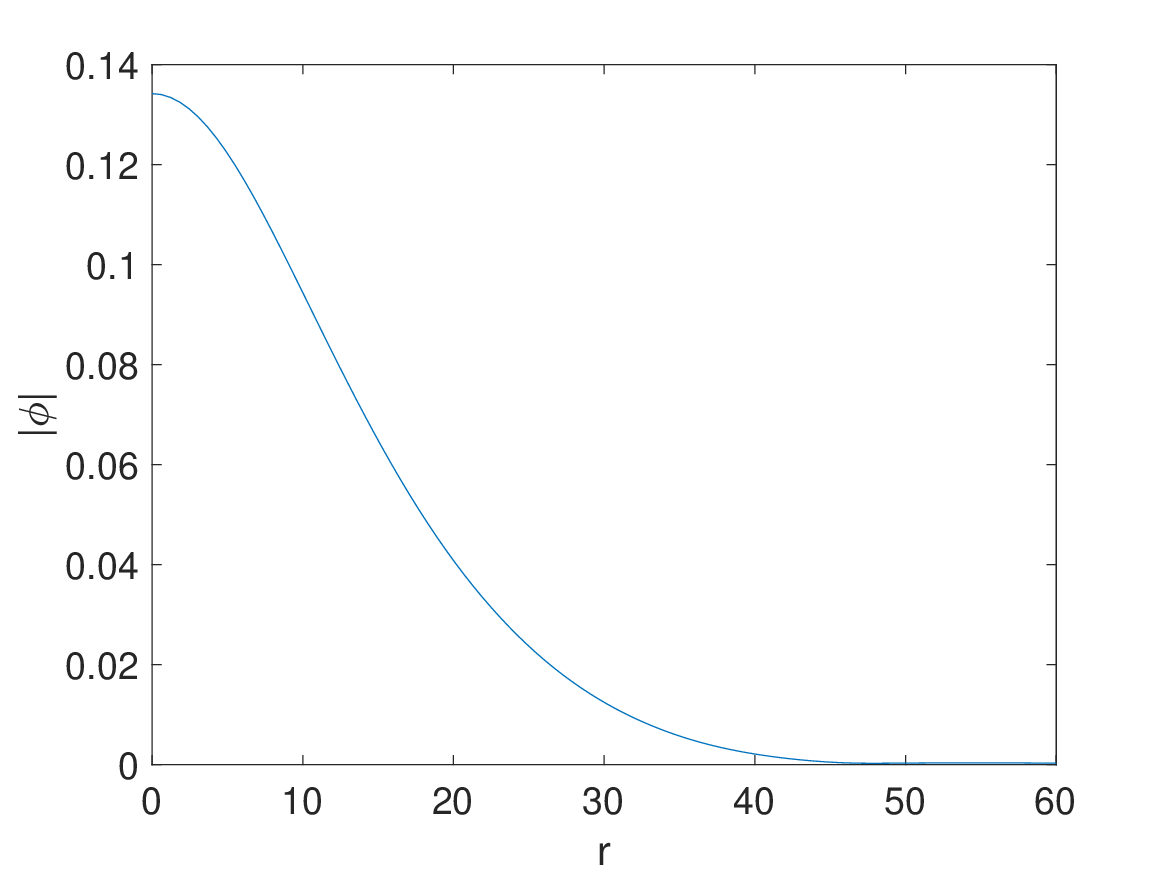}
  \includegraphics[width=0.49\textwidth]{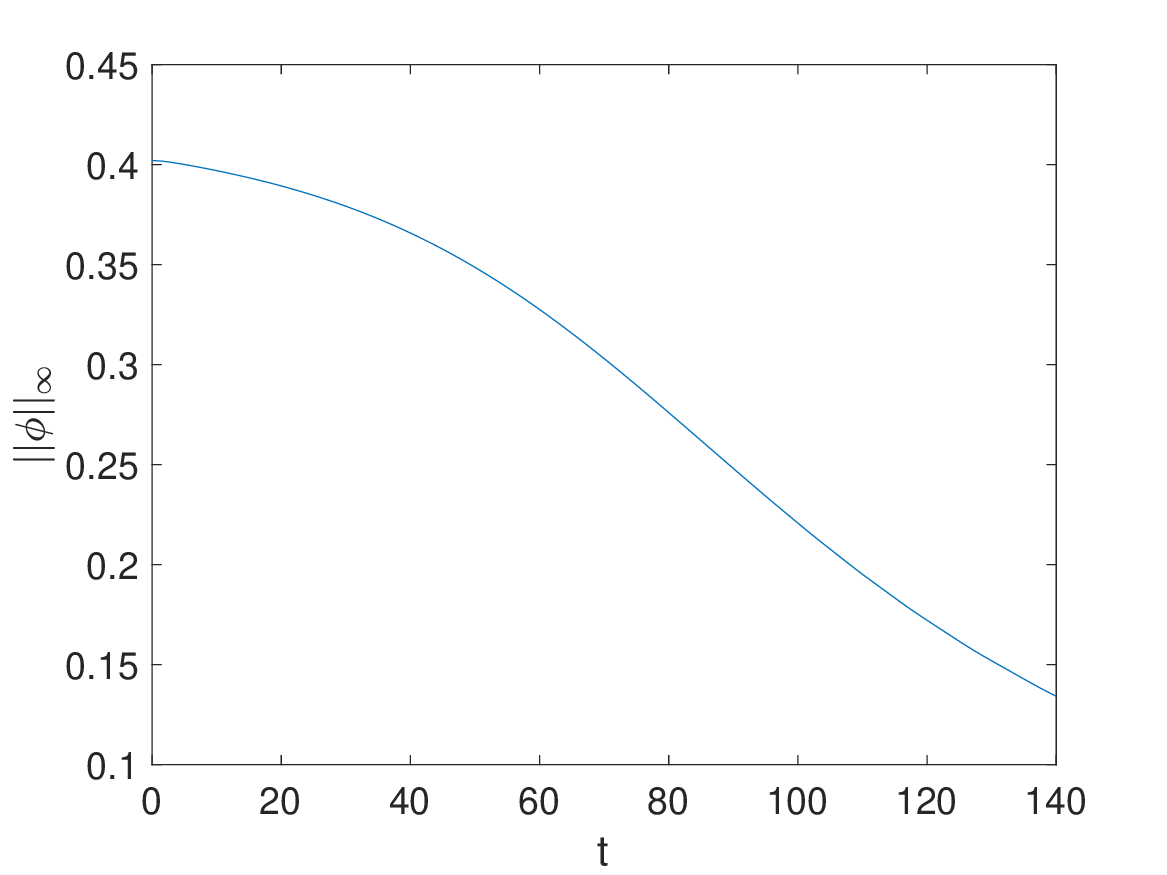}
 \caption{Solution to equation (\ref{SSNLt}) for initial data 
 (\ref{solpert}) with $\lambda=0.99$: on the left the 
  solution for $t=140$, on the 
 right the $L^{\infty}$ norm of the solution in dependence of time.}
 \label{NLS3Dsolom01_099max}
\end{figure}

This interpretation is confirmed by the $L^{\infty}$ norm of the 
solution in Fig.~\ref{NLS3Dsolom01_099max} on the right. The solution
at the final recorded time is shown on the left of the same figure. 

However, if the same solitary wave is perturbed such that the initial 
data have a larger mass than the solitary wave, the final state 
appears to be a solitary wave on the stable branch plus radiation.  
\begin{figure}[htb!]
  \includegraphics[width=0.7\textwidth]{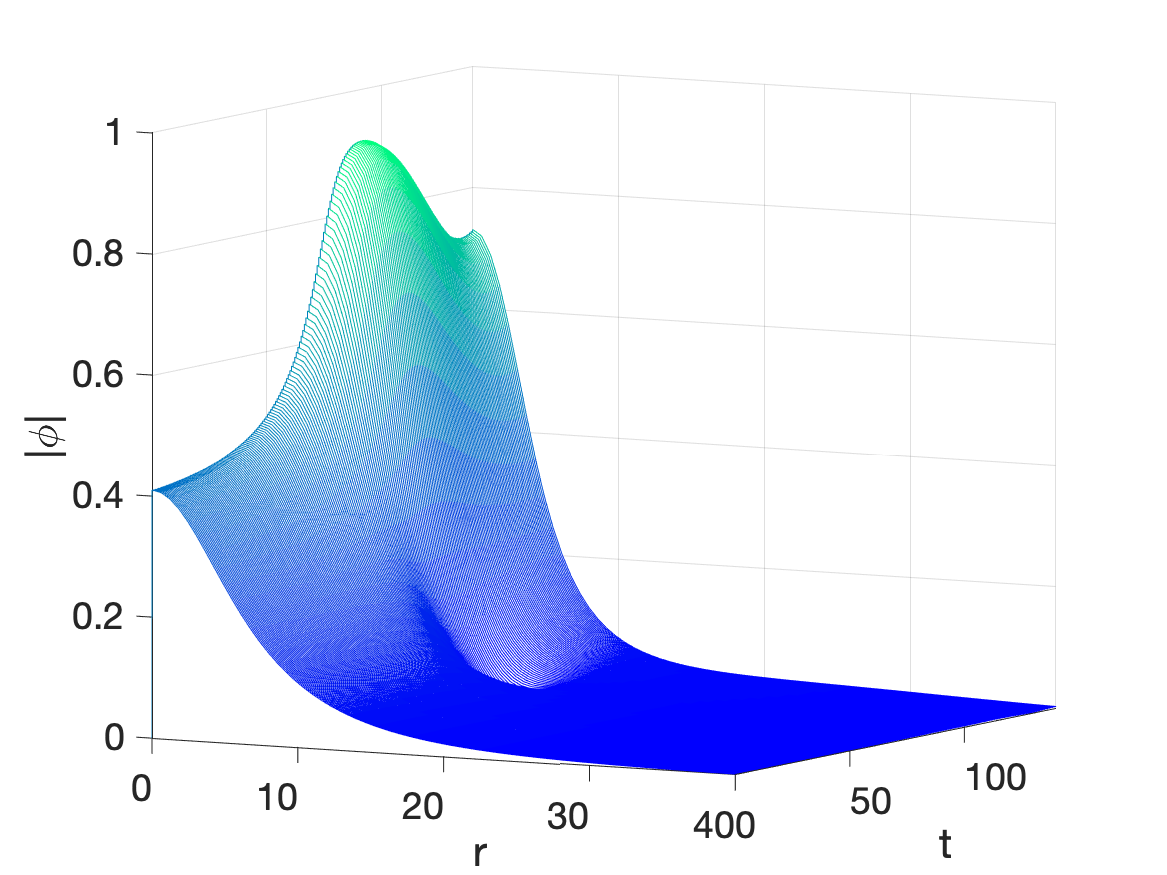}
 \caption{Solution to equation (\ref{SSNLt}) for initial data 
 (\ref{solpert}) with $\lambda=1.01$.}
 \label{NLS3Dsolom01_101water}
\end{figure}

\begin{figure}[htb!]
  \includegraphics[width=0.49\textwidth]{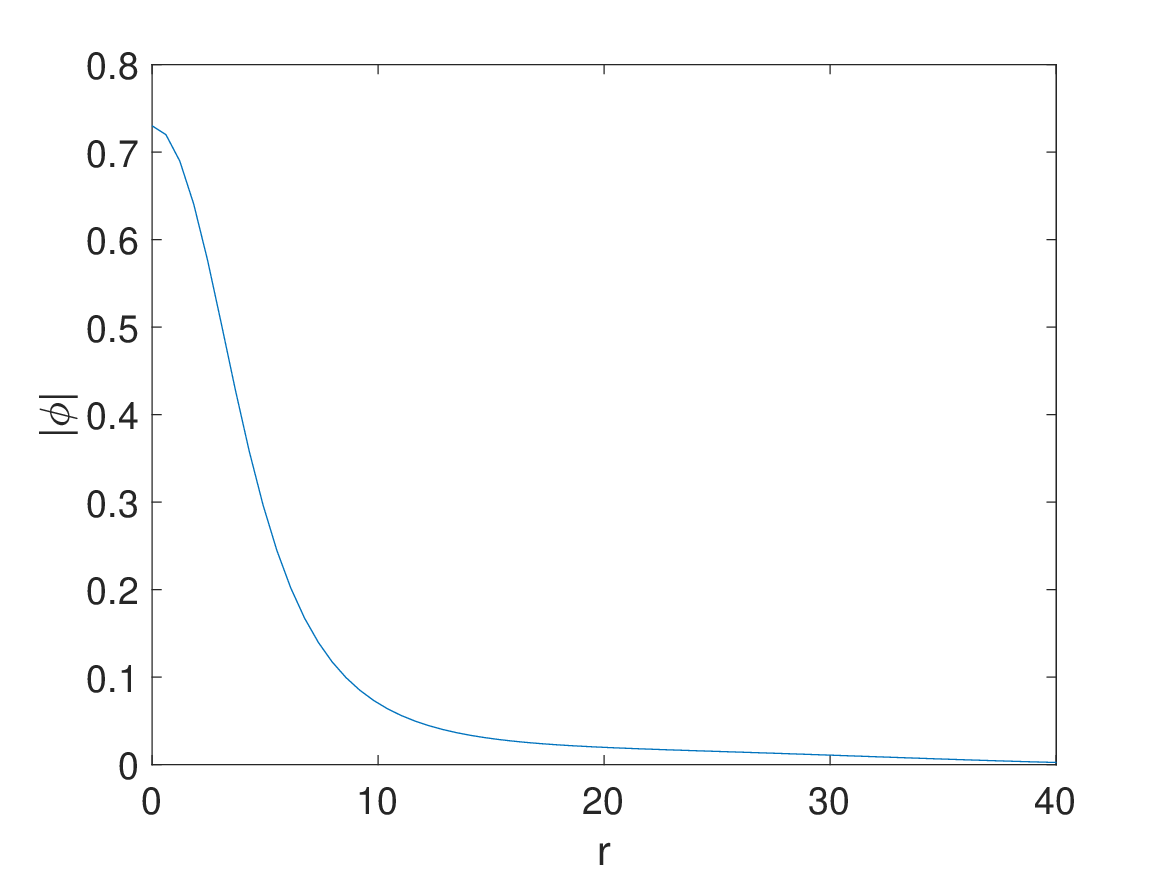}
  \includegraphics[width=0.49\textwidth]{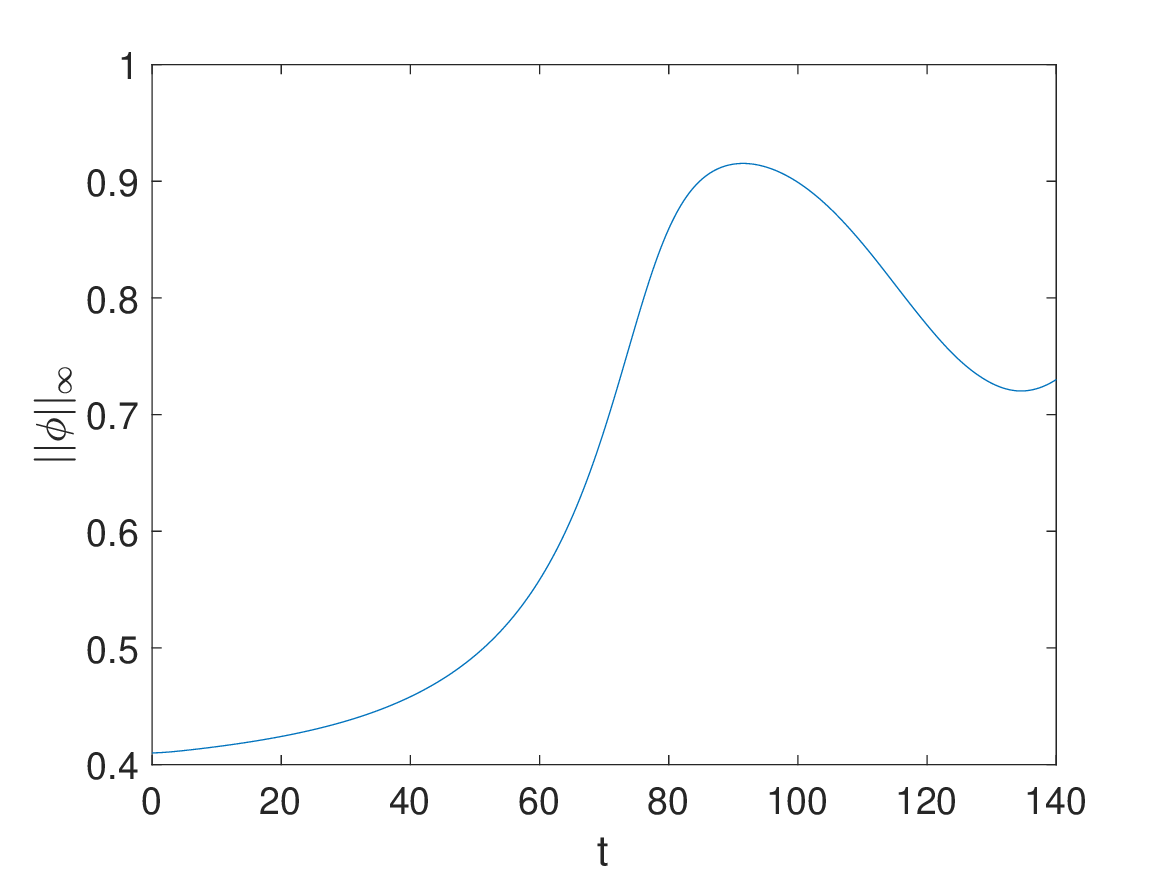}
 \caption{Solution to equation (\ref{SSNLt}) for initial data 
 (\ref{solpert}) with $\lambda=1.01$: on the left the 
  solution for $t=140$, on the 
 right the $L^{\infty}$ norm of the solution in dependence of time.}
 \label{NLS3Dsolom01_101max}
\end{figure}

The solution at 
the final recorded time on the left of Fig.~\ref{NLS3Dsolom01_101max}
indicates that it approaches  a stable ground state. There is still strong
radiation near the central peak which means that it will take much 
longer times to reach the pure solitary wave. This is confirmed by 
the $L^{\infty}$ norm of the solution on the right of the same figure 
which appears to oscillate around an asymptotic value. 

\subsection{Localized initial data}
In this section we study the time evolution of rapidly decreasing 
initial data of the form
\begin{equation}
	\phi(r,0) = c \exp(-r^{2}/s_{1}),\quad 0<c< 1
	\label{gauss}.
\end{equation}

\begin{figure}[htb!]
  \includegraphics[width=0.7\textwidth]{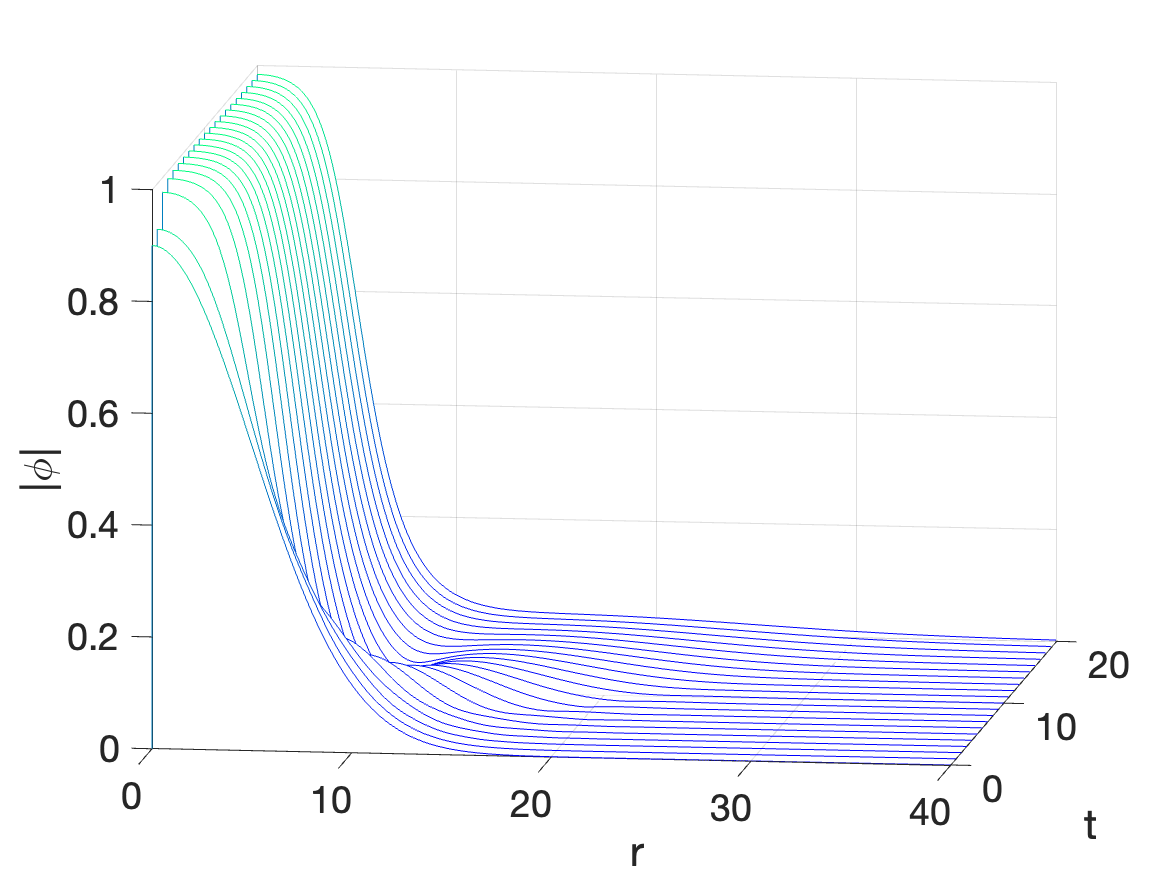}
 \caption{Solution to equation (\ref{SSNLt}) for initial data 
 (\ref{gauss}) with $c=0.9$ and $s_{m}=50$.}
 \label{NLS3D09gausswater}
\end{figure}

We use the same numerical parameters as in the previous subsection for 
$c=0.9$ and $s_{1}=50$, the initial iterate in the construction 
of the ground state for $\omega=0.1$. The resulting solution is shown
in Fig.~\ref{NLS3D09gausswater}. It can be seen that some
radiation is emitted, but that the solution seems to approach a 
ground state for long times. 

\begin{figure}[htb!]
  \includegraphics[width=0.49\textwidth]{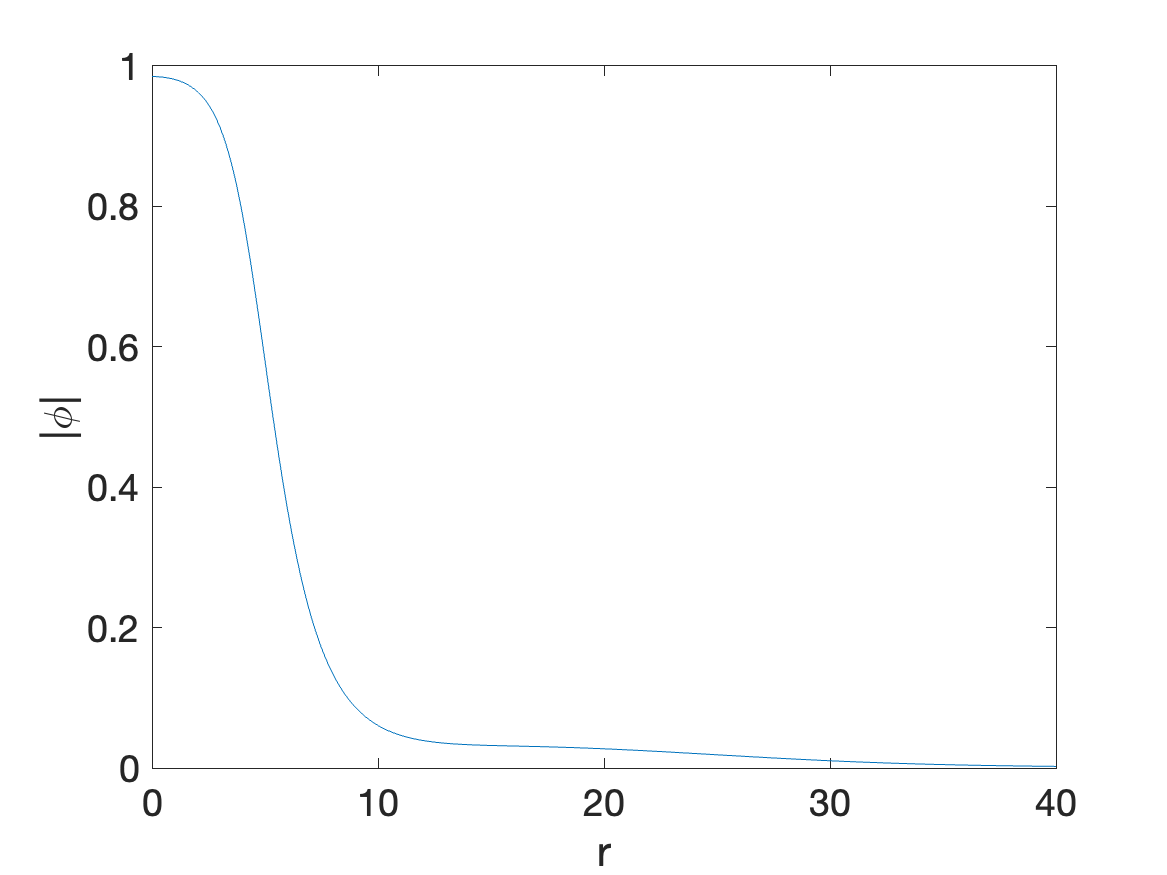}
  \includegraphics[width=0.49\textwidth]{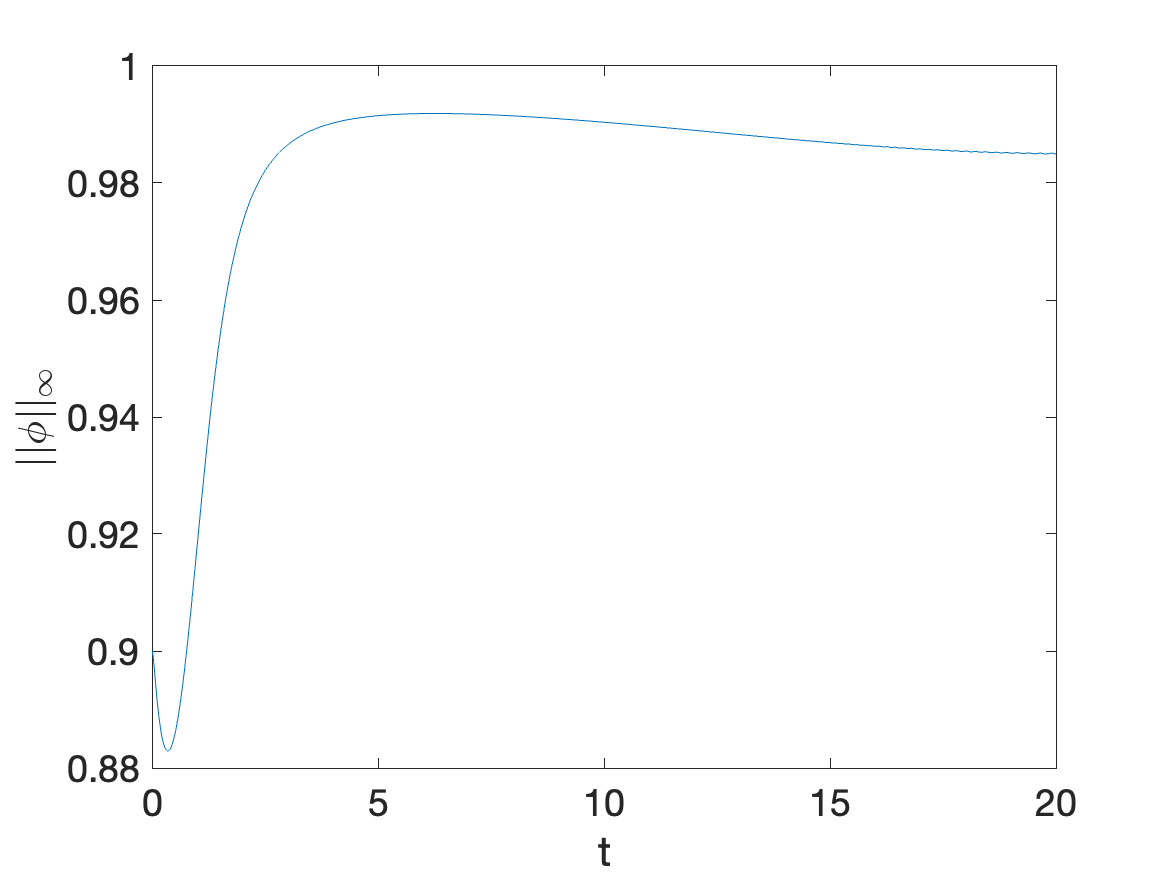}
 \caption{Solution to equation (\ref{SSNLt}) for initial data 
 (\ref{gauss}) with $c=0.9$ and $s_{1}=50$: on the left the solution for 
 $t=20$, on the right the $L^{\infty}$ norm.}
 \label{NLS3D09gaussinf}
\end{figure}

The solution at the final time is shown on the left of 
Fig.~\ref{NLS3D09gaussinf}. It can be recognized that there is a
slower fall-off than in the ground states, which indicates
that there is still some radiation to be emitted before the ground 
state fully emerges. The interpretation of this final state as a 
ground state is in accordance with the $L^{\infty}$ norm of the 
solution shown on the right of the same figure. It appears to slowly reach a 
plateau. 

\begin{figure}[htb!]
  \includegraphics[width=0.7\textwidth]{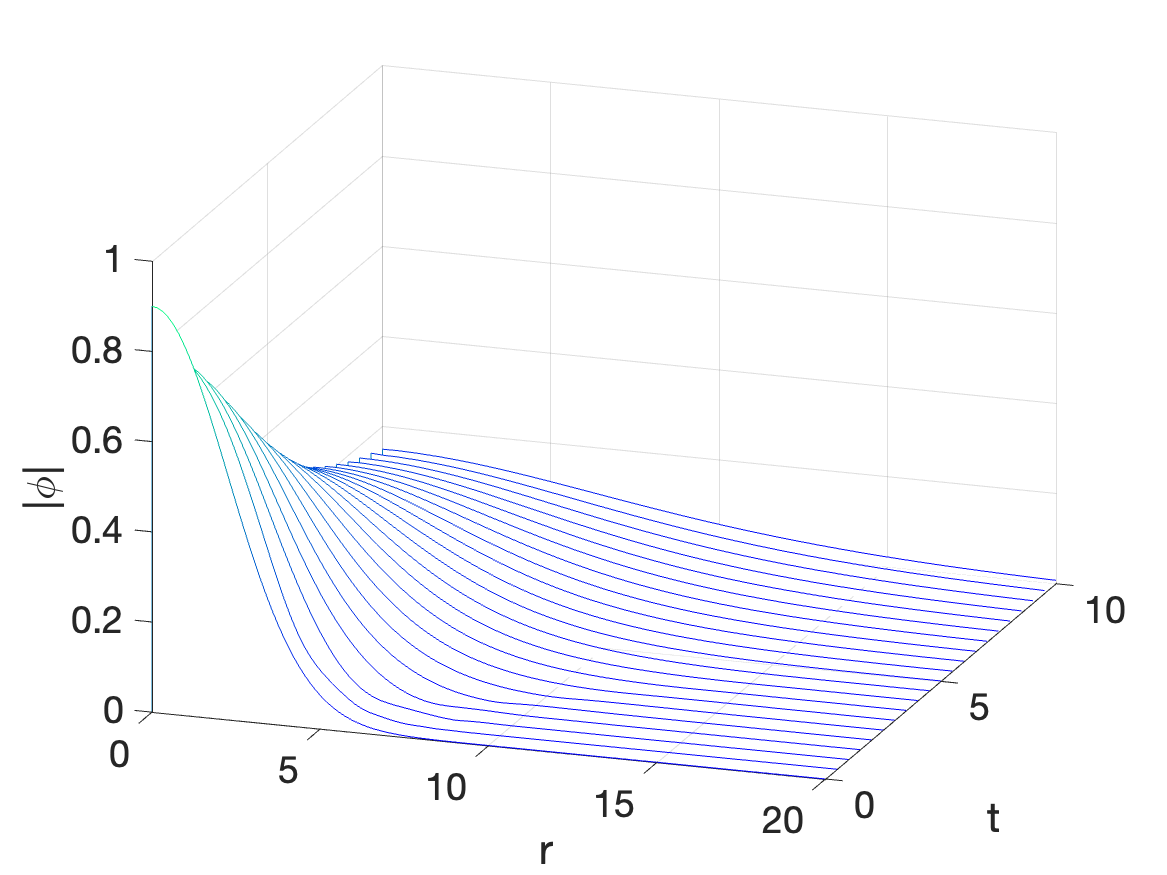}
 \caption{Solution to equation (\ref{SSNLt}) for initial data 
 (\ref{gauss}) with $c=0.9$ and $s_{1}=10$.}
 \label{NLS3D09gauss10}
\end{figure}

The situation is somewhat different if initial data of the form 
(\ref{gauss}) are considered that are more localized, for instance
$c=0.9$ and $s_{1}=10$. We use the same numerical parameters, but on 
shorter time scales, $t\leq 10$. The solution is shown in 
Fig.~\ref{NLS3D09gauss10}. The initial hump is rapidly radiated away,
no stable structure seems to emerge. 

\begin{figure}[htb!]
  \includegraphics[width=0.49\textwidth]{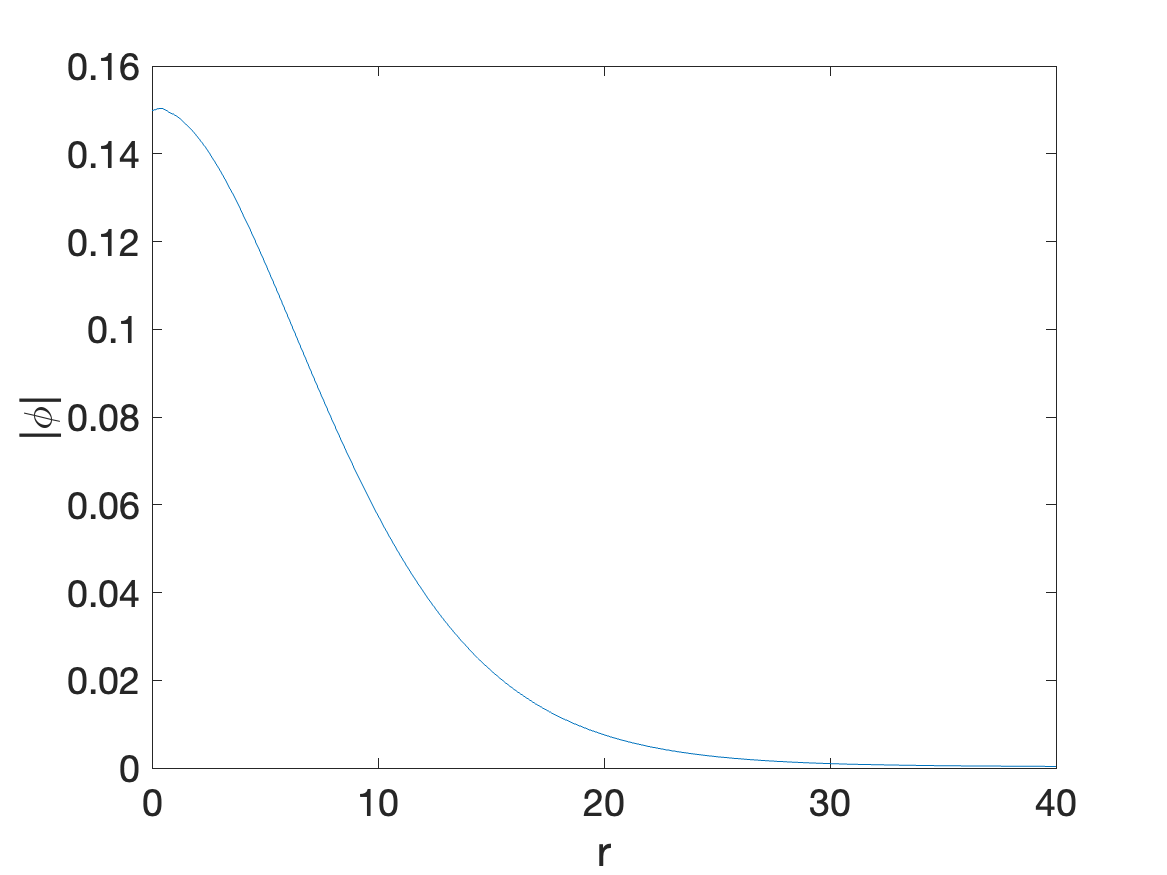}
  \includegraphics[width=0.49\textwidth]{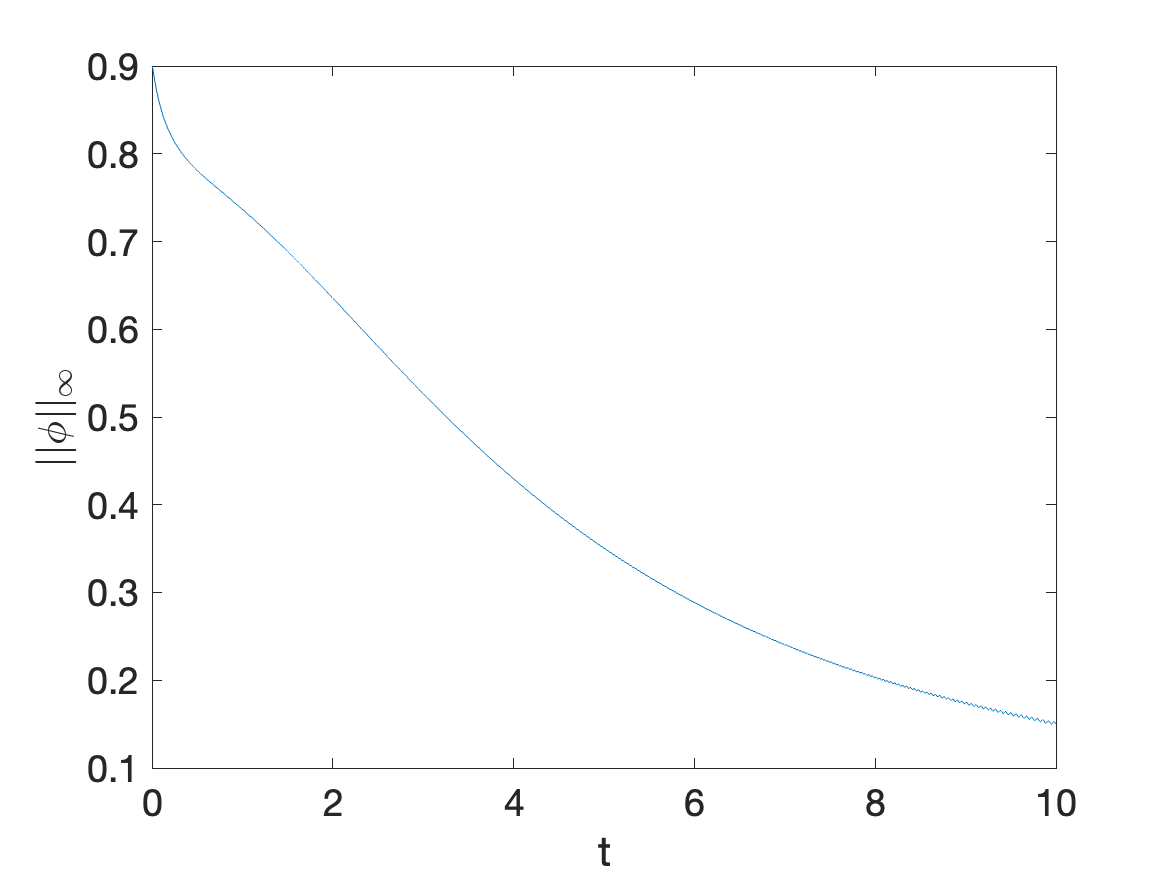}
 \caption{Solution to equation (\ref{SSNLt}) for initial data 
 (\ref{gauss}) with $c=0.9$ and $s_{1}=10$: on the left the solution for 
 $t=10$, on the right the $L^{\infty}$ norm.}
 \label{NLS3D09gauss10inf}
\end{figure}

The solution at the final time is shown in 
Fig.~\ref{NLS3D09gauss10inf} on the left. The central peak appears to
quickly decrease in height and broaden, it just seems to be radiated 
away on large time scales. This is in accordance with $L^{\infty}$ 
norm on the right of the same figure which appears to decrease 
monotonically. Whereas no condition on the initial data is known for 
ground states to appear, it seems to be that the data need to have 
sufficient bulk near the maximum. 


\appendix

\section{Proof of Theorem~\ref{thm:asymptotics1}}

For any $\omega\in (0,\omega^*)$ and for any $z\in (0,+\infty)$, we define 
\begin{equation*}
  f(\omega,z)=\left( 1 + z^2 \right)^{-1/\alpha}\left( \left( 1-\frac{\omega}{\omega^*} \right) +z^2 \right)^{-1/2}
\end{equation*} and we recall that 
$M(\omega)=\frac{1}{\alpha} \, \frac{\omega^{1/\alpha - 1/2}}{(\omega^*)^{1/\alpha}} 
\int_0^{+\infty} f(\omega,z) \, dz$. For any $\omega\in (0,\omega^*)$, $f(\omega,\cdot)$ is clearly integrable. Moreover, for $\omega$ small enough, $f(\omega,\cdot)$ can be uniformly bounded by an integrable function. Hence the asymptotics as $\omega\to 0^+$ follow from  Lebesgue's dominated convergence theorem.

In the limit $\omega\to \omega^*$, the function $f$ has an non-integrable singularity at the origin. More precisely, 
\begin{align*}
  \int_0^{1} f(\omega,z) \, dz &= \int_0^{1} \left( 1 + z^2 \right)^{-1/\alpha}\left( \left( 1-\frac{\omega}{\omega^*} \right) +z^2 \right)^{-1/2} \, dz \\
  &= \int_0^{1} \left( \left( 1-\frac{\omega}{\omega^*} \right) +z^2 \right)^{-1/2} \, dz + \int_0^{1} \left(\left( 1 + z^2 \right)^{-1/\alpha}-1\right)\left( \left( 1-\frac{\omega}{\omega^*} \right) +z^2 \right)^{-1/2} \, dz.
\end{align*}
Now, the integrand in the last integral is uniformly bounded by constant. Using again Lebesgue's dominated convergence theorem, we deduce that 
\begin{equation*}
  0\le \lim_{\omega\to\omega^*}\int_0^{1} \left(1-\left( 1 + z^2 \right)^{-1/\alpha}\right)\left( \left( 1-\frac{\omega}{\omega^*} \right) +z^2 \right)^{-1/2} \, dz<+\infty.
\end{equation*}
The first integral can be computed explicitly. We have
\begin{align*}
  \int_0^{1} \left( \left( 1-\frac{\omega}{\omega^*} \right) +z^2 \right)^{-1/2} \, dz &=\arcsinh\left(\left( 1-\frac{\omega}{\omega^*} \right)^{-1/2}\right)\\
  &=\ln\left(\left( 1-\frac{\omega}{\omega^*} \right)^{-1/2}\right)+\ln\left(1+\sqrt{1+\left( 1-\frac{\omega}{\omega^*} \right)^{1/2}}\right)\underset{\omega \to \omega^*}{\sim}-\frac{1}{2}\ln(\omega^*-\omega).
\end{align*}
Finally, for any $\omega\in (0,\omega^*)$, $f(\omega,\cdot)$ can be uniformly bounded by an integrable function on $(1,+\infty)$. As a conclusion, 
\begin{equation*}
  M(\omega) \underset{\omega \to \omega^*}{\sim}-\frac{1}{2\alpha \sqrt{\omega^*}}\ln(\omega^*-\omega).
\end{equation*}
\medskip
Next, the function $\omega \in (0,\omega^*)\mapsto f(\omega,z)$ is $\mathcal C^2$ and 
\begin{align*}
  &\partial_\omega f(\omega,z)=\frac{1}{2\omega^*}\left( 1 + z^2 \right)^{-1/\alpha}\left( \left( 1-\frac{\omega}{\omega^*} \right) +z^2 \right)^{-3/2}=\frac{1}{2\omega^*}f(\omega,z)\left( \left( 1-\frac{\omega}{\omega^*} \right) +z^2 \right)^{-1},\\
  &\partial^2_\omega f(\omega,z)=\frac{3}{4(\omega^*)^2}\left( 1 + z^2 \right)^{-1/\alpha}\left( \left( 1-\frac{\omega}{\omega^*} \right) +z^2 \right)^{-5/2}=\frac{3}{2\omega^*}\partial_\omega f(\omega,z)\left( \left( 1-\frac{\omega}{\omega^*} \right) +z^2 \right)^{-1}.
\end{align*}
If $\alpha=2$, using the same arguments as above, we have
\begin{equation*}
  \partial_\omega M(\omega)=\frac{1}{4(\omega^*)^{3/2}}\int_0^{+\infty}\left( 1 + z^2 \right)^{-1/2}\left( \left( 1-\frac{\omega}{\omega^*} \right) +z^2 \right)^{-3/2}\,dz >0
\end{equation*}
and
\begin{equation*}
  \partial^2_\omega M(\omega)=\frac{3}{8(\omega^*)^{5/2}}\int_0^{+\infty}\left( 1 + z^2 \right)^{-1/2}\left( \left( 1-\frac{\omega}{\omega^*} \right) +z^2 \right)^{-5/2}\,dz >0
\end{equation*}
for all $\omega\in (0,\omega^*)$. Hence $M(\omega)$ is strictly increasing and convex on $(0,\omega^*)$.  
If $\alpha\neq 2$, we obtain, for any $\omega\in (0,+\infty)$,
\begin{align*}
  &\partial_\omega M(\omega)=\left(\frac{1}{\alpha}-\frac12\right)\frac{1}{\omega}M(\omega)+\frac{\omega^{1/\alpha-1/2}}{2\alpha(\omega^*)^{1+1/\alpha}}\int_0^{+\infty}\left( 1 + z^2 \right)^{-1/\alpha}\left( \left( 1-\frac{\omega}{\omega^*} \right) +z^2 \right)^{-3/2}\,dz
\end{align*}
which is clearly strictly positive whenever $\left(\frac{1}{\alpha}-\frac12\right)>0$, \emph{i.e.} $\alpha < 2$. 

If $\alpha>2$, a straightforward computation leads to  
\begin{align*}
  \partial_\omega M(\omega)&=\frac{\omega^{1/\alpha-3/2}\omega^*}{\alpha (\omega^*)^{1/\alpha}}\left[
  \left(\frac{1}{\alpha}-\frac12\right)\int_0^{+\infty} \frac{1}{\omega_*}f(\omega,z)\, dz+\frac{\omega}{\omega^*}\int_0^{+\infty} \partial_{\omega}f(\omega,z)\, dz\right]\\
  &=\frac{\omega^*}{\alpha (\omega^*)^{1/\alpha}}\omega^{1/\alpha-3/2}
  \int_0^{+\infty} \partial_{\omega}f(\omega,z)\left(\frac{2-\alpha}{\alpha}\left( \left( 1-\frac{\omega}{\omega^*} \right) +z^2 \right)+\frac{\omega}{\omega^*}\right)\, dz
\end{align*}
where we use the fact that $\frac{1}{\omega_*}f(\omega,z)=2 \partial_{\omega}f(\omega,z)\left( \left( 1-\frac{\omega}{\omega^*} \right) +z^2 \right)$. As a consequence, 
\begin{align*}
  \partial^2_\omega M(\omega)=&\,\frac{\omega^*}{\alpha (\omega^*)^{1/\alpha}}\omega^{1/\alpha-5/2}\left[\left(\frac{1}{\alpha}-\frac32\right)\int_0^{+\infty} \partial_{\omega}f(\omega,z)\left(\frac{2-\alpha}{\alpha}\left( \left( 1-\frac{\omega}{\omega^*} \right) +z^2 \right)+\frac{\omega}{\omega^*}\right)\, dz\right.\\
  &\left.+ \omega\int_0^{+\infty} \partial^2_{\omega}f(\omega,z)\left(\frac{2-\alpha}{\alpha}\left( \left( 1-\frac{\omega}{\omega^*} \right) +z^2 \right)+\frac{\omega}{\omega^*}\right)\, dz+\frac{2(\alpha-1)}{\alpha}\frac{\omega}{\omega^*}\int_0^{+\infty} \partial_{\omega}f(\omega,z)\, dz\right]\\
  =&\,\frac{(\omega^*)^2}{\alpha (\omega^*)^{1/\alpha}}\omega^{1/\alpha-5/2}\left[\int_0^{+\infty} \frac{1}{\omega^*}\partial_{\omega}f(\omega,z)\left(\frac{3\alpha-2}{2\alpha}\frac{\alpha-2}{\alpha}\left( \left( 1-\frac{\omega}{\omega^*} \right) +z^2 \right)+\frac{\alpha-2}{2\alpha}\frac{\omega}{\omega^*}\right)\, dz+\right.\\
  &\left.+ \frac{\omega}{\omega^*}\int_0^{+\infty} \partial^2_{\omega}f(\omega,z)\left(\frac{2-\alpha}{\alpha}\left( \left( 1-\frac{\omega}{\omega^*} \right) +z^2 \right)+\frac{\omega}{\omega^*}\right)\, dz\right]\\
  =&\,\frac{(\omega^*)^2}{\alpha (\omega^*)^{1/\alpha}}\omega^{1/\alpha-5/2}\left[\int_0^{+\infty} \partial^2_{\omega}f(\omega,z)\left(\frac{3\alpha-2}{3\alpha}\frac{\alpha-2}{\alpha}\left( \left( 1-\frac{\omega}{\omega^*} \right) +z^2 \right)^2\right)\, dz+\right.\\
  &\left.+ \int_0^{+\infty} \partial^2_{\omega}f(\omega,z)\left(-\frac{2(\alpha-2)}{3\alpha}\left( \left( 1-\frac{\omega}{\omega^*} \right) +z^2 \right)\frac{\omega}{\omega^*}+\left(\frac{\omega}{\omega^*}\right)^2\right)\, dz\right]\\
\end{align*}
where we use $\frac{1}{\omega_*}\partial_\omega f(\omega,z)=\frac{2}{3} \partial^2_{\omega}f(\omega,z)\left( \left( 1-\frac{\omega}{\omega^*} \right) +z^2 \right)$. To conclude we prove that the function 
\begin{equation*}
  G_\alpha(\omega,z):= \frac{3\alpha-2}{3\alpha}\frac{\alpha-2}{\alpha}\left( \left( 1-\frac{\omega}{\omega^*} \right) +z^2 \right)^2 -\frac{2(\alpha-2)}{3\alpha}\left( \left( 1-\frac{\omega}{\omega^*} \right) +z^2 \right)\frac{\omega}{\omega^*}+\left(\frac{\omega}{\omega^*}\right)^2>0
\end{equation*}
for all $\alpha>2$, $\omega\in (0,\omega^*)$ and $z\in (0,+\infty)$. Indeed, $G_\alpha(\omega,z)$ can be written as 
\begin{align*}
  G_\alpha(\omega,z)=\left(\frac{\omega}{\omega^*}-\frac{\alpha-2}{3\alpha}\left( \left( 1-\frac{\omega}{\omega^*} \right) +z^2 \right)\right)^2+\frac{\alpha-2}{3\alpha}\left(\frac{8\alpha-4}{3\alpha}\right)\left( \left( 1-\frac{\omega}{\omega^*} \right) +z^2 \right)^2
\end{align*} 
which is strictly positive whenever $\alpha>2$. This concludes the proof of Theorem~\eqref{thm:asymptotics1}.


\begin{thebibliography}{10}

\bibitem{BahhiPhD-25}
{\sc M.~Bahhi}, {\em Étude mathématique d'équations effectives en dynamique
  quantique relativiste.}, phd theis, Université Bourgogne Europe, Dijon,
  France, 2025.

\bibitem{birem}
{\sc M.~Birem and C.~Klein}, {\em Multidomain spectral method for
  schr{\"o}dinger equations}, Adv. Comput. Math., 42 (2016), pp.~395--423.

\bibitem{CKS}
{\sc R.~Carles, C.~Klein, and C.~Sparber}, {\em On ground state (in-)stability
  in multi-dimensional cubic-quintic schrödinger equations}, ESAIM: M2AN, 57
  (2023), pp.~423--443.

\bibitem{CC}
{\sc C.~Clenshaw and A.~Curtis}, {\em A method for numerical integration on an
  automatic computer.}, Numer. Math., 2 (1960), pp.~197--205.

\bibitem{BieGenRot-15}
{\sc S.~{De Bi\`evre}, F.~{Genoud}, and S.~{Rota Nodari}}, {\em Orbital
  stability: analysis meets geometry}, in Nonlinear optical and atomic systems.
  At the interface of physics and mathematics, Cham: Springer; Lille: Centre
  europ\'een pour les math\'ematiques, la physique et leurs interactions
  (CEMPI), 2015, pp.~147--273.

\bibitem{BieRot-19}
{\sc S.~De~Bi{\`e}vre and S.~Rota~Nodari}, {\em Orbital stability via the
  energy-momentum method: the case of higher dimensional symmetry groups},
  Arch. Rational Mech. Anal., 231 (2019), pp.~233--284.

\bibitem{EstRot-12}
{\sc M.~J. Esteban and S.~{Rota Nodari}}, {\em Symmetric ground states for a
  stationary relativistic mean-field model for nucleons in the nonrelativistic
  limit}, Rev. Math. Phys., 24 (2012), pp.~1250025--1250055.

\bibitem{EstRot-13}
\leavevmode\vrule height 2pt depth -1.6pt width 23pt, {\em Ground states for a
  stationary mean-field model for a nucleon}, Ann. Henri Poincar{\'e}, 14
  (2013), pp.~1287--1303.

\bibitem{GenRot-24}
{\sc F.~Genoud and S.~Rota~Nodari}, {\em Standing wave solutions of a
  quasilinear schr\"odinger equation in the small frequency limit}, 2024.

\bibitem{GriShaStr-87}
{\sc M.~Grillakis, J.~Shatah, and W.~Strauss}, {\em Stability theory of
  solitary waves in the presence of symmetry. {I}}, J. Funct. Anal., 74 (1987),
  pp.~160--197.

\bibitem{GriShaStr-90}
\leavevmode\vrule height 2pt depth -1.6pt width 23pt, {\em Stability theory of
  solitary waves in the presence of symmetry. {II}}, J. Funct. Anal., 94
  (1990), pp.~308--348.

\bibitem{kenig2004cauchy}
{\sc C.~E. Kenig, G.~Ponce, and L.~Vega}, {\em The cauchy problem for
  quasi-linear {S}chr{\"o}dinger equations}, Inventiones mathematicae, 158
  (2004).

\bibitem{etna}
{\sc C.~Klein}, {\em Fourth order time-stepping for low dispersion korteweg-de
  vries and nonlinear schrödinger equations.}, Electron. Trans. Numer. Anal.,
  29 (2007), pp.~116--135.

\bibitem{KRN}
{\sc C.~Klein and S.~Rota~Nodari}, {\em On a nonlinear schr{\"o}dinger equation
  for nucleons in one space dimension}, ESAIM: M2AN,  (2021), pp.~409--427.

\bibitem{TreRot-13}
{\sc L.~{Le Treust} and S.~{Rota Nodari}}, {\em Symmetric excited states for a
  mean-field model for a nucleon}, Journal of Differential Equations, 255
  (2013), pp.~3536--3563.

\bibitem{LewRot-15}
{\sc M.~Lewin and S.~{Rota Nodari}}, {\em Uniqueness and non-degeneracy for a
  nuclear nonlinear {S}chr{\"o}dinger equation}, NoDEA Nonlinear Differential
  Equations Appl., 22 (2015), pp.~673--698.

\bibitem{LewRot-20}
{\sc M.~Lewin and S.~Rota~Nodari}, {\em The double-power nonlinear
  {S}chr\"odinger equation and its generalizations: uniqueness, non-degeneracy
  and applications}, Calc. Var. Partial Differential Equations, 59 (2020),
  pp.~Paper No. 197, 49.

\bibitem{MorMur-14}
{\sc V.~Moroz and C.~B. Muratov}, {\em Asymptotic properties of ground states
  of scalar field equations with a vanishing parameter}, J. Eur. Math. Soc.
  (JEMS), 16 (2014).

\bibitem{Rota-HDR}
{\sc S.~{Rota Nodari}}, {\em Contributions to the mathematical study of models
  from classical and quantumphysics}, habilitation à diriger des recherches,
  Universit{\'e} de Bourgogne, February 2021.

\bibitem{trefethen}
{\sc L.~N. Trefethen}, {\em Spectral Methods in MATLAB}, Society for Industrial
  and Applied Mathematics, 2000.

\bibitem{WR}
{\sc J.~A. Weideman and S.~C. Reddy}, {\em A matlab differentiation matrix
  suite}, ACM Trans. Math. Softw., 26 (2000), p.~465–519.

\bibitem{Weinstein-85}
{\sc M.~I. Weinstein}, {\em Modulational stability of ground states of
  nonlinear {S}chr{\"o}dinger equations}, SIAM J. Math. Anal., 16 (1985),
  pp.~472--491.

\end{thebibliography}
\end{document}